\documentclass[]{amsart}
\usepackage[
  top=0.95in,
  bottom=0.85in,
  left=0.90in,
  right=0.90in,
  bindingoffset=0.25in,
  heightrounded,
]{geometry}
\usepackage{babel}
\usepackage[utf8]{inputenc}
\usepackage{amsmath}
\usepackage{amsthm}
\usepackage{amsfonts}
\usepackage{amssymb}
\usepackage{enumerate}
\usepackage{stackrel}
\usepackage{mathrsfs}
\usepackage{graphicx}
\usepackage{tikz-cd}
\usepackage[hidelinks]{hyperref}
\usepackage{fancyhdr}
\usepackage{cleveref}
\usetikzlibrary{babel}
\author{Iván Chércoles Cuesta}
\title{The space of Conradian left-preorders}

\hypersetup{
     colorlinks=true,
     linkcolor=blue,
     citecolor=red,
     }

\providecommand{\df}[3]{#1:#2 \longrightarrow #3}
\providecommand{\natp}{\mathbb{N}_{>0}}
\providecommand{\nat}{\mathbb{N}}
\providecommand{\ent}{\mathbb{Z}}
\providecommand{\real}{\mathbb{R}}
\providecommand{\co}{\circ}
\providecommand{\sii}{\quad \Longleftrightarrow \quad}

\providecommand{\cont}{\subseteq}

\providecommand{\inv}{^{-1}}

\providecommand{\tq}{\;/\;}

\providecommand{\abs}[1]{|#1|}
\providecommand{\gen}[1]{\langle #1\rangle}
\providecommand{\semgen}[1]{\langle #1\rangle^{+}}

\providecommand{\rela}{\mathcal{R}}

\providecommand{\nor}{\trianglelefteq}
\providecommand{\card}[1]{\mathrm{Card}\left(#1\right)}

\providecommand{\lleq}{\preceq}
\providecommand{\lle}{\prec}
\providecommand{\lleqb}{\lleq_{\bullet}}

\providecommand{\settminus}{-}

\providecommand{\cord}{\mathcal{CO}}

\providecommand{\gr}{G}
\providecommand{\sgr}{H}
\providecommand{\sgrb}{K}
\providecommand{\sgrc}{L}
\providecommand{\grb}{H}
\providecommand{\grc}{K}
\providecommand{\ga}{\alpha}
\providecommand{\gb}{\beta}
\providecommand{\gc}{\gamma}
\providecommand{\gd}{\delta}
\providecommand{\gy}{\theta}
\providecommand{\gx}{\xi}
\providecommand{\gz}{\zeta}
\providecommand{\ce}{C}
\providecommand{\cce}{{\ce_0}}
\providecommand{\sce}{S_0}

\providecommand{\uu}{u}
\providecommand{\vv}{v}
\providecommand{\ww}{w}
\providecommand{\mmen}{^{\pm 1}}
\providecommand{\gee}{\phi}
\providecommand{\urep}{\uu_0}
\providecommand{\vrep}{\vv_0}
\providecommand{\wrep}{\ww_0}
\providecommand{\espacio}{X}
\providecommand{\men}{<}
\providecommand{\hacheaux}{\mathcal{H}}

\providecommand{\defrelativeorder}[2]{left-preorder on $#1$ relative to $#2$}

\newtheorem{ej}{Example}[section]
\newtheorem{teo}[ej]{Theorem}
\newtheorem{prop}[ej]{Proposition}
\newtheorem{lema}[ej]{Lemma}
\newtheorem{coro}[ej]{Corollary}
\theoremstyle{definition}
\newtheorem{defi}[ej]{Definition}
\newtheorem{obs}[ej]{Remark}
\newtheorem*{nota}{Notation}

\usepackage{blindtext}

\begin{document}

\address{Instituto de Ciencias Matemáticas (CSIC-UAM-UC3M-UCM), Consejo Superior de Investigaciones Científicas, C/ Nicolás Cabrera, 13–15, Campus de Cantoblanco UAM, 28049 Madrid, Spain.}
\email{ivan.chercoles@gmail.com}
\date{\today}

\let\thefootnote\relax
\footnotetext{The author is supported by the grant CEX2019-000904-S-20-1 funded by: MCIN/AEI/ 10.13039/501100011033. The author is also partially supported by the grant PID2021-126254NB-I00 of the Ministry of Science and Innovation of Spain.}

\maketitle

\begin{abstract}
We define Conradian left-preorders and the space of Conradian left-preorders. We show that this space is either finite or uncountable. We describe conditions that are equivalent to say that the space of Conradian left-preorders is finite. We provide some characterizations of Conradian left-preorders.
\end{abstract}


\section{Introduction}
The main objective of this paper is defining and studying the space of Conradian left-preorders. Conradian left-preorders arise as a generalization of the notion of Conradian left-order, which was defined in \cite{conrad} and has been studied so many times (see for intance \cite{god}, \cite{dynamics-navas}, \cite{caracterizacionconrad-navasrivasclay} and \cite{Conradianfinito-rivas}). As a relevant property, we recall that a group admits a Conradian left-order if and only if it is locally indicable (see Corollary 9.21 from \cite{orderedgroupstopology-clay}).

\begin{defi}
Let $\gr$ be a group. A \textit{left-preorder} $\lleq$ on $\gr$ is a total transitive reflexive binary relation that is left-invariant, meaning that $$\ga_1\lleq \ga_2\mbox{ implies }\ga \ga_1\lleq \ga \ga_2\mbox{  for all }\ga,\ga_1,\ga_2\in \gr.$$ For convenience, assume that the trivial left-preorder is not a left-preorder, where the \textit{trivial left-preorder} is defined by $\ga_1\lleq \ga_2$ for all $\ga_1,\ga_2\in \gr$. We say that  $\ga_1\lle \ga_2$ when  $\ga_1\lleq \ga_2$ but $\ga_2\not\lleq \ga_1$.
\end{defi}

\begin{defi}
A left-preorder $\lleq$ is \textit{Conradian} if for all $\ga,\gb\in \gr$ such that $1\lle \ga$ and $1\lle \gb$ there is $n\in\nat$ satisfying $1\lle \ga\inv\gb\ga^n$. We define $\cord(\gr)$ to be the \textit{set of Conradian left-preorders} on $\gr$ and $\cord_\ce(\gr)$ to be \textit{the set of Conradian left-preorders} on $\gr$ \textit{relative to} $\ce$ proper subgroup of $\gr$.
\end{defi}

We have a first result, which is a corollary of \Cref{teorema lista de equivalencias absolutas con unicidad introduccion}.

\begin{coro}[\Cref{numero de absolutos es finito o incontable}]
Let $\gr$ be a group. Then $\cord(\gr)$ is either finite or uncountable.
\end{coro}

Any left-preorder $\lleq$ on $\gr$ give us a unique decomposition $\gr=P\sqcup P\inv \sqcup \ce$ where $P$ is the positive cone of $\lleq$ and $\ce$ is a proper subgroup of $\gr$, and we say that $\lleq$ is a left-preorder on $\gr$ relative to $\ce$. We can associate each $\lleq$ left-preorder on $\gr$ with the map sending elements in $P$ to $1$, elements in $P\inv$ to $-1$ and elements in $\ce$ to $0$. Hence, we can see $\cord(\gr)$ as a subset of the set $\{\pm 1,0\}^\gr$. We consider on $\{\pm 1,0\}$ the discrete topology and we consider on $\{\pm 1,0\}^\gr$ the product topology. We define the \textit{space of Conradian left-preorders} on $\gr$ as on the set $\cord(\gr)$ endowed with the subspace topology of $\{\pm 1,0\}^\gr$ with the product topology. We can also consider the subspace topology on $\cord_\ce(\gr)$ as it is a subset of $\cord(\gr)$.\\

The main theorem gives us equivalent conditions for $\cord(\gr)$ being finite. Before, we need a definition.

\begin{defi}
Let $\gr$ be a group. A \textit{rational series} starting on $\gr$ is a finite series of the form
\begin{equation}\label{eq racional introduccion}
\gr_0 \nor \gr_1\nor \cdots \nor \gr_n=\gr
\end{equation}
satisfying that each $\gr_i$ is normal on $\gr_{i+1}$ and the group $\gr_{i+1}/\gr_i$ is torsion-free abelian of abelian rank $1$ for all $0\leq i<n$. Given a rational series starting on $\gr$ as in \Cref{eq racional introduccion}, we say that it is of maximal length if for any other $$\grb_0 \nor \grb_1\nor \cdots \nor \grb_k=\gr$$ rational series starting on $\gr$, we have $k\leq n$.
\end{defi}

\begin{teo}[\Cref{teorema lista de equivalencias absolutas con unicidad}]\label{teorema lista de equivalencias absolutas con unicidad introduccion}
Let $\gr$ be a group. Assume that $\cord(\gr)$ is non-empty. Let $\cce=\bigcap_{\sgr\in \hacheaux}\sgr $ for $\hacheaux=\{\sgr \mbox{ subgroup with } \cord_\sgr(\gr)\neq \emptyset \}$. Then, the following conditions are equivalent:
\begin{itemize}
\item[a)] All the elements on $\cord(\gr)$ are isolated.
\item[b)] There exists a finite rational series starting on $\gr$ as in \Cref{eq racional introduccion} of maximal length, and this series has no abelian jumps and satisfies that $\cord(\gr_0)=\emptyset$.
\item[c)] There exists a unique finite rational series starting on $\gr$ as in \Cref{eq racional introduccion} of maximal length, and this series has no abelian jumps and satisfies that $\cord(\gr_0)=\emptyset$.
\item[d)] There exists a finite rational series finite rational series starting on $\gr$ as in \Cref{eq racional introduccion} for $\gr_0=\cce$, and this series has no abelian jumps.
\item[e)] There exists a unique finite rational series starting on $\gr$ as in \Cref{eq racional introduccion} for $\gr_0=\cce$, and this series has no abelian jumps.
\item[f)] $\cord_\cce(\gr)$ is finite.
\item[g)] $\cord(\gr)$ is finite.
\end{itemize}
Moreover, under these equivalent conditions, we have that $\card{\cord(\gr)}=2^{n+1}-2$ and $\card{\cord_\ce(\gr)}=2^n$ for $n$ the length of any series as in b), c), d) and e).
\end{teo}

Also, as a corollary, when we are on the case $\cce=\{1\}$ on the last theorem, we have a description for groups admitting a Conradian left-order. A partial version of this corollary was proved by Rivas on \cite{Conradianfinito-rivas} (theorems A and B on that paper). This is because the space of Conradian left-orders on $\gr$ is equal to $\cord_{\{1\}}(\gr)$ by following the definitions.\\

In order to prove \Cref{teorema lista de equivalencias absolutas con unicidad introduccion}, we prove a theorem that itself also has interest as we provide equivalent conditions for $\cord_\ce(\gr)$ being finite. 

\begin{teo}[\Cref{teorema lista de equivalencias relativas con unicidad}]\label{teorema lista de equivalencias relativas con unicidad introduccion}
Let $\gr$ be a group and $\ce$ be a proper subgroup. Then, the following conditions are equivalent:
\begin{itemize}
\item[a)] All the elements on $\cord_\ce(\gr)$ are isolated and $\cord_\ce(\gr)$ is non-empty.
\item[b)] There exists an element on $\cord_\ce(\gr)$ that is isolated.
\item[c)] There exists a finite rational series on $\gr$ starting on $\gr$ as in \Cref{eq racional introduccion} for $\gr_0=\ce$ and without abelian jumps.
\item[d)] There exists a unique finite rational series on $\gr$ starting on $\gr$ as in \Cref{eq racional introduccion} for $\gr_0=\ce$, and this series does not have abelian jumps.
\item[e)] $\cord_\ce(\gr)$ is finite and non-empty.
\end{itemize}
Moreover, under these equivalent conditions, we have that $\card{\cord_\ce(\gr)}=2^n$ for $n$ the length of any series as in c) and d).
\end{teo}

The space of left-preorders, which is an space containing the space of Conradian left-preorders, has been studied in some cases. For example, has been studied by Decaup and Rond in \cite{valuated-decauprond} for the abelian case. The space of left-preorders is also used as an auxiliary tool in \cite{morris-antolinrivas} with the name of \textit{space of relative orders}. The space of left-preorders can be seen as a generalization of the space of left-orders and, since Conradian left-orders and the space of Conradian left-orders are so relevant, it is reasonable to consider a left-preorders version of them. This leads to our definition of Conradian left-preorder and the space of Conradian left-preorders.\\

In this paper we also check some properties of Conradian left-preorders that are analogous to some basic properties of Conradian left-orders. Conrad characterized Conradian left-orders as left-orders satisfying that "locally" the group is isomorphic to a subgroup of the additive group of the real numbers, and such isomorphism is order preserving. In \cite{Conradianoparados-jimenez}, it is shown that, in the definition of Conradian left-orders, one can take $n=2$. We provide more general versions of these facts for the context of Conradian left-preorders. In \cite{caracterizacionconrad-navasrivasclay} it is provided a dynamical version of Conradian left-orders, as they are identified as actions on the group without crossings. We prove that there is a dynamical version of Conradian left-preorders by identifying them with actions  without crossings on the set of cosets modulo the corresponding subgroup. These kind of actions play a fundamental role on \cite{morris-antolinrivas}. We state it without the definition of crossing, which can be founded on the corresponding section. We now condense in a theorem all the mentioned characterizations of Conradian left-preorders.

\begin{teo}[\Cref{Conradiano como cocientes de saltos} \& \Cref{Conradiano como dinamica}]
Let $\lleq$ be a left-preorder on $\gr$ relative to $\ce$. The following conditions are equivalent:
\begin{itemize}
\item[a)] For every $(\sgrb,\sgrc)$ $\lleq$-convex jump (that is, a pair of distinct convex subgroups with no other convex subgroup between them), we have that $\sgrb$ is normal in $\sgrc$ and the group $\sgrc/\sgrb$ is isomorphic to a subgroup of the additive group of the real numbers, and such isomorphism is order preserving.
\item[b)] For every pair of positive elements $\ga,\gb\in \gr$ we have $\ga\lle\gb\ga^2$.
\item[c)] The action of $\gr$ on $\gr/\ce$ by left-multiplication admits no crossing.
\item[d)] $\lleq$ is Conradian (for every pair of positive elements $\ga,\gb\in \gr$, one has that $\ga\lle\gb \ga^n$ for some $n\in\nat$).
\end{itemize}
\end{teo}

We now examine the distribution of this paper. Section 1 is devoted to the study of some basic properties of left-preorders and its positive cones. In Section 2 we define Conradian left-preorders and we state some properties they have. In Section 3 we define the space of Conradian left-preorders and we prove some topological properties. In Section 4 we focus on proving \Cref{teorema lista de equivalencias relativas con unicidad introduccion} and other related results about $\cord_\ce(\gr)$ following the ideas on \cite{Conradianfinito-rivas}. In Section 5 we prove \Cref{teorema lista de equivalencias absolutas con unicidad introduccion} and other technical results about $\cord(\gr)$, as well as some direct consequences. Finally, in Section 6, we briefly state and prove the dynamical version of Conradian left-preorders we have announced before.

\section{Left-preorders and positive cones}

The notion of left-preorder can be also found in other works. For example, it is present in \cite{valuated-decauprond} and \cite{locallymoving} with the name of \textit{left-invariant preorder}. We use this shorter terminology to abbreviate. This notion will be equivalent to the notion of relative order in \cite{morris-antolinrivas}. We avoid considering the trivial left-preorder. In this section we detail some definitions and basic properties related to left-preorders.\\

A \textit{preorder} $\lleq$ on a set $X$ is a transitive reflexive binary relation that is total (for all $x_1,x_2\in X$ we have that $x_1\lleq x_2$ or that $x_2\lleq x_1$). Denote by $[x]_{\lleq}$ the set of $y\in X$ such that $x\lleq y$ and $y\lleq x$. A preorder $\lleq$ on $X$ is \textit{trivial} if $x_1\lleq x_2$ for all $x_1,x_2\in X$. An \textit{order} on $X$ is a preorder that is antisymmetric.\\

Given a group $\gr$ acting on a set $X$, a relation $\rela$ on $X$ is \textit{invariant under the action of} $\gr$ if $x_1\rela x_2$ implies $\ga x_1\rela \ga x_2$ for all $x_1,x_2\in X$ and all $\ga\in \gr$. A left-preorder $\lleq$ on a group $\gr$ is a \textit{left-preorder} on $\gr$ that is non-trivial and that is invariant under the action of $\gr$ on itself by left multiplication $$\ga_1\lleq \ga_2\mbox{ implies }\ga \ga_1\lleq \ga \ga_2\mbox{  for all }\ga,\ga_1,\ga_2\in \gr.$$ A \textit{left-order} on $\gr$ is an antisymmetric left-preorder on $\gr$.\\

We now state some basic properties related to left-preorders. We left the proofs as execises since they are routinary.

\begin{lema}
Let $\gr$ be a group and $\lleq$ a left-preorder on $\gr$. Then $[1]_{\lleq}$ is a proper subgroup of $\gr$.
\end{lema}

\begin{lema}\label{lema preorden como orden en cocientes}
Let $\gr$ be a group. Given $\lleq$ a left-preorder on $\gr$, the relation $\lleqb$ on $\gr/[1]_{\lleq}$ defined by $$\ga[1]_{\lleq}\lleqb \gb[1]_{\lleq}\; \Leftrightarrow\; \ga\lleq \gb $$ is a total order invariant under left multiplication. Reciprocally, given $\ce$ a proper subgroup of $\gr$ and $\lleqb$ a total order on $\gr/\ce$ invariant under left multiplication, the relation $\lleq$ on $\gr$ defined by $$\ga\lleq \gb \; \Leftrightarrow\; \ga \ce\lleqb \gb \ce$$ is a left-preorder on $\gr$ satisfying $\ce=[1]_{\lleq}$.
\end{lema}
 
Notice that a left-preorder $\lleq$ on $\gr$ is a left-order on $\gr$ if and only if $[1]_{\lleq}=\{1\}$. In particular, last lemma is trivial for left-orders as one can identify $\gr$ with $\gr/\{1\}$. Moreover, when the respective subgroup $[1]_{\lleq}$ or $\ce$ is normal, then a total order invariant under left multiplication on $\gr/[1]_{\lleq}$ or on $\gr/\ce$ respectively corresponds to a left-order on such quotient group.

\begin{nota}
Through all this paper we will identify each left-preorder on $\gr$ with a total left-invariant order on the set of cosets over a subgroup of $\gr$. Given a left preorder $\lleq$ on $\gr$, the expression $\ga [1]_{\lleq}\lle \gb [1]_{\lleq}$ will mean that $\ga [1]_{\lleq}\lleq \gb [1]_{\lleq}$ and $\ga [1]_{\lleq}\neq \gb [1]_{\lleq}$.
\end{nota}

\begin{defi}
Let $\gr$ be a group and let $\ce$ be a proper subgroup of $\gr$. A left-preorder on $\gr$ \textit{relative to} $\ce$ is a left preorder $\lleq$ satisfying that $\ce=[1]_{\lleqb}$.
\end{defi}

We now present positive cones as a central notion on our work. We skip the proof of the following three results as they are routinary, but the ideas can be founded at the discussion of Definition 4 on \cite{relativelyconvex-antolinwarrenzoran}.

\begin{lema}\label{lema auxiliar productos previo a orden relativo como semigrupo}
Let $\gr$ be a group, let $\ce $ be a subgroup of $\gr$ and let $P$ is a subsemigroup of $\gr$. Suppose that $\gr=P\cup P\inv\cup \ce $ and that $P$ and $\ce $ are disjoint. Then, we have $ \ce P\subseteq P$ and $P\ce \subseteq P$.
\end{lema} 

\begin{lema}\label{lema previo orden relativo como semigrupo}
Let $\gr$ be a group, let $\ce $ be a subgroup of $\gr$ and let $P$ is a subsemigroup of $\gr$. Suppose that $\gr=P\cup P\inv\cup \ce $. Then, $P$ and $\ce $ are disjoint if and only if $P$, $P\inv $ and $\ce $ are pairwise disjoint.
\end{lema}

\begin{prop}\label{orden relativo como semigrupo}
Given $\lleq$ a \defrelativeorder{\gr}{\ce }, the set $P=\{\ga\in \gr\tq \ce \lle \ga \ce \}$ satisfies both
\begin{enumerate}
\item[a)] $P$ is a subsemigroup of $\gr$.
\item[b)] $\gr=P\cup P\inv\cup \ce $. 
\item[c)] $P$, $P\inv $ and $\ce $ are pairwise disjoint.
\end{enumerate}
and we call $P$ the positive cone of $\lle$. Reciprocally, if $\ce $ is a proper subgroup of $\gr$ and $P\cont \gr$ is a subset in such a way the last three conditions are satisfied, defining $$\ga\ce \lleq \gb\ce  \Leftrightarrow \ga\inv \gb\in P\cup \ce $$ give us a left-preorder on $\gr$ relative to $\ce $  whose positive cone is $P$.
\end{prop}

It is important also the notion of convexity.

\begin{defi}
Let $\gr$ be a group and $\ce$ be a proper subgroup. Let $\lleq$ be a left-preorder on $\gr$ relative to $\ce$. A subgroup $\sgr$ of $\gr$ is $\lleq$\textit{-convex} if for each $\ga,\gb,\gc\in\gr$ such that $\ga\ce\lleq\gb\ce\lleq\gc\ce$ and such that $\ga,\gc\in\sgr$, then $\gb\in\sgr$. In particular, $\ce\subseteq \sgr$.
\end{defi}

\begin{prop}\label{convexos estan linealmente ordenados}
Let $\gr$ be a group and $\ce$ be a proper subgroup. Let $\lleq$ be a left-preorder on $\gr$ relative to $\ce$. Then, the set of $\lleq$-convex subgroups is totally ordered by the inclusion.
\begin{proof}
It is routinary (analogous to the well known proof for the left-orders case).
\end{proof}
\end{prop}

We now define the functions $\mu_{\gr,\ce,\sgr}$ and $\rho_{\gr,\ce,\sgr}$. The function $\mu_{\gr,\ce,\sgr}$ give us a way to induce a left-preorder from two left-preorders. Reciprocally, $\rho_{\gr,\ce,\sgr}$ splits left-preorders on two left-preorders. We also study the behavior of these functions with the notion of convexity.

\begin{nota}
Let $A,B$ be two sets. We denote $A\sqcup B$ the union $A\cup B$ if we also have that $A$ and $B$ are disjoint.
\end{nota}

\begin{prop}\label{ordenar desde dos ordenes}
Let $\gr$ be a group and $\ce$ be a proper subgroup. Let $\sgr$ be a proper subgroup of $\gr$ such that $\ce\subsetneq\sgr$. Let $\lleq_1$ be a left-preorder on $\sgr$ relative to $\ce$ with positive cone $P_1$ and let $\lleq_2$ be a left preorder on $\gr$ relative to $\sgr$ with positive cone $P_2$. Then, the set $P_1\cup P_2$ defines a positive cone of a left-preorder on $\gr$ relative to $\ce$, which we denote by $\mu_{\gr,\ce,\sgr}(\lleq_1,\lleq_2)$.
\begin{proof}
By \Cref{orden relativo como semigrupo}, in order to prove that $P_1\cup P_2$ defines a positive cone of a left-preorder on $\gr$ relative to $\ce$, we need to prove that $P_1\cup P_2$ is a subsemigroup of $\gr$ and that $\gr=(P_1\sqcup P_2)\sqcup \ce \sqcup(P_1\sqcup P_2)\inv$. We are going to prove that each condition hold one by one.\\

Firstly, we prove that $P_1\cup P_2$ is a subsemigroup of $\gr$. Since  $\lleq_2$ is a left preorder on $\gr$ relative to $\sgr$ with positive cone $P_2$, combining \Cref{orden relativo como semigrupo}, \Cref{lema auxiliar productos previo a orden relativo como semigrupo} and the fact that $P_1\subseteq \sgr$, we deduce that
\begin{equation*}
 P_1 P_2\subseteq \sgr P_2\subseteq P_2 \quad \mbox{ and }\quad P_2 P_1 \subseteq P_2\sgr \subseteq P_2.
\end{equation*}
From this and the fact that $P_1$ and $P_2$ are both subsemigroups of $\gr$, we deduce that 
$$(P_1\cup P_2)(P_1\cup P_2)\subseteq P_1P_1\cup P_2P_2\cup P_1P_2\cup P_2P_1 \subseteq P_1\cup P_2.$$ This proves that $P_1\cup P_2$ is a subsemigroup of $\gr$, as we wanted.\\

Secondly, we show that $\gr=(P_1\sqcup P_2)\sqcup \ce \sqcup(P_1\sqcup P_2)\inv$. By the definition of $P_1$ and $P_2$ as positive cones, we have $$\sgr=P_1\sqcup \ce\sqcup P_1\inv\;\mbox{ and }\; \gr=P_2\sqcup \sgr\sqcup P_2\inv.$$ Hence, we have  $$\gr=P_2\sqcup \sgr\sqcup P_2\inv=P_2\sqcup P_1\sqcup \ce\sqcup P_1\inv\sqcup P_2\inv=(P_1\sqcup P_2)\sqcup \ce\sqcup (P_1\sqcup P_2)\inv.$$
\end{proof}
\end{prop}

\begin{prop}\label{convexos arriba ordenar desde dos ordenes}
Let $\gr$ be a group and $\ce$ be a proper subgroup. Let $\sgr$ be a proper subgroup of $\gr$ such that $\ce\subsetneq\sgr$. Let $\lleq_1$ be a left-preorder on $\sgr$ relative to $\ce$ with positive cone $P_1$ and let $\lleq_2$ be a left preorder on $\gr$ relative to $\sgr$ with positive cone $P_2$. Let $\lleq=\mu_{\gr,\ce,\sgr}(\lleq_1,\lleq_2)$, which is a left-preorder on $\gr$ relative to $\ce$ with positive cone $P_1\cup P_2$. Let $\sgrb$ be a subgroup of $\gr$ such that $\sgr\subseteq \sgrb$. Then, $\sgrb$ is $\lleq$-convex if and only if $\sgrb$ is $\lleq_2$-convex.

\begin{proof}
Suppose that $\sgrb$ is $\lleq$-convex. We are going to prove that $\sgrb$ is $\lleq_2$-convex. Indeed, consider $\ga,\gc\in\sgrb$ and $\gb\in\gr$ such that $$\ga\sgr\lleq_2\gb\sgr\lleq_2\gc\sgr.$$ Then we have three cases: $\ga\sgr=\gb\sgr$, $\gb\sgr=\gc\sgr$ or $\ga\sgr\lle_2\gb\sgr\lle_2\gc\sgr$. We are going to prove that, in all these cases, we have $\gb\in\sgrb$.\\

If we are in the case $\ga\sgr=\gb\sgr$, then $\ga\inv\gb$ is in $\sgr\subseteq \sgrb$, and since $\ga\in\sgrb$, we obtain $\gb\in \sgrb$. The case $\gb\sgr=\gc\sgr$ is analogous, so we also obtain $\gb\in \sgrb$.\\
 
Assume now we are in the case $\ga\sgr\lle_2\gb\sgr\lle_2\gc\sgr.$ As $P_2$ defines the positive cone of $\lleq_2$, we deduce that $\ga\inv\gb$ and $\gb\inv\gc$ are both in $P_2$. Since $P_1\cup P_2$ is the positive cone of $\lleq$, we deduce $$\ga\ce\lle\gb\ce\lle\gc\ce.$$ Using that $\sgrb$ is $\lleq$-convex by hypothesis, we deduce that $\gb\in\sgrb$. Therefore, we have $\gb\in\sgrb$ in the three cases, so we conclude that $\sgrb$ is $\lleq_2$-convex.\\

Reciprocally, assume that $\sgrb$ is $\lleq_2$-convex. We are going to prove that that $\sgrb$ is $\lleq$-convex. Indeed, take $\ga,\gc\in\sgrb$ and $\gb\in\gr$ such that $$\ga\ce\lleq\gb\ce\lleq\gc\ce.$$ From last equation, since $P_1\cup P_2$ defines the positive cone of $\lleq$, we obtain $\ga\inv\gb,\gb\inv\gc\in P_1\cup P_2\cup \ce$. Then, noticing that $P_1\cup \ce\subseteq \sgr$, we deduce that $\ga\inv\gb,\gb\inv\gc\in P_2\cup \sgr$. Since $P_2$ is the positive cone of the left-preorder $\lleq_2$ relative to $\sgr$, we deduce that $$\ga\sgr\lleq_2\gb\sgr\lleq_2\gc\sgr.$$ Using that $\sgrb$ is $\lleq_2$-convex by hypothesis and that $\ga,\gc\in\sgrb$, we obtain that $\gb\in\sgrb$. We conclude that $\sgrb$ is a $\lleq$-convex subgroup, as we wanted.
\end{proof}
\end{prop}

As a consequence, since $\sgr$ is clearly $\lleq_2$-convex, we have the following.

\begin{coro}\label{convexo entre medias ordenar desde dos ordenes}
Let $\gr$ be a group and $\ce$ be a proper subgroup. Let $\sgr$ be a proper subgroup of $\gr$ such that $\ce\subsetneq\sgr$. Let $\lleq_1$ be a left-preorder on $\sgr$ relative to $\ce$ with positive cone $P_1$ and let $\lleq_2$ be a left preorder on $\gr$ relative to $\sgr$ with positive cone $P_2$. Let $\lleq=\mu_{\gr,\ce,\sgr}(\lleq_1,\lleq_2)$, which is a left-preorder on $\gr$ relative to $\ce$ with positive cone $P_1\cup P_2$. Then, $\sgr$ is $\lleq$-convex.
\end{coro}

\begin{nota}
Let $\gr$ be a group and $\ce$ be a proper subgroup. Let $\lleq$ be a left-preorder on $\gr$ relative to $\ce$.  Let $\sgr$ be a subgroup of $\gr$. We denote by $\lleq_{|\sgr}$ the restriction of $\lleq$ to $\sgr$. The relation $\lleq_{|\sgr}$ on $\sgr$ is a left-preorder on $\gr$ relative to $\ce\cap\sgr$. Equivalently, as in \Cref{lema preorden como orden en cocientes}, when seeing $\lleq$ as a relation on $\gr/\ce$, $\lleq_{|\sgr}$ corresponds to the natural order induced by $\lleq$ on $\sgr/(\ce\cap\sgr)$. By definition, if $\lleq$ is Conradian, then $\lleq_{|\sgr}$ is Conradian.
\end{nota}

\begin{lema}\label{lema una de las proyecciones es la restriccion}
Let $\gr$ be a group and $\ce$ be a proper subgroup. Let $\sgr$ be a proper subgroup of $\gr$ such that $\ce\subsetneq\sgr$. Let $\lleq_1$ be a left-preorder on $\sgr$ relative to $\ce$ with positive cone $P_1$ and let $\lleq_2$ be a left preorder on $\gr$ relative to $\sgr$ with positive cone $P_2$. Let $\lleq=\mu_{\gr,\ce,\sgr}(\lleq_1,\lleq_2)$, which is a left-preorder on $\gr$ relative to $\ce$ with positive cone $P_1\cup P_2$. Then, $\lleq_1$ is equal to the restriction $\lleq_{|\sgr}$.
\begin{proof}
We only need to prove that, for every $\ga,\gb\in\sgr$, we have $\ga\ce\lleq \gb\ce$ if and only if $\ga\ce\lleq_1 \gb\ce$. Take $\ga,\gb\in\sgr$. Notice that $\ga\inv\gb\in \sgr$ and that $\sgr\cap P_2=\emptyset$. Knowng this, we have the chain of equivalences $$\ga\ce\lleq_1 \gb\ce \quad \Longleftrightarrow \quad \ga\inv\gb\in P_1\cup \ce \quad \Longleftrightarrow \quad \ga\inv\gb\in P_1\cup P_2\cup \ce \quad \Longleftrightarrow \quad \ga\ce\lleq \gb\ce.$$
\end{proof}
\end{lema}

\begin{prop}\label{convexos abajo ordenar desde dos ordenes}
Let $\gr$ be a group and $\ce$ be a proper subgroup. Let $\sgr$ be a proper subgroup of $\gr$ such that $\ce\subsetneq\sgr$. Let $\lleq_1$ be a left-preorder on $\sgr$ relative to $\ce$ with positive cone $P_1$ and let $\lleq_2$ be a left preorder on $\gr$ relative to $\sgr$ with positive cone $P_2$. Let $\lleq=\mu_{\gr,\ce,\sgr}(\lleq_1,\lleq_2)$, which is a left-preorder on $\gr$ relative to $\ce$ with positive cone $P_1\cup P_2$. Let $\sgrc$ be a subgroup of $\gr$ such that $\sgrc\subseteq \sgr$. Then, $\sgrc$ is $\lleq$-convex if and only if $\sgrc$ is $\lleq_1$-convex.

\begin{proof}
The fact that being $\lleq$-convex implies being $\lleq_1$-convex is a direct consequence of the fact that, by \Cref{lema una de las proyecciones es la restriccion}, $\lleq_1$ is equal to the restriction $\lleq_{|\sgr}$. Reciprocally, suppose that $\sgrc$ is $\lleq_1$-convex. We are going to prove that that $\sgrc$ is $\lleq$-convex. In order to do this, take $\ga,\gc\in\sgrc$ and $\gb\in\gr$ such that
\begin{equation*}
\ga\ce\lleq\gb\ce\lleq\gc\ce.
\end{equation*}
By \Cref{convexo entre medias ordenar desde dos ordenes}, $\sgr$ is $\lleq$-convex, so $\ga,\gc$ in $\sgrc\subseteq \sgr$ implies $\gb\in\sgr$. Hence, we have
\begin{equation*}
\ga\ce\lleq_{|\sgr}\gb\ce\lleq_{|\sgr}\gc\ce.
\end{equation*}
By \Cref{lema una de las proyecciones es la restriccion} we have that $\lleq_1$ is equal to $\lleq_{|\sgr}$, so applying that $\sgrc$ is $\lleq_1$-convex by hypothesis, we obtain $\gb\in\sgrc$. We conclude that $\sgrc$ is $\lleq$-convex.
\end{proof}
\end{prop}

\begin{prop}\label{dos ordenes desde orden}
Let $\gr$ be a group and $\ce$ be a proper subgroup. Let $\sgr$ be a proper subgroup of $\gr$ such that $\ce\subsetneq\sgr$. Let $\lleq$ be a left-preorder on $\sgr$ relative to $\ce$ with positive cone $P$. Suppose that $\sgr$ is $\lleq$-convex. Then,
\begin{itemize}
\item The set $P\cap \sgr$ defines a positive cone of a left-preorder $\lleq_1$ on $\sgr$ relative to $\ce$.
\item The set $P\settminus \sgr$ defines a positive cone of a left-preorder $\lleq_2$ on $\gr$ relative to $\sgr$.
\end{itemize}
We denote $\rho_{\gr,\ce,\sgr}(\lleq)=(\lleq_1,\lleq_2)$.
\begin{proof}
We are going to prove each of the two conditions of the statement separately. The first condition is rutinary using \Cref{orden relativo como semigrupo}.\\

We continue by proving the second condition, that is, we want to prove that $P\settminus \sgr$ defines a positive cone of a left-preorder $\lleq_2$ on $\gr$ relative to $\sgr$. In order to achieve this, by \Cref{orden relativo como semigrupo}, we only need to show that $P\settminus \sgr$ is a subsemigroup of $\gr$, that $\gr=(P\settminus \sgr)\cup \sgr \cup(P\settminus \sgr)\inv$, and that $P\settminus \sgr$, $(P\settminus \sgr)\inv$ and $\sgr$ are pairwise disjoint. We check each of this conditions one by one. The fact that $P\settminus \sgr$, $(P\settminus \sgr)\inv$ and $\sgr$ are pairwise disjoint is direct. The fact $$\gr=(P\settminus \sgr)\cup \sgr \cup(P\settminus \sgr)\inv$$ follows by $$\gr=P\cup\ce\cup P\inv\subseteq P\cup\sgr\cup P\inv=(P\settminus \sgr)\cup \sgr \cup(P\settminus \sgr)\inv,$$ as the other inclusion is trivial. The last condition we need to check is that $P\settminus \sgr$ is a subsemigroup. We are going to prove it. Take $\ga,\gb \in P\settminus \sgr$. We want to prove that $\ga\gb \in P\settminus \sgr$. In order to do that, we are going to show that $\ga\gb \in P$ but $\ga\gb \notin \sgr$. The fact that $$\ga\gb \in P$$ follows directly by the fact that $P$ is a semigroup and $$\ga \in P\;\mbox{ and }\; \gb \in P$$ Hence, we only need to prove that $\ga\gb \notin \sgr$. We do it by contradiction, so assume $\ga\gb \in \sgr$. Since $P$ is the positive cone of $\lleq$ and since $\ga,\gb \in P\settminus \sgr$, we deduce $$\ce\lle \ga\ce \;\mbox{ and }\; \ce\lle\gb\ce.$$ Then, we obtain $$\ce\lle \ga\ce\lle \ga\gb\ce.$$ Applying to this equation the fact that $\sgr$ is $\lleq$-convex by hypothesis, $1\in\sgr$ obviously and $\ga\gb \in \sgr$ by hypothesis, we deduce that $\ga\in\sgr$, a contradiction with $\ga\notin \sgr$. Therefore, we have checked that $P\settminus \sgr$ is a subsemigroup. Hence, we have proved all the needed conditions to show that $P\settminus \sgr$ defines a positive cone of a left-preorder $\lleq_2$ on $\gr$ relative to $\sgr$, as we wanted.
\end{proof}
\end{prop}

We finish the section by introducing some technicalities we will need in other sections.

\begin{lema}\label{convexo en convexo}
Let $\gr$ be a group and $\ce$ be a proper subgroup. Let $\lleq$ be a Conradian left-preorder on $\gr$ relative to $\ce$. Let $\sgr_1,\sgr_2$ be proper subgroups of $\gr$ such that $$\ce \subseteq \sgr_1\subseteq \sgr_2 \subseteq \gr.$$ If $\sgr_1$ is $\lleq_{|\sgr_2}$-convex and $\sgr_2$ is $\lleq$-convex, then $\sgr_1$ is $\lleq$-convex.
\begin{proof}
This follows from the definitions.
\end{proof}
\end{lema}

\begin{defi}
Let $\gr$ be a group and $\ce$ be a proper subgroup. Let $\lleq$ be a left-preorder on $\gr$ relative to $\ce$ with positive cone $P$. Let $\sgrb$ and $\sgrc$ be $\lleq$-convex subgroups of $\gr$ such that $\sgrb\subsetneq\sgrc$. We define the \textit{left-preorder induced by} $\lleq$ on $\sgrc$ relative to $\sgrb$ to be the left-preorder on $\sgrc$ relative to $\sgrb$ having positive cone $(P\cap \sgrc)\settminus\sgrb$. This is well defined by applying statement 1 from \Cref{dos ordenes desde orden} and then statement 2 from \Cref{dos ordenes desde orden}.
\end{defi}

\begin{obs}\label{observacion ordenes inducidos}
Let $\gr$ be a group and $\ce$ be a proper subgroup. Let $\lleq$ be a left-preorder on $\gr$ relative to $\ce$. Let $\sgrb$ and $\sgrc$ be $\lleq$-convex subgroups of $\gr$ such that $\sgrb\subsetneq\sgrc$. Notice that the notion of left-preorder induced by $\lleq$ on $\sgrc$ relative to $\sgrb$ is well defined. Indeed, applying statement 1 from \Cref{dos ordenes desde orden}, we deduce that  $P\cap \sgrc$ defines the positive cone of a left-preorder on $\sgrc$ relative to $\ce$, and applying statement 2 from \Cref{dos ordenes desde orden}, we conclude that $(P\cap \sgrc)\settminus\sgrb$ defines the positive cone of a left-preorder induced by $\lleq$ on $\sgrc$ relative to $\sgrb$. \\

In general, if we denote by $\hat{\lleq}$ the left-preorder induced by $\lleq$ on $\sgrc$ relative to $\sgrb$, then for each $\ga,\gb\in \sgrb$ we have  $$\ga\sgrb\hat{\lle}\gb\sgrb \quad \Longrightarrow \quad \ga\ce\lle\gb\ce.$$ This is a direct consequence of the fact that $(P\cap \sgrc)\settminus\sgrb\subseteq P$. In an equivalent way, using the properties of total orders, for each $\ga,\gb\in \sgrb$ we have $$\ga\ce\lleq\gb\ce \quad \Longrightarrow \quad \ga\sgrb\hat{\lleq}\gb\sgrb.$$
\end{obs}

\begin{lema}\label{ordenes inducidos ordenar desde dos ordenes}
Let $\gr$ be a group and $\ce$ be a proper subgroup. Let $\sgr$ be a proper subgroup of $\gr$ such that $\ce\subsetneq\sgr$. Let $\lleq_1$ be a left-preorder on $\sgr$ relative to $\ce$ with positive cone $P_1$ and let $\lleq_2$ be a left preorder on $\gr$ relative to $\sgr$ with positive cone $P_2$. Let $\lleq=\mu_{\gr,\ce,\sgr}(\lleq_1,\lleq_2)$, which is a left-preorder on $\gr$ relative to $\ce$ with positive cone $P_1\cup P_2$. Let $\sgrb$ and $\sgrc$ be $\lleq$-convex subgroups of $\gr$ such that $\sgrb\subsetneq\sgrc$. 
\begin{enumerate}
\item If $\sgrb\subsetneq\sgrc\subseteq \sgr$, then the left-preorder induced by $\lleq$ on $\sgrc$ relative to $\sgrb$ and the left-preorder induced by $\lleq_1$ on $\sgrc$ relative to $\sgrb$ are equal.
\item If $\sgr\subseteq\sgrb\subsetneq \sgrc$, then the left-preorder induced by $\lleq$ on $\sgrc$ relative to $\sgrb$ and the left-preorder induced by $\lleq_2$ on $\sgrc$ relative to $\sgrb$ are equal.
\end{enumerate}
\begin{proof}
We are going to prove each of the affirmations. We start by proving the first one, so assume that $\sgrb\subsetneq\sgrc\subseteq \sgr$. We are going to prove that the left-preorder induced by $\lleq$ on $\sgrc/\sgrb$ and the left-preorder induced by $\lleq_1$ on $\sgrc/\sgrb$ are equal by showing that both left-preorders have the same positive cone. Since $\sgrc\subseteq \sgr$ and $P_2$ is the positive cone of a left-preorder on $\gr$ relative to $\sgr$, we deduce that $\sgrc\cap P_2=\emptyset$. Then, this implies $$\left( (P_1\cup P_2)\cap \sgrc \right)\settminus \sgrb =\left( P_1\cap\sgrc \right)\settminus \sgrb .$$ Notice that, in the previous equation, the left term is the positive cone of the left-preorder induced by $\lleq$ on $\sgrc/\sgrb$ and the right term is the positive cone of the left-preorder induced by $\lleq_1$ on $\sgrc/\sgrb$, as we wanted to prove.\\

We continue by proving the second affirmation, so suppose that $\sgr\subseteq\sgrb\subsetneq \sgrc$. Since $\sgr\subseteq \sgrb$ and $P_1$ is the positive cone of a left-preorder on $\sgr$ relative to $\ce$, we deduce that $P_1\subseteq \sgrb$. Using ths fact, we have $$\left( (P_1\cup P_2)\cap \sgrc \right)\settminus \sgrb =\left( (P_1\cup P_2) \settminus \sgrb\right)\cap \sgrc =\left( P_2 \settminus \sgrb\right)\cap\sgrc=\left( P_2\cap\sgrc \right)\settminus \sgrb .$$ The left term of the last equation is the positive cone of the left-preorder induced by $\lleq$ on $\sgrc/\sgrb$ and the right term of the last equation is the positive cone of the left-preorder induced by $\lleq_2$ on $\sgrc/\sgrb$, as we wanted.
\end{proof}
\end{lema}

\section{Conradian left-preorders}

We now present the definition of Conradian left-preorder. Notice that, for $\ce=\{1\}$, we recover the classical definition of Conradian left-order. For a reference discussing Conradian left-orders, see Chapter 9 of \cite{orderedgroupstopology-clay}.

\begin{defi}
Let $\gr$ be a group and $\ce$ be a proper subgroup. A left-preorder $\lleq$ on $\gr$ relative to $\ce$ is \textit{Conradian} if for all $\ga,\gb\in \gr$ such that $\ce\lle\ga\ce$ and $\ce\lle\gb\ce$ there is $n\in\nat$ satisfying $\ga\ce\lle\gb\ga^n\ce$. We define $\cord_\ce(\gr)$ to be the \textit{set of Conradian left-preorders on} $\gr$ \textit{relative to} $\ce$ and $\cord(\gr)$ to be the \textit{set of Conradian left-preorders on} $\gr$.
\end{defi}

One can easily prove that the last definition is independent of either $0$ is considered to be in $\nat$ or not. For us $0\in\nat$.\\

We now want an equivalent definition of Conradian left-preorder using convex jumps. To consider the theorem with the equivalence, we need some auxiliary results. 

\begin{defi}
Let $\gr$ be a group and $\ce$ be a proper subgroup. Let $\lleq$ be a left-preorder relative to $\ce$. A $\lleq$\textit{-convex jump} is a pair $(\sgrb,\sgrc)$ of $\lleq$-convex subgroups of $\gr$ such that $\sgrb\subsetneq\sgrc$ and such that there is no $\lleq$-convex subgroup strictly between them.
\end{defi}

\begin{prop}\label{convexo maximal y minimal}
Let $\gr$ be a group and $\ce$ be a proper subgroup. Let $\lleq$ be a left-preorder relative to $\ce$. Consider $\ga\in \gr$. Define $$\sgrb_{\ga}=\bigcup_{\ga\notin \sgrb \mbox{ convex}}\sgrb \quad\mbox{and}\quad \sgrc_{\ga}=\bigcap_{\ga\in \sgrc \mbox{ convex}}\sgrc.$$ Then, $\sgrb_{\ga}$ and $\sgrc_{\ga}$ are $\lleq$-convex subgroups and $(\sgrb_{\ga},\sgrc_{\ga})$ is a $\lleq$-convex jump.
\begin{proof}
The set $\sgrb_{\ga}$ is a subgroup as it is an union of a linearly ordered family of subgroups, as a consequence of \Cref{convexos estan linealmente ordenados}. The fact that it is $\lleq$-convex is elementary and left as an exercise. It is also an elementary exercise to prove that $\sgrc_{\ga}$ is a $\lleq$-convex subgroup.\\

We finally prove that $(\sgrb_{\ga},\sgrc_{\ga})$ is a $\lleq$-convex jump. Since we know that $\sgrb_{\ga}$ and $\sgrc_{\ga}$ are $\lleq$-convex subgroups and since $\ga\notin \sgrb_{\ga}$ and $\ga\in\sgrc_{\ga}$, applying \Cref{convexos estan linealmente ordenados} we obtain that $$\sgrb_{\ga}\subsetneq \sgrc_{\ga}.$$ Hence, in order to show that $(\sgrb_{\ga},\sgrc_{\ga})$ is a $\lleq$-convex jump we only need to show that there is no $\lleq$-convex subgroup strictly between them. This will follow if we prove that any $\lleq$-convex subgroup $\sgr$ satisfies that $\sgr\subseteq\sgrb_{\ga}$ or that $\sgrc_{\ga}\subseteq\sgr$. However, this is clear since, by definition, we have $\sgr\subseteq\sgrb_{\ga}$ if $\ga\notin\sgr$ and $\sgrc_{\ga}\subseteq\sgr$ if $\ga\in\sgr$. We conclude that $(\sgrb_{\ga},\sgrc_{\ga})$ is a $\lleq$-convex jump.
\end{proof}
\end{prop}

\begin{coro}\label{lema convexo maximal}
Let $\gr$ be a group and $\ce$ be a proper subgroup. Let $\lleq$ be a Conradian left-preorder on $\gr$ relative to $\ce$. If $\gr$ is finitely generated, then there exists a maximal proper $\lleq$-convex subgroup.
\begin{proof}
Assume that $\gr$ is finitely generated, say $\gr=\gen{\ga_1,\dots,\ga_k}$ for $k\in\natp$. By possibly taking inverses, we can assume that $\ce \lleq \ga_i\ce$ for all $i=1,\dots,k$. Moreover, by possibly changing the indexes, we can assume that 
\begin{equation}\label{lema convexo maximal eq1}
\ce\lleq\ga_1\ce\lleq\dots\lleq\ga_k\ce.
\end{equation}
Define $$\sgrb=\bigcup_{\ga_k\notin \sgrb \mbox{ convex}}\sgrb \quad\mbox{and}\quad \sgrc=\bigcap_{\ga_k\in \sgrc \mbox{ convex}}\sgrc.$$  By \Cref{convexo maximal y minimal} we know that $(\sgrb,\sgrc)$ is a $\lleq$-convex jump. This means that $\sgrb$ and $\sgrc$ are $\lleq$-convex subgroups such that $\sgrb\subsetneq\sgrc$ and such that there is no $\lleq$-convex subgroups strictly between them. In particular, if we prove that $\sgrc=\gr$, we have that $\sgrb$ is a maximal proper $\lleq$-convex subgroup, as we wanted. Hence, we only need to show that $\sgrc=\gr$. We prove this by proving that $\ga_i\in\sgrc$ for $i=1,\dots,k$. Indeed, take $i\in\{1,\dots,k\}$. Notice that $\ga_k\in\sgrc$ by definition of $\sgrc$. From \Cref{lema convexo maximal eq1} we have $$\ce\lleq\ga_i\ce\lleq\ga_k\ce.$$ This, combined with the fact that $\ga_k\in\sgrc$ and that $\sgrc$ is $\lleq$-convex, implies that $\ga_i\in\sgrc$, as we wanted.
\end{proof}
\end{coro}

The condition of being finitely generated in the previous corollary cannot be avoided, as we see in next example.

\begin{ej}\label{ejemplo no fin gen no convexo maximal}
Consider the group $\gr=\bigoplus_{i\in \nat}\ent$. We can write $\gr=\bigoplus_{i\in \nat}\gen{\gx_i}$ where $\gen{\gx_i}$ is isomorphic to $\ent$. Define $$(\gx_i^{j_i}) \lle (\gx_i^{k_i}) \quad \Longleftrightarrow \quad j_{i_0}<k_{i_0} \;\mbox{ for }\, i_0 \mbox{ the 
maximum such that } j_{i_0}\neq k_{i_0}.$$ It is routine to check that $\lleq$ is a left-order on $\gr$. In other words, $\lleq$ is a left-preorder on $\gr$ relative to $\{1\}$. Since $\gr$ is abelian, we have that $\lleq$ is Conradian. One can check that the set of $\lleq$ convex subgroups is $$\{\gen{\gx_1,\dots,\gx_k} \tq k\in\nat\}\cup \gr.$$ Clearly, there is no maximal proper $\lleq$-convex subgroup.
\end{ej}

From now, we are going to adapt some classical ideas from the theory of orderable groups (we follow Chapter 9 of \cite{orderedgroupstopology-clay} as a reference). The following lemma follows the discussion on Lemma 9.7 from \cite{orderedgroupstopology-clay}.

\begin{lema}\label{proposicion Conradiano 1}
Let $\gr$ be a group and $\ce$ be a proper subgroup. Let $\lleq$ be a left-preorder relative to $\ce$. Let $\sgr$ be a maximal $\lleq$-convex proper subgroup of $\gr$. If $\lleq$ is Conradian, then for every $\ga,\gb\in\gr\settminus\sgr$ such that $\ce\lle\ga\ce$ and $\ce\lle\gb\ce$ there exists $n\in\nat$ such that $\gb\ce\lle\ga^n\ce$.
\begin{proof}
Assume that $\lleq$ is Conradian. Consider $\ga\in\gr\settminus\sgr$ such that $\ce\lle\ga\ce$. Define $$\sce=\{\gd\in\gr\tq \exists n\in\nat \mbox{ such that } \gd\ce\lle\ga^n\ce\}.$$ Notice that the proof will conclude if we show that every $\gb\in\gr\settminus\sgr$ such that $\ce\lle\gb\ce$ satisfies $\gb\in\sce$. In order to do that, we are going to consider the action of $\gr$ on the power set of $\gr$ by left-multiplication and we are going to show that the stabilizer of $S_0$ is $\lleq$-convex properly containing $\sgr$. Notice that, by definition of $\sce$, we have 
\begin{equation}\label{proposicion Conradiano 1 eq1}
\ga\sce=\sce.
\end{equation}
In other words, $\ga$ stabilizes $\sce$ under the left-multiplication action.\\

We start by showing that
\begin{equation}\label{proposicion Conradiano 1 eq2}
\sce\subseteq \gy\sce \mbox{ for all } \gy\in\gr \mbox{ such that } \ce\lleq\gy\ce,
\end{equation}
as we will use this fact many times along the proof. Indeed, consider $\gy\in\gr$ such that $\ce\lleq\gy\ce$. Notice that $\ce\lleq\ga^m\ce$ for all $m\in\nat$. Since $\lleq$ is Conradian, for all $m\in\nat$ there exists $n_m\in\nat$ such that $\ga^m\ce\lle\gy(\ga^m)^{n_m}\ce$, and hence $\gy\inv\ga^m\ce\lle(\ga^m)^{n_m}\ce$. Now, if we take $\gx\in\sce$, by definition of $\sce$, there is $m\in\nat$ such that $\gx\ce\lle\ga^m\ce$, so then $\gy\inv\gx\ce\lle\gy\inv\ga^m\ce$. Now, since we know that $\gy\inv\ga^m\ce\lle(\ga^m)^{n_m}\ce$, we deduce that $\gy\inv\gx\ce\lle(\ga^m)^{n_m}\ce$, and hence $\gy\inv\gx\in\sce$. Therefore, we have shown that $\gy\inv\gx\in\sce$ for all  $\gx\in\sce$, so then $\gx=\gy(\gy\inv\gx)\in\gy\sce$. This proves \cref{proposicion Conradiano 1 eq2}, as we wanted.\\

We consider $\gr_{\sce}$ the stabilizer of $\sce$ by the action on the power set of $\gr$ by left multiplication. We are now going to prove that $\gr_{\sce}$ is a $\lleq$-convex subgroup. Indeed, consider $\gx,\gy,\gz\in\gr$ such that $\gx\ce\lleq\gy\ce\lleq\gz\ce$ and such that $\gx,\gz\in\gr_{\sce}$. From $\gx\ce\lleq\gy\ce$ we deduce that $\ce\lleq\gx\inv\gy\ce$. Then, using \Cref{proposicion Conradiano 1 eq2} we deduce that $\sce\subseteq (\gx\inv\gy)\sce$. From this we deduce $\gx\sce\subseteq\gy\sce$, but $\gx\sce=\sce$ since $\gx\in\gr_{\sce}$, so then $\sce\subseteq\gy\sce$. Analogously, we can deduce that $\gy\sce\subseteq\sce$ from the facts $\gy\ce\lleq\gz\ce$ and $\gz\in\gr_{\sce}$. Therefore, we obtain $\gy\sce=\sce$, and hence $\gy\in\gr_{\sce}$. This shows that $\gr_{\sce}$ is a $\lleq$-convex subgroup, as we wanted.\\

We are now going to prove that $\sgr\subseteq \gr_{\sce}$. Indeed, let $\gee\in\sgr$. Recall that, at the beginning of the proof, we considered an element $\ga\in\gr\settminus\sgr$ such that $\ce\lle\ga\ce$. Since $\gee\in\sgr$ and $\sgr$ is $\lleq$-convex, we deduce $\ce\lleq\gee\mmen\ga\ce$. Then, applying \Cref{proposicion Conradiano 1 eq2}, we obtain $\sce\subseteq \gee\mmen\ga\sce$. Then, $\gee\mmen\sce\subseteq \ga\sce$, and combining this with \Cref{proposicion Conradiano 1 eq1}, we deduce $\gee\mmen\sce\subseteq \sce$. We have shown that $\gee\sce\subseteq \sce$ and $\gee\inv\sce\subseteq \sce$, or equivalently $\gee\sce\subseteq \sce$ and $\sce\subseteq \gee\sce$, so we deduce $\sce=\gee\sce$. This shows that $ \gee\in\gr_{\sce}$. Hence, we conclude that $\sgr\subseteq \gr_{\sce}$, as we wanted.\\

We now focus on showing that $\gr=\gr_{\sce}$. At this point, we have shown that $\gr_{\sce}$ is a $\lleq$-convex subgroup and that $\sgr\subseteq \gr_{\sce}$. Since $\sgr$ be a maximal $\lleq$-convex proper subgroup of $\gr$ by hypothesis, if we show that $\sgr\neq \gr_{\sce}$, we can conclude that $\gr=\gr_{\sce}$. Indeed, we know that $\alpha\notin \sgr$, but $\alpha\in \gr_{\sce}$ as a consequence of \Cref{proposicion Conradiano 1 eq1}, so this shows that $\sgr\neq \gr_{\sce}$. We conclude that $\gr=\gr_{\sce}$, as we wanted.\\

Finally, we are now able to complete the proof from the proven fact $\gr=\gr_{\sce}$. Take $\gb\in\gr\settminus\sgr$ such that $\ce\lle\gb\ce$. Since $\gr=\gr_{\sce}$, we deduce $\gb\in\gr_{\sce}$, so then $\gb\sce=\sce$. Since $1\ce=\ce\lle\gb\ce$, we obtain $1\in\sce$, so $\gb\in\gb\sce$, and hence $\gb\in\sce$. Therefore, by definition of $\sce$, we conclude that there exists $n\in\nat$ such that $\gb\ce\lle\ga^n\ce$. Since the element $\ga$ was fixed at the beginning as an arbitrary $\ga\in\gr\settminus\sgr$ such that $\ce\lle\ga\ce$, we have shown that for every $\ga,\gb\in\gr\settminus\sgr$ such that $\ce\lle\ga\ce$ and $\ce\lle\gb\ce$ there exists $n\in\nat$ such that $\gb\ce\lle\ga^n\ce$.
\end{proof}
\end{lema}

\begin{coro}\label{corolario Conradiano 1}
Let $\gr$ be a group and $\ce$ be a proper subgroup. Let $\lleq$ be a left-preorder relative to $\ce$ and let $(\sgrb,\sgrc)$ be a $\lleq$-convex jump. If $\lleq$ is Conradian, then for every $\ga,\gb\in\sgrc\settminus\sgrb$ such that $\ce\lle\ga\ce$ and $\ce\lle\gb\ce$ there exists $n\in\nat$ such that $\gb\ce\lle\ga^n\ce$.
\end{coro}

The following lemma is bassed on Lemma 9.8 from \cite{orderedgroupstopology-clay}.

\begin{lema}\label{lema Conradiano 2}
Let $\gr$ be a group and $\ce$ be a proper subgroup. Let $\lleq$ be a left-preorder relative to $\ce$. If $\lleq$ is Conradian and $(\sgrb,\sgrc)$ is a $\lleq$-convex jump, then $\sgrb$ is normal in $\sgrc$.
\begin{proof}
We want to show that $\gc \sgrb \gc\inv =\sgrb$ for all $\gc\in\sgrc$. We are going to show that $\gc \sgrb \gc\inv =\sgrb$ for all $\gc\in\sgrc\settminus\sgrb$ such that $\ce\lle \gc \ce$, as for $\gc \ce\lle \ce$ follows by taking the inverse and for $\gc\in\sgrb$ is trivial. \\

We start by proving that $\gc \sgrb \gc\inv \subseteq \sgrb$ by contradiction. Indeed, suppose that there is $\ga\in\sgrb$ such that $\gc\ga\gc\inv\notin\sgrb$. In particular, $\gc\ga\gc\inv\notin\ce$ as $\ce\subseteq \sgrb$, so we deduce that $\ce\lle \gc\ga\gc\inv\ce$ or $\gc\ga\gc\inv\ce\lle \ce$. If we are in the case $\ce\lle \gc\ga\gc\inv\ce$, by \Cref{corolario Conradiano 1} knowing that  $\ce\lle \gc \ce$, we deduce that there exists $n\in\nat$ such that $\gc \ce\lle (\gc\ga\gc\inv)^n\ce$, so then $\gc \ce\lle \gc\ga^n\gc\inv\ce$, and hence $\ga^{-n}\ce\lle\gc\inv\ce$, but since $\gc\inv\ce\lle \ce$ as $\ce\lle \gc \ce$, we obtain that $$\ga^{-n}\ce\lle\gc\inv\ce\lle \ce,$$ but $\ga^{-n},1\in\sgrb$ and $\sgrb$ is convex, so $\gc\inv\in\sgrb$, and then $\gc\in\sgrb$, a contradiction. If we are in the case $\gc\ga\gc\inv\ce\lle\ce$, then applying the previous case for $(\gc\ga\gc\inv)\inv=\gc\ga\inv\gc\inv$ we also conclude  $\gc\in\sgrb$, a contradiction. Therefore, we have shown that $\gc \sgrb \gc\inv \subseteq \sgrb$, as we wanted.\\

To finish the proof, we only need to show that $\sgrb \subseteq \gc \sgrb\gc\inv$. In order to do that, we claim that $\gc\inv\gd\gc\ce\lle\gc \ce$ for all $\gd\in\sgrb$. Consider $\gd\in\sgrb$. As we know that $\gc\in\sgrc\settminus\sgrb$ and $\ce\lle \gc \ce$, applying that $\sgrb$ is convex, $\sgrb\subseteq\sgrc$ and $\gd\in\sgrb$, we deduce that $\gd\ce\lle \gc\ce$, so then $\ce\lle\gd\inv\gc\ce$. Applying \Cref{corolario Conradiano 1} noticing that $\gd\inv\gc\in \sgrc\settminus\sgrb$, we deduce that there is $n\in\nat$ such that $\gc\ce\lle(\gd\inv\gc)^n\ce$. Then, $\ce\lle\gc\inv (\gd\inv\gc)^n\ce$, but we know that $\ce\lle\gc\ce$, so multiplying this two $\lleq$-positive elements we obtain $\ce\lle\gc\inv(\gd\inv\gc)^n\gc\ce$, that is, $\ce\lle(\gc\inv(\gd\inv\gc)\gc)^n\ce$. As $\lleq$-positivity is preserved by positive roots (note that $n>0$), we obtain $\ce\lle\gc\inv(\gd\inv\gc)\gc\ce$. Left-multiplying by $(\gc\inv\gd\inv\gc)\inv=\gc\inv\gd\gc$, we obtain $\gc\inv\gd\gc\ce\lle\gc \ce$, as we claimed. \\

We are now able to prove that $\sgrb \subseteq \gc \sgrb\gc\inv$ to conclude the proof using the proven claim saying that $\gc\inv\gd\gc\ce\lle\gc \ce$ for all $\gd\in\sgrb$. Indeed, take $\gd\in\sgrb$. Then, we have three cases: $\ce\lle \gc\inv\gd\gc\ce$, $\gc\inv\gd\gc\ce\lle \ce$ or $\gc\inv\gd\gc\ce=\ce$. Assume that we are in the case $\ce\lle \gc\inv\gd\gc\ce$. As $\gd \in\sgrb$ implies $\gd^l \in\sgrb$ for all $l\in\natp$, and applying the conclusion of the previous paragraph we  deduce that $\gc\inv\gd^l\gc\ce\lle\gc \ce$ for all $l\in\natp$. This is equivalent to $(\gc\inv\gd\gc)^l\ce\lle\gc \ce$ for all $l\in\natp$, so using the fact that $\ce\lle \gc\inv\gd\gc\ce$ and \Cref{corolario Conradiano 1}, we conclude that $\gc\inv\gd\gc\in \sgrb$. Then $\gd\in \gc \sgrb\gc\inv$, as we wanted. If we are in the case $\gc\inv\gd\gc\ce\lle \ce$, we can just apply the previous case for $(\gc\inv\gd\gc)\inv=\gc\inv\gd\inv\gc$, so we obtain $\gd\inv\in \gc \sgrb\gc\inv$, and hence $\gd\in \gc \sgrb\gc\inv$. Finally, if we are in the case $\gc\inv\gd\gc\ce=\ce$, as $\ce\subseteq\sgrb$, we have $\gc\inv\gd\gc\in \sgrb$, and hence $\gd\in \gc \sgrb\gc\inv$. Therefore, we have shown that $\sgrb \subseteq \gc \sgrb\gc\inv$, as we needed.
\end{proof}
\end{lema}

\begin{nota}
Let $\gr$ be a group. Let $\lleq$ be total order on $\gr$. An injective map $\df{\tau}{\gr}{\real}$ is said to be $\lleq$-preserving if for all $\ga,\gb\in\gr$ we have $$\ga\lleq\gb \quad\Longleftrightarrow \quad\tau(\ga)\leq\tau(\gb).$$  As a consequence of total orders' properties and injectivity, this is equivalent to say that for all $\ga,\gb\in\gr$ we have $$\ga\lle\gb \quad\Longleftrightarrow \quad\tau(\ga)<\tau(\gb).$$ Notice that implicitely we are considering in $\real$ the usual order.
\end{nota}

Also, we previously need the following well known theorem (for a reference see Theorem 2.6 of \cite{orderedgroupstopology-clay}).

\begin{teo}[Hölder's theorem]\label{teorema de holder}
Let $\gr$ be a group. Let $\lleq$ be a left-order on $\gr$. Assume that $\lleq$ is archimedean, meaning that for all $\ga,\gb\in\gr\settminus \{1\}$ such that $1\lle\ga$ and  $1\lle\gb$ there exists $n\in\nat$ such that $\gb\lle \ga^n$. Then, there exists a $\lleq$-preserving injective group morphism $\df{\tau}{\gr}{\real}$.
\end{teo}

We can now state the theorem with the equivalent version of being a Conradian left-preorder using convex jumps. We also prove that, in the definition of Conradian left-preorder, one can take $n=2$. This last fact was proven for Conradian left-orders on \cite{Conradianoparados-jimenez} but using other kind of arguments. The theorem is proved using the ideas on Theorem 9.5 from \cite{orderedgroupstopology-clay}.

\begin{teo}\label{Conradiano como cocientes de saltos}
Let $\gr$ be a group and $\ce$ be a proper subgroup. Let $\lleq$ be a left-preorder relative to $\ce$. Then, the following conditions are equivalent:
\begin{itemize}
\item[a)] For every $(\sgrb,\sgrc)$ $\lleq$-convex jump, $\sgrb$ is normal in $\sgrc$ and there exists a $\hat{\lleq}$-preserving injective morphism $\df{\tau_{\sgrb,\sgrc}}{\sgrc/\sgrb}{\real}$ for $\hat{\lleq}$ the left-preorder induced by $\lleq$ on $\sgrc$ relative to $\sgrb$.
\item[b)] For all $\ga,\gb\in \gr$ such that $\ce\lle\ga\ce$ and $\ce\lle\gb\ce$ we have $\ga\ce\lle\gb\ga^2\ce$.
\item[c)] $\lleq$ is Conradian.
\end{itemize}
\begin{proof}
We are going to prove that each condition implies the next one. We start by proving that a) implies b). Assume that for every $(\sgrb,\sgrc)$ $\lleq$-convex jump, $\sgrb$ is normal in $\sgrc$ and there exists a $\hat{\lleq}$-preserving injective morphism $\df{\tau_{\sgrb,\sgrc}}{\sgrc/\sgrb}{\real}$ for $\hat{\lleq}$ the left-preorder induced by $\lleq$ on $\sgrc$ relative to $\sgrb$. Let $\ga,\gb\in \gr$ be such that $\ce\lle\ga\ce$ and $\ce\lle\gb\ce$. Define the subgroups $$\sgrb_{\ga}=\bigcup_{\ga\notin \sgrb \mbox{ convex}}\sgrb \quad\mbox{and}\quad \sgrc_{\ga}=\bigcap_{\ga\in \sgrc \mbox{ convex}}\sgrc.$$ By \Cref{convexo maximal y minimal}, $\sgrb_{\ga}$ and $\sgrc_{\ga}$ are $\lleq$-convex subgroups and $(\sgrb_{\ga},\sgrc_{\ga})$ is a $\lleq$-convex jump. One can analogously define $\sgrb_{\gb}$ and $\sgrc_{\gb}$, so we also have that $\sgrb_{\gb}$ and $\sgrc_{\gb}$ are $\lleq$-convex and $(\sgrb_{\gb},\sgrc_{\gb})$ is a $\lleq$-convex jump. Using the hypothesis, $\sgrb_{\ga}$ is normal in $\sgrc_{\ga}$, $\sgrb_{\gb}$ is normal in $\sgrc_{\gb}$ and we can consider the $\hat{\lleq}$-preserving injective morphisms $$\df{\tau_{\ga}=\tau_{\sgrb_{\ga},\sgrc_{\ga}}}{\sgrc_{\ga}/\sgrb_{\ga}}{\real}\qquad\df{\tau_{\gb}=\tau_{\sgrb_{\gb},\sgrc_{\gb}}}{\sgrc_{\gb}/\sgrb_{\gb}}{\real}.$$ Applying \Cref{convexos estan linealmente ordenados}, we have that $\sgrb_{\ga}\subsetneq\sgrc_{\ga}\subseteq \sgrb_{\gb}\subsetneq\sgrc_{\gb}$, or that $\sgrb_{\gb}\subsetneq\sgrc_{\gb}\subseteq \sgrb_{\ga}\subsetneq\sgrc_{\ga}$ or that $\sgrb_{\ga}=\sgrb_{\gb}\subsetneq\sgrc_{\ga}=\sgrc_{\gb}$. Hence, we can consider $(\sgrb,\sgrc)$ the $\lleq$-convex jump that is bigger or equal than the other, and $\df{\tau_{\sgrb,\sgrc}}{\sgrc/\sgrb}{\real}$ the corresponding $\hat{\lleq}$-preserving injective morphism given by the hypothesis, where $\hat{\lleq}$ is the left-preorder induced by $\lleq$ on $\sgrc$ relative to $\sgrb$. Notice that $\ga,\gb\in \sgrc$. Since $\ce\lle\ga\ce$ and $\ce\lle\gb\ce$, by \Cref{observacion ordenes inducidos}, we deduce that $\sgrb\hat{\lleq}\ga\sgrb$ and $\sgrb\hat{\lleq}\gb\sgrb$. Then, as $\tau$ is order preserving, we deduce that $\tau(\ga\sgrb)\geq 0$ and $\tau(\gb\sgrb)\geq 0$. By definition, it is clear that $\ga\notin \sgrb_{\ga}$ and $\gb\notin \sgrb_{\gb}$, so then $\ga\notin \sgrb$ or $\gb\notin \sgrb$, and hence $\tau(\ga\sgrb)\neq 0$ or $\tau(\gb\sgrb)\neq 0$, but we know that $\tau(\ga\sgrb)\geq 0$ and $\tau(\gb\sgrb)\geq 0$, so then $\tau(\ga\sgrb)+\tau(\gb\sgrb)> 0$. Hence, we have $$\tau(\ga\inv\gb\ga^2\sgrb)=-\tau(\ga\sgrb)+\tau(\gb\sgrb)+2\tau(\ga\sgrb)=\tau(\ga\sgrb)+\tau(\gb\sgrb)> 0.$$ Then, using that $\tau$ is $\hat{\lleq}$-preserving, we obtain that $$\sgrb\hat{\lle}\ga\inv\gb\ga^2\sgrb.$$ By \Cref{observacion ordenes inducidos}, we conclude that $$\ce\lle\ga\inv\gb\ga^2\ce.$$ Therefore, we have $$\ga\ce\lle\gb\ga^2\ce.$$ This shows that a) implies b).\\

Clearly b) implies c). We now prove that c) implies a). Assume that $\lleq$ is Conradian. Let $(\sgrb,\sgrc)$ be a $\lleq$-convex jump. By \Cref{lema Conradiano 2}, we know that $\sgrb$ is normal in $\sgrc$. Hence, we can consider the group $\sgrc/\sgrb$. We can consider $\hat{\lleq}$ the left-preorder induced by $\lleq$ on $\sgrc$ relative to $\sgrb$, which is in fact a left-order on the group $\sgrc/\sgrb$. We claim that $\hat{\lleq}$ is an archimedean left-order. Indeed, take $\ga,\gb\in\sgrc\settminus\sgrb$ such that $\sgrb\hat{\lle}\ga\sgrb$ and $\sgrb\hat{\lle}\gb\sgrb$. By \Cref{observacion ordenes inducidos}, we have $\ce\lle\ga\ce$ and $\ce\lle\gb\ce$. Then, applying \Cref{corolario Conradiano 1}, we deduce that there exists $n\in\nat$ such that $\gb\ce\lle\ga^n\ce$. Then, by \Cref{observacion ordenes inducidos}, we deduce that there exists $n\in\nat$ such that $\gb\sgrb\hat{\lleq}\ga^n\sgrb$. From this by changing $n$ by $n+1$, since we know that $\sgrb\hat{\lle}\ga\sgrb$, we deduce that there exists $n\in\nat$ such that $\gb\sgrb\hat{\lle}\ga^n\sgrb$. This proves that $\hat{\lleq}$ is an archimedean left-order, as we claimed. Applying \Cref{teorema de holder}, there exists a $\hat{\lleq}$-preserving injective morphism $\df{\tau_{\sgrb,\sgrc}}{\sgrc/\sgrb}{\real}$, as we wanted. This proves that c) implies a).
\end{proof}
\end{teo}

We can use the previous theorem to show that the functions $\mu_{\gr,\ce,\sgr}$ and $\rho_{\gr,\ce,\sgr}$ preserve being Conradian. Before that, we prove that they preserve convex jumps.

\begin{lema}\label{lema saltos convexos separando en dos ordenes}
Let $\gr$ be a group and $\ce$ be a proper subgroup. Let $\sgr$ be a proper subgroup of $\gr$ such that $\ce\subsetneq\sgr$. Let $\sgrb, \sgrc$ be two subgroups of $\gr$. Let $\lleq_1$ be a left-preorder on $\sgr$ relative to $\ce$ with positive cone $P_1$ and let $\lleq_2$ be a left preorder on $\gr$ relative to $\sgr$ with positive cone $P_2$. Let $\lleq=\mu_{\gr,\ce,\sgr}(\lleq_1,\lleq_2)$, which is a left-preorder on $\gr$ relative to $\ce$ with positive cone $P_1\cup P_2$. Then, $(\sgrb,\sgrc)$ is a $\lleq_i$-convex jump for some $i=1,2$ if and only if $(\sgrb,\sgrc)$ is a $\lleq$-convex jump.
\begin{proof}
We prove each implication. Firstly, assume that $(\sgrb,\sgrc)$ a $\lleq_i$-convex jump for some $i=1,2$, and we are going to prove that $(\sgrb,\sgrc)$ a $\lleq$-convex jump. Indeed, the fact that $(\sgrb,\sgrc)$ is a $\lleq_i$-convex jump implies that $\sgrb$ and $\sgrc$ are $\lleq_i$-convex. Then, by \Cref{convexos abajo ordenar desde dos ordenes} or \Cref{convexos arriba ordenar desde dos ordenes}, we deduce that $\sgrb$ and $\sgrc$ are $\lleq$-convex. Notice that there is not a $\lleq$-convex subgroup strictly between them $\sgrb$ and $\sgrc$, as \Cref{convexos abajo ordenar desde dos ordenes} or \Cref{convexos arriba ordenar desde dos ordenes} would imply that there is a $\lleq_i$-convex subgroup strictly between $\sgrb$ and $\sgrc$, a contradiction with the fact that $(\sgrb,\sgrc)$ a $\lleq_i$-convex jump. Hence, we have shown that $(\sgrb,\sgrc)$ is a $\lleq$-convex jump.\\

Secondly, suppose that $(\sgrb,\sgrc)$ is a $\lleq$-convex jump, and we are going to show that $(\sgrb,\sgrc)$ is a $\lleq_i$-convex jump for some $i=1,2$. Since $(\sgrb,\sgrc)$ is a $\lleq$-convex jump, we know that there is no $\lleq$-convex subgroup strictly between $\sgrb$ and $\sgrc$. Then, as $\sgr$ is $\lleq$-convex by \Cref{convexo entre medias ordenar desde dos ordenes}, we deduce that $$\sgr\subseteq \sgrb\subsetneq\sgrc\quad \mbox{ or } \quad \sgrb\subsetneq\sgrc\subseteq \sgr.$$ If we are in the case $$\sgr\subseteq \sgrb\subsetneq\sgrc,$$ then \Cref{convexos abajo ordenar desde dos ordenes} implies that $\sgrb$ and $\sgrc$ are $\lleq_1$-convex. Moreover, in this case $(\sgrb,\sgrc)$ is a $\lleq_1$-convex jump as, if there is a $\lleq_1$-convex subgroup strictly between them, by \Cref{convexos abajo ordenar desde dos ordenes} we would obtain a $\lleq$-convex subgroup strictly between them, which is impossible. If we are in the case $$\sgrb\subsetneq\sgrc\subseteq \sgr,$$ we can do an analogous argumentation as in the other case but using \Cref{convexos arriba ordenar desde dos ordenes}. Therefore, we have shown that  $(\sgrb,\sgrc)$ is a $\lleq_i$-convex jump for some $i=1,2$, as we wanted.
\end{proof}
\end{lema}

\begin{prop}\label{Conradiano en extension y separacion de ordenes}
Let $\gr$ be a group and $\ce$ be a proper subgroup. Let $\sgr$ be a proper subgroup of $\gr$ such that $\ce\subsetneq\sgr$. Let $\lleq_1$ be a left-preorder on $\sgr$ relative to $\ce$ with positive cone $P_1$ and let $\lleq_2$ be a left preorder on $\gr$ relative to $\sgr$ with positive cone $P_2$. Let $\lleq=\mu_{\gr,\ce,\sgr}(\lleq_1,\lleq_2)$, which is a left-preorder on $\gr$ relative to $\ce$ with positive cone $P_1\cup P_2$. Then, $\lleq_1$ and $\lleq_2$ are Conradian if and only if $\lleq$ is Conradian.
\begin{proof}
We are going to prove each implication. Assume that $\lleq_1$ and $\lleq_2$ are Conradian. We are going to prove that $\lleq$ is Conradian. In order to do that, by \Cref{Conradiano como cocientes de saltos}, we want to prove that for every $(\sgrb,\sgrc)$ $\lleq$-convex jump, $\sgrb$ is normal in $\sgrc$ and there exists a $\hat{\lleq}$-preserving injective morphism $\df{\tau_{\sgrb,\sgrc}}{\sgrc/\sgrb}{\real}$ where $\hat{\lleq}$ is the left-preorder induced by $\lleq$ on $\sgrc$ relative to $\sgrb$. Take $(\sgrb,\sgrc)$ a $\lleq$-convex jump. By \Cref{lema saltos convexos separando en dos ordenes}, $(\sgrb,\sgrc)$ is a $\lleq_i$-convex jump for $i=1$ or $i=2$. Since $\lleq_i$ is Conradian by hypothesis and since we know that $(\sgrb,\sgrc)$ is a $\lleq_i$-convex jump, using \Cref{Conradiano como cocientes de saltos}, we deduce that $\sgrb$ is normal in $\sgrc$ and there exists an $\hat{\lleq_i}$-preserving injective group morphism $\df{\tau_{\sgr,\sgrb}}{\sgrc/\sgrb}{\real}$ where is $\hat{\lleq_i}$ is the left-preorder induced by $\lleq_i$ on $\sgrc$ relative to $\sgrb$. Applying \Cref{ordenes inducidos ordenar desde dos ordenes}, we deduce that $\sgrb$ is normal in $\sgrc$ and there exists an $\hat{\lleq}$-preserving injective group morphism $\df{\tau_{\sgr,\sgrb}}{\sgrc/\sgrb}{\real}$ for $\hat{\lleq}$ the left-preorder induced by $\lleq$ on $\sgrc$ relative to $\sgrb$, as we wanted. This proves that $\lleq$ is Conradian.\\

Conversely, we now prove the other implication. Suppose that $\lleq$ is Conradian. We are going to prove that $\lleq_1$ and $\lleq_2$ are Conradian. In other words, if we take an arbitrary $i\in\{1,2\}$, we are going to show that $\lleq_i$ is Conradian. In order to deduce that, by \Cref{Conradiano como cocientes de saltos}, we just need to prove that for every $(\sgrb,\sgrc)$ $\lleq_i$-convex jump, $\sgrb$ is normal in $\sgrc$ and there exists an $\hat{\lleq_i}$-preserving injective group morphism $\df{\tau_{\sgr,\sgrb}}{\sgrc/\sgrb}{\real}$ where $\hat{\lleq_i}$ is the left-preorder induced by $\lleq_i$ on $\sgrc$ relative to $\sgrb$. Indeed, take $(\sgrb,\sgrc)$ a $\lleq_i$-convex jump. From \Cref{lema saltos convexos separando en dos ordenes}, we deduce that $(\sgrb,\sgrc)$ is a $\lleq$-convex jump. As $\lleq$ is Conradian by hypothesis and $(\sgrb,\sgrc)$ is a $\lleq$-convex jump, by \Cref{Conradiano como cocientes de saltos}, we deduce that $\sgrb$ is normal in $\sgrc$ and there exists an $\hat{\lleq}$-preserving injective group morphism $\df{\tau_{\sgr,\sgrb}}{\sgrc/\sgrb}{\real}$ for $\hat{\lleq}$ the left-preorder induced by $\lleq$ on $\sgrc$ relative to $\sgrb$. Using \Cref{ordenes inducidos ordenar desde dos ordenes}, we deduce that, in fact, $\df{\tau_{\sgr,\sgrb}}{\sgrc/\sgrb}{\real}$ is a $\hat{\lleq_i}$-preserving injective group morphism where is $\hat{\lleq_i}$ is the left-preorder induced by $\lleq_i$ on $\sgrc$ relative to $\sgrb$. With this, we can conclude that $\lleq_i$ is Conradian, as we wanted.
\end{proof}
\end{prop}

We can now describe $\mu_{\gr,\ce,\sgr}$ and $\rho_{\gr,\ce,\sgr}$ as mutually bijective functions defined on sets consisting of Conradian left-preorders.

\begin{prop}\label{biyeccion Conradiana general n=2}
Let $\gr$ be a group and $\ce$ be a proper subgroup. Let $\sgr$ be a proper subgroup of $\gr$ such that $\ce\subsetneq\sgr$. Consider the following maps described, identifying each left-preorder with its positive cone, by 
\begin{equation*}
\begin{array}{r@{\hspace{0pt}}c@{\hspace{0pt}}
c@{\hspace{4pt}}l}
\mu_{\gr,\ce,\sgr}\colon & \cord_\ce(\sgr)\times \cord_\sgr(\gr) &\longrightarrow& \{\lleq\in \cord_\ce(\gr) \tq \sgr \mbox{ is }\lleq-\mbox{convex in }\gr\} \\

& (P_1,P_2)\ \ &\mapsto& P_1\cup P_2
\end{array}
\end{equation*}
\begin{equation*}
\begin{array}{r@{\hspace{0pt}}c@{\hspace{0pt}}
c@{\hspace{4pt}}l}
\rho_{\gr,\ce,\sgr} \colon & \{\lleq\in \cord_\ce(\gr) \tq \sgr \mbox{ is }\lleq-\mbox{convex in }\gr\} &\longrightarrow& \cord_\ce(\sgr)\times \cord_\sgr(\gr) \\

& P\ \ &\mapsto& (P\cap \sgr, P\settminus \sgr)
\end{array}
\end{equation*}
Then, these two maps $\mu_{\gr,\ce,\sgr}$ and $\rho_{\gr,\ce,\sgr}$ are bijective and $\mu_{\gr,\ce,\sgr}=(\rho_{\gr,\ce,\sgr})\inv$.
\begin{proof}
The map $\mu_{\gr,\ce,\sgr}$ is well defined by \Cref{ordenar desde dos ordenes} and \Cref{Conradiano en extension y separacion de ordenes}, and $\rho_{\gr,\ce,\sgr}$ is well defined by \Cref{dos ordenes desde orden} and \Cref{Conradiano en extension y separacion de ordenes}. In order to prove that $\mu_{\gr,\ce,\sgr}$ and $\rho_{\gr,\ce,\sgr}$ are bijective and $\mu_{\gr,\ce,\sgr}=(\rho_{\gr,\ce,\sgr})\inv$, we only need to show that $\mu_{\gr,\ce,\sgr}\co\rho_{\gr,\ce,\sgr}$ and $\rho_{\gr,\ce,\sgr}\co \mu_{\gr,\ce,\sgr}$ are the identity maps on the corresponding sets. We are going to prove it.\\

We first prove that $\mu_{\gr,\ce,\sgr}\co\rho_{\gr,\ce,\sgr}$ is the identity map. Indeed, given $P$ the positive cone of an element in $\cord_\ce(\gr)$ having $\sgr$ as a convex subgroup, we have $$\left( \mu_{\gr,\ce,\sgr}\co\rho_{\gr,\ce,\sgr}\right)(P)=\mu_{\gr,\ce,\sgr}(P\cap \sgr, P\settminus \sgr)=(P\cap \sgr) \cup (P\settminus \sgr)= P.$$ Therefore, $\mu_{\gr,\ce,\sgr}\co\rho_{\gr,\ce,\sgr}$ is the identity map.\\

We lastly prove that  $\rho_{\gr,\ce,\sgr}\co \mu_{\gr,\ce,\sgr}$ is the identity map. Indeed, take $P_1$ the positive cone of an element in $\cord_\ce(\sgr)$ and $P_2$ the positive cone of an element in $\cord_\sgr(\gr)$. Notice that $$P_1\subseteq \sgr\;\mbox{ and }\; P_2\cap \sgr=\emptyset.$$ Using this and following the definitions we have $$\left(\rho_{\gr,\ce,\sgr}\co \mu_{\gr,\ce,\sgr}\right)(P_1, P_2)=\rho_{\gr,\ce,\sgr}(P_1\cup P_2)=\left( (P_1\cup P_2)\cap \sgr, (P_1\cup P_2)\settminus \sgr \right)=(P_1,P_2).$$ Therefore, $\rho_{\gr,\ce,\sgr}\co \mu_{\gr,\ce,\sgr}$ is the identity map.
\end{proof}
\end{prop}

We can generalize the last proposition inductively.

\begin{prop}\label{biyeccion Conradiana general para todo n}
Let $\gr$ be a group and $\ce$ be a proper subgroup. Let $$\ce=\gr_0 \subsetneq \gr_1\subsetneq \dots \subsetneq \gr_n=\gr$$ be a finite chain of proper subgroups of $\gr$ for $n\in\natp$. Consider the following maps, described by identifying each left-preorder with its positive cone, which are defined by
\begin{equation*}
\begin{array}{r@{\hspace{0pt}}c@{\hspace{0pt}}
c@{\hspace{4pt}}l}
\mu\colon & \prod_{j=1}^{n}\cord_{\gr_{j-1}}(\gr_j) &\longrightarrow& \{\lleq\in \cord_\ce(\gr) \tq \gr_i \mbox{ are }\lleq-\mbox{convex for all } i\} \\

& (P_1,\dots, P_n)\ \ &\mapsto& \bigcup_{j=1}^n P_j
\end{array}
\end{equation*}
\begin{equation*}
\begin{array}{r@{\hspace{0pt}}c@{\hspace{0pt}}
c@{\hspace{4pt}}l}
\rho\colon &  \{\lleq\in \cord_\ce(\gr) \tq \gr_i \mbox{ is }\lleq-\mbox{convex for all } i\} &\longrightarrow& \prod_{j=0}^{n-1}\cord_{\gr_{j-1}}(\gr_j) \\

& P\ \ &\mapsto&  \left((P\cap \gr_{i})\settminus \gr_{i-1}\right)_{i=1}^n
\end{array}
\end{equation*}
Then, these two maps $\mu$ and $\rho$ are bijective and $\mu=\rho\inv$.
\end{prop}

\section{The space of Conradian left-preorders}

We now introduce the space of Conradian left-preorders. The space of Conradian left-preorders is a subspace of the space of left-preorders. For the defition of the space of left-preorders, see \cite{morris-antolinrivas} for instance. The space of left-preorders plays an important role in \cite{locallymoving}, where it is shown that a certain natural quotient on a product of the space of left-preorders it is homeomorphic to a space of normalized harmonic actions. It is remarkable that space of normalized harmonic actions used in other important results: in \cite{deroinhurtado} it is shown that an irreducible lattice in a real semi-simple Lie group of real rank at least two and finite center is not left-orderable. Also, the space of left-preorders of certain abelian groups, which coincides with the space of Conradian left-preorders, has some connections with algebraic geometry (the \textit{Zariski-Riemann space}), as it is considered in \cite{alggeo1} to obtain completions of real fans, and also in \cite{alggeo2} to study toric varieties.\\

The definition of the space of Conradian left-preorders is analogous to the one given for the space of left-preorders (see \cite{morris-antolinrivas}), which follows the ideas of the space of left-orders considered by Sikora in \cite{topology-sikora}. 

\begin{nota}
Given $\gr$ a group and $\ce$ a proper subgroup, \Cref{lema previo orden relativo como semigrupo} provides us a way to identify left-preorders with elements of $\{\pm 1,0\}^\gr$. Specifically, given $\lleq$ a left-preorder on $\gr$ relative to $\ce $ and denoting by $P$ its positive cone, we can define a $\phi_\lleq\in \{\pm1,0\}^\gr$ by:
\begin{equation*}
\begin{array}{r@{\hspace{0pt}}c@{\hspace{0pt}}
c@{\hspace{4pt}}l}
\phi_\lleq\colon &\gr &\longrightarrow& \{\pm1,0\} \\

&\ga \ \ &\mapsto&\begin{cases}1 \quad \textrm{if} \quad \ga \in P \\
 -1 \quad \textrm{if}\quad \ga \in P\inv \\
 0 \quad \textrm{otherwise} \end{cases}
\end{array}
\end{equation*}
In particular, we can also consider this identification restricted for Conradian left-preorders, so we can identify the set $\cord(\gr)$ of Conradian left-preorders on $\gr$ with a subset of $\{\pm1,0\}^\gr$. The \textit{space of Conradian left-preorders of} $\gr$ is the set $\cord(\gr)$ endowed with the subspace topology of $\{\pm 1,0\}^\gr$ with the product topology. The \textit{space of Conradian left-preorders of} $\gr$ \textit{relative to} $\ce $ is the set $\cord_\ce (\gr)$ endowed with the subspace topology of $\cord(\gr)$.
\end{nota}

\begin{obs}
A basis of open neighborhoods of $\phi\in\cord(\gr)$ is given by the sets $$U_{g_1,\dots,g_n}(\phi)=\{\psi\in\cord(\gr) \tq \psi(g_i)=\phi(g_i)\}$$ for all $g_1,\dots, g_n\in \gr$ and all $n\in \natp$. Notice that the sets $$\overline{U}_{g_1,\dots,g_n}(\phi)=\{\psi\in\cord(\gr) \tq \psi(g_i)\neq\phi(g_i)\}$$ are also open. This remark is also true changing $\cord(\gr)$ by $\cord_\ce (\gr)$.
\end{obs}

The following two propositions are direct as these properties are true on $\{\pm1,0\}^\gr$ and they are hereditary for subspaces.

\begin{prop}\label{espacio de ordenes es totalmente disconexo}
The spaces $\cord(\gr)$ and $\cord_\ce (\gr)$ are totally disconnected and Hausdorff.
\end{prop}

\begin{prop}\label{espacio de ordenes en contable es espacio metrico}
If $\gr$ is countable, then the space $\{\pm1,0\}^\gr$ endowed with the product topology is metrizable. In particular, if $\gr$ is countable, then $\cord(\gr)$ and $\cord_\ce(\gr)$ are metrizable spaces.
\end{prop}

\begin{nota}
Define the following subsets of $\{\pm1,0\}^\gr$
\begin{equation*}
V_{\ga }^{r}=\{\psi\in \{\pm1,0\}^\gr \tq \psi(\ga )=r\} \quad \textrm{for each} \; \ga \in \gr\; \textrm{and each }\; r\in\{\pm1,0\},
\end{equation*}
\begin{equation*}
\overline{V}_{\ga }^{r}=\{\psi\in \{\pm1,0\}^\gr\tq \psi(\ga )\neq r\}=\bigcup_{k\in\{\pm 1, 0\}\settminus\{r\}}V_{\ga }^{k} \quad \textrm{for each} \; \ga \in \gr\; \textrm{and each }\; r\in\{\pm1,0\}.
\end{equation*}
Notice that they are simultaneously open and closed sets.
\end{nota}

\begin{prop}\label{espacio de ordenes son cerrados y compactos 1}
The space $\cord_\ce (\gr)$ is a closed subspace of $\{\pm1,0\}^\gr$. In particular, $\cord_\ce (\gr)$ is a compact space.
\begin{proof}
We prove that $\cord_\ce (\gr)$ is closed by describing it as intersections and finite unions of open sets. Consider the open sets $V_{g}^{r}$ and $\overline{V}_{g}^{r}$ of $\{\pm1,0\}^\gr$ as above. In order to obtain this description we will use \Cref{orden relativo como semigrupo} and \Cref{Conradiano como cocientes de saltos}. Indeed, from this we deduce that having an element of $\cord_\ce (\gr)$ is equivalent to give a decomposition of the form $\gr=P\sqcup P\inv \sqcup \ce $ where $P$ is a subsemigroup of $\gr$ and where $\ga,\gb\in P$ implies $\ga\inv\gb\ga^2\in P$. This is equivalent to give an element $\phi$ of $\{\pm1,0\}^\gr$ satisfying the following conditions:
\begin{enumerate}
\item $\left(\phi\inv(\{1\})\right)\inv=\phi\inv(\{-1\})$.
\item $\phi\inv(\{0\})=\ce $.
\item $\phi\inv(\{1\})$ is a subsemigroup of $\gr$.
\item For all $\ga,\gb\in \gr$, we have that $\ga,\gb\in\phi\inv(\{1\})$ implies $\ga\inv\gb\ga^2\in\phi\inv(\{1\})$.
\end{enumerate}
This gives us the description
\begin{equation*}
\begin{split}
\{\pm1,0\}^\gr \settminus \cord_\ce (\gr)=&\left((\bigcup_{\ga \in \gr}V_{\ga \inv}^{1}\cap \overline{V}_{\ga }^{-1}))\cup (\bigcup_{\ga \in \gr}\overline{V}_{\ga \inv}^{1}\cap V_{\ga }^{-1}))\right)\cup\left((\bigcup_{\ga \in \ce }\overline{V}_{\ga }^{0}))\cup (\bigcup_{\ga \in \ce \settminus \gr} V_{\ga }^{0}))\right)\cup\\ & \cup\left(\bigcup_{\ga ,\gb \in \gr}V_{\ga }^{1}\cap V_{\gb }^{1}\cap \overline{V}_{\ga \gb }^{1}\right)\cup\left(\bigcup_{\ga ,\gb \in \gr}V_{\ga }^{1}\cap V_{\gb }^{1}\cap \overline{V}_{\ga\inv\gb\ga^2 }^{1}\right)
\end{split}
\end{equation*}
where each of the four sets between parentheses corresponds to the respective negation of each condition bellow. Therefore, we have shown that $\cord_\ce (\gr)$ is a closed subset.\\

Finally, as a consequence of Tychonoff's theorem, $\{\pm1,0\}^\gr$ is compact, so its closed subspace $\cord_\ce (\gr)$ is compact.
\end{proof}
\end{prop}

\begin{prop}\label{espacio de ordenes son cerrados y compactos 2}
Suppose that $\gr$ is finitely generated. Then, the space $\cord(\gr)$ is a closed subspace of $\{\pm1,0\}^\gr$. In particular, $\cord(\gr)$ is a compact space.
\begin{proof}
We are going to prove it in an analogous way as in \Cref{espacio de ordenes son cerrados y compactos 1}. Since $\gr$ is finitely generated, we can take $\gr=\gen{g_1,\dots,g_n}$. Consider the open sets $V_{g}^{r}$ and $\overline{V}_{g}^{r}$ of $\{\pm1,0\}^\gr$ as above.\\

By \Cref{orden relativo como semigrupo} and \Cref{Conradiano como cocientes de saltos}, having a element of $\cord(\gr)$ is equivalent to give a decomposition of the form $\gr=P\sqcup P\inv \sqcup \ce $ where $P$ is a subsemigroup of $\gr$ satisfying that  $\ga,\gb\in P$ implies $\ga\inv\gb\ga^2\in P$ and where $\ce $ is a proper subgroup. This is equivalent to give an element $\phi\in \{\pm1,0\}^\gr$ satisfying the following:
\begin{enumerate}
\item $\left(\phi\inv(\{1\})\right)\inv=\phi\inv(\{-1\})$.
\item $\phi\inv(\{0\})$ is a subgroup.
\item $\phi\inv(\{1\})$ is a subsemigroup of $\gr$.
\item For all $\ga,\gb\in \gr$, we have that $\ga,\gb\in\phi\inv(\{1\})$ implies $\ga\inv\gb\ga^2\in\phi\inv(\{1\})$.
\item $g_i\notin\phi\inv(\{0\})$ for some $i\in\{1,\dots, n\}$.
\end{enumerate}
This gives us the description
\begin{equation*}
\begin{split}
\{\pm1,0\}^\gr \settminus \cord_\ce (\gr)=&\left((\bigcup_{\ga \in \gr}V_{\ga \inv}^{1}\cap \overline{V}_{\ga }^{-1}))\cup (\bigcup_{\ga \in \gr}\overline{V}_{\ga \inv}^{1}\cap V_{\ga }^{-1}))\right)\cup \left(\bigcup_{\ga ,\gb \in \gr}V_{\ga }^{0}\cap V_{\gb }^{0}\cap \overline{V}_{\ga \gb\inv }^{0}\right) \cup\\ & \cup\left(\bigcup_{\ga ,\gb \in \gr}V_{\ga }^{1}\cap V_{\gb }^{1}\cap \overline{V}_{\ga \gb }^{1}\right)\cup\left(\bigcup_{\ga ,\gb \in \gr}V_{\ga }^{1}\cap V_{\gb }^{1}\cap \overline{V}_{\ga\inv\gb\ga^2 }^{1}\right)\cup  \left(\bigcap_{i=1, \dots, n}V_{g_i}^{0}  \right)
\end{split}
\end{equation*}
Therefore $\cord(\gr)$ is closed on $\{\pm1,0\}^\gr$. In particular, $\cord(\gr)$ is a compact space.
\end{proof}
\end{prop}

\begin{obs}
In the previous proposition we need the condition of $\gr$ being finitely generated, as $\cord(\bigoplus_{i\in \nat}\ent)$ is not closed on $\{\pm1,0\}^{\bigoplus_{i\in \nat}\ent}$. Indeed, since $\bigoplus_{i\in \nat}\ent$ is abelian, the space of Conradian left-preorders is equal to the space of left-preorders, and in \cite{morris-antolinrivas} it is shown that the space of left-preorders of $\bigoplus_{i\in \nat}\ent$ is not closed on $\{\pm1,0\}^{\bigoplus_{i\in \nat}\ent}$.
\end{obs}

\begin{obs}\label{observacion cantor}
Corollary 2-98 from \cite{topology-hockingyoung}, we know that a compact, Hausdorff and totally disconnected metric space that has no isolated points is isomorphic to the Cantor set. Then, by the previous results, if $\gr$ is finitely generated and $\cord(\gr)$ has no isolated elements, then $\cord(\gr)$ is a Cantor set. In the same way, if $\gr$ is countable and $\cord_\ce (\gr)$ has no isolated elements, then $\cord_\ce (\gr)$ is a Cantor set.
\end{obs}

If $\gr$ is a group and $\ce$ is a proper normal subgroup, then the space $\cord_\ce (\gr)$ of Conradian left-orders on $\gr$ relative to $\ce$ is equal to the space of Conradian left-orders on $\gr/\ce$. Indeed, as we justified in \Cref{lema preorden como orden en cocientes}, recall that we think of $\cord_\ce (\gr)$ as the set of total orders on $\gr/\ce$ that are invariant under left multiplication and that satisfy the Conradian condition, and this is precisely the definition of the set of Conradian left-orders on $\gr/\ce$. Notice that the respective endowed topology are also build in the same way, so both spaces are the same.

\begin{ej}\label{ejemplo estructura en cocientes abelianos}
Let $\gr$ be a group and $\ce$ be a proper normal subgroup such that $\gr/\ce$ is torsion-free abelian. Then, $\cord_\ce (\gr)$ is equal to the space of Conradian left-orders on $\gr/\ce$, but every left-order on an abelian group is Conradian, so then $\cord_\ce (\gr)$ is equal to the space of left-orders on the group $\gr/\ce$. By \cite{topology-sikora}, we know that the space of left-orders on a non-trivial torsion-free abelian group is a discrete space with two elements if the abelian rank is $1$ and it is a Cantor set (having no isolated elements)  if the abelian rank is greater than $1$. Hence, we conclude the following:
\begin{enumerate}
\item If $\gr/\ce$ is an torsion-free abelian group of abelian rank $1$, then $\cord_{\ce}(\gr)$ is a discrete space with two elements.
\item If $\gr/\ce$ is an torsion-free abelian group of abelian rank greater than $1$, then $\cord_{\ce}(\gr)$ has no isolated orders and it is a Cantor set.
\end{enumerate}
\end{ej}

We conclude the section with a technical lemma we will need in the next section.

\begin{lema}\label{lema aislado con convexos implica aislado en un cociente}
Let $\gr$ be a group and $\ce$ be a proper subgroup. Let $\gr_1,\gr_2$ be proper subgroups of $\gr$ such that $$\ce \subsetneq \gr_1\subsetneq \gr_2 \subsetneq \gr.$$ If $\cord_\ce(\gr)$ has an isolated element $\lleq$ such that $\gr_1$ and $\gr_2$ are $\lleq$-convex subgroups, then $\cord_{\gr_1}(\gr_2)$ has an isolated element.
\begin{proof}
We are going to prove this by contraposition. In this way, assume that $\cord_{\gr_1}(\gr_2)$ has no isolated elements. We are going to prove that there is no isolated element on $\cord_\ce(\gr)$ having $\gr_1$ and $\gr_2$ as convex subgroups. In order to do that, take $\lleq$ a Conradian left-preorder on $\gr$ relative to $\ce$ having $\gr_1$ and $\gr_2$ as convex subgroups, and denote by $P$ its positive cone. We now show that $\lleq$ is not isolated in $\cord_\ce(\gr)$, that is, we now prove that for every finite collection of elements belonging to $P$, there is a Conradian left-preorder $\bar{\lleq}$ on $\gr$ relative to $\ce$ such that all this collection is on the positive cone of $\bar{\lleq}$ but $\bar{\lleq}\neq\lleq$. We prove now this fact.\\

Take $\ga_1,\dots,\ga_k\in P$ a finite collection. We are going to consider the mutually inverse bijections $\mu$ and $\rho$ from \Cref{biyeccion Conradiana general para todo n} for the finite chain of subgroups $$\ce \subsetneq \gr_1\subsetneq \gr_2 \subsetneq \gr.$$ Since $\gr_1$ and $\gr_2$ are $\lleq$-convex subgroups, we can consider $$\rho(\lleq)=(\lleq_1,\lleq_2,\lleq_3) \;\mbox{ where }\; \lleq_1\in \cord_\ce(\gr_1)\;\mbox{, }\; \lleq_2\in \cord_{\gr_1}(\gr_2) \;\mbox{ and }\; \lleq_3\in \cord_{\gr_2}(\gr).$$ Consider $P_1$ the positive cone of $\lleq_1$, $P_2$ the positive cone of $\lleq_2$ and $P_3$ the positive cone of $\lleq_3$. By the definition the mutually inverse bijections $\mu$ and $\rho$ given in \Cref{biyeccion Conradiana general para todo n}, we have
\begin{equation}\label{lema aislado con convexos implica aislado en un cociente eq1}
P=P_1\cup P_2\cup P_3.
\end{equation}
By possibly reordering the elements $\ga_1,\dots,\ga_k\in P$, we can assume $$\ga_1,\dots,\ga_m\in P_2\;\mbox{ and }\;\ga_{m+1},\dots,\ga_k\in P\settminus  P_2 .$$ Applying \Cref{lema aislado con convexos implica aislado en un cociente eq1}, we obtain 
\begin{equation}\label{lema aislado con convexos implica aislado en un cociente eq2}
\ga_1,\dots,\ga_m\in P_2\;\mbox{ and }\;\ga_{m+1},\dots,\ga_k\in P_1\cup P_3.
\end{equation}
 Since $\ga_1,\dots,\ga_m$ is a finite (possibly empty) collection of $\lleq_2$-positive elements where $\lleq_2\in\cord_{\gr_1}(\gr_2)$, and since $\cord_{\gr_1}(\gr_2)$ has no isolated elements by hypothesis, we can take $\bar{\lleq_2}\in\cord_{\gr_1}(\gr_2)$ with positive cone $\bar{P_2}$ and satisfying that 
\begin{equation}\label{lema aislado con convexos implica aislado en un cociente eq3}
\ga_1,\dots,\ga_m\in \bar{P_2}
\end{equation}
and $\bar{\lleq_2}$ is not equal to $\lleq_2$. Define $\bar{\lleq}=\mu(\lleq_1,\bar{\lleq_2},\lleq_3)$, which belongs to $\cord_\ce(\gr)$. As $\mu$ and $\rho$ mutually inverse bijections by \Cref{biyeccion Conradiana general para todo n}, we obtain 
\begin{equation}\label{lema aislado con convexos implica aislado en un cociente eq4}
\rho(\bar{\lleq})=(\lleq_1,\bar{\lleq_2},\lleq_3)
\end{equation}
By the definition the mutually inverse bijections $\mu$ and $\rho$ given in \Cref{biyeccion Conradiana general para todo n}, we know that  
\begin{equation}\label{lema aislado con convexos implica aislado en un cociente eq5}
\bar{P}=P_1\cup  \bar{P_2}\cup P_3.
\end{equation}
From \Cref{lema aislado con convexos implica aislado en un cociente eq2} and \Cref{lema aislado con convexos implica aislado en un cociente eq3} we obtain $$\ga_1,\dots,\ga_k\in P_1\cup  \bar{P_2}\cup P_3.$$ Applying \Cref{lema aislado con convexos implica aislado en un cociente eq5}, we deduce that 
\begin{equation}\label{lema aislado con convexos implica aislado en un cociente eq6}
\ga_1,\dots,\ga_k\in\bar{P}.
\end{equation}
The map $\rho$ is bijective by \Cref{biyeccion Conradiana general para todo n}, and we know that $\bar{\lleq_2}$ is not equal to $\lleq_2$, so from \Cref{lema aislado con convexos implica aislado en un cociente eq3} and \Cref{lema aislado con convexos implica aislado en un cociente eq3}, we deduce that $\bar{\lleq}$ is not equal to $\lleq$. Therefore, $\bar{\lleq}$ is a Conradian left-preorder on $\gr$ relative to $\ce$ such that $\bar{\lleq}$ is not equal to $\lleq$ and, by \Cref{lema aislado con convexos implica aislado en un cociente eq6}, the elements $\ga_1,\dots,\ga_k$ are on the positive cone of $\bar{\lleq}$, as we wanted to prove. Hence, $\lleq$ is not isolated in $\cord_\ce(\gr)$, as we wanted.
\end{proof}
\end{lema}

\section{Groups having finitely many Conradian left-preorders relative to a given subgroup}

The main objective of this section is proving \Cref{teorema lista de equivalencias relativas con unicidad introduccion} from the introduction, which is stated here as \Cref{teorema lista de equivalencias relativas con unicidad}. We start by some technical results, and then we prove \Cref{teorema lista de equivalencias relativas}, which is a partial version of the theorem we have mention. Then, we improve the partial version with other auxiliary results, obtaining the needed theorem. We conclude the section with two corollaries of the main theorem. The section uses the ideas from \cite{Conradianfinito-rivas}.

\begin{defi}
Let $\gr$ be a group. A finite series $$\gr_0 \nor \gr_1\nor \cdots \nor \gr_n$$ on $\gr$ is said to be \textit{rational} if each $\gr_i$ is normal on $\gr_{i+1}$ and the group $\gr_{i+1}/\gr_i$ is torsion-free abelian of abelian rank $1$ for all $0\leq i<n$. Also, we say that such a series starts in $\gr_0$, ends in $\gr_n$ and has length $n$. 
\end{defi}

\begin{defi}
Let $\gr$ be a group and $\ce$ be a proper subgroup. Let $$\gr_0 \nor \gr_1\nor \cdots \nor \gr_n$$ be a finite rational series on $\gr$. An \textit{abelian jump} on this series is a subseries the form $$\gr_i\nor\gr_{i+1}\nor\gr_{i+2}$$ for $i\in\{0,\dots, n-2\}$ such that $\gr_i$ is normal on $\gr_{i+2}$ and the group $\gr_{i+2}/\gr_i$ is abelian.
\end{defi}

\begin{lema}\label{lema cociente no abeliano}
Let $\gr$ be a group. Take $$\gr_0 \nor \gr_1\nor \cdots \nor \gr_n=\gr$$ a rational series on $\gr$ without abelian jumps. Then, for each $0\leq i<n-2$ there exists $\ga\in\gr_{i+1}$ and $\gb\in\gr_{i+2}\settminus \gr_{i+1}$ such that $$\gb\inv\ga\gb\gr_i\neq \ga\gr_i$$ on the group $\gr_{i+1}/\gr_i$.
\begin{proof}
Take $0\leq i<n-2$. We are going to prove it by contradiction, so assume that 
\begin{equation}\label{ecuacion 1 lema cociente no abeliano}
\gb\inv\ga\gb\gr_i= \ga\gr_i    \mbox{ for all } \ga\in\gr_{i+1} \mbox{ and all } \gb\in\gr_{i+2}\settminus \gr_{i+1}.
\end{equation}
Notice that this equation makes sense on the group $\gr_{i+1}/\gr_i$ because $\gb\inv\ga\gb\in\gr_{i+1}$ as $\gr_{i+1}$ is normal on $\gr_{i+2}$. We are going to prove that $$\gr_i\nor\gr_{i+1}\nor\gr_{i+2}$$ is an abelian jump. In order to do that, we are going to prove that $\gr_i$ is normal on $\gr_{i+2}$ and that the group $\gr_{i+2}/\gr_i$ is abelian.\\

We first prove that $\gr_i$ is normal on $\gr_{i+2}$. We are going to prove it by showing that for all $\gc\in\gr_i$ and all $\gb\in\gr_{i+2}$ we have $\gb\inv\gc\gb\in \gr_i$. Take $\gc\in\gr_i$ and $\gb\in\gr_{i+2}$. Since $\gr_i\subseteq\gr_{i+1}$, we have $\gc\in\gr_{i+1}$, so by \Cref{ecuacion 1 lema cociente no abeliano}, we deduce that $\gb\inv\gc\gb\gr_i= \gc\gr_i$. Notice that $\gc\gr_i=\gr_i$ as $\gc\in\gr_i$, so then $\gb\inv\gc\gb\gr_i=\gr_i$. Hence, $\gb\inv\gc\gb\in \gr_i$. Therefore, we have shown that $\gr_i$ is normal on $\gr_{i+2}$.\\

We now prove that the group $\gr_{i+2}/\gr_i$ is abelian. We are going to prove it by showing that for all $\gb_1,\gb_2\in\gr_{i+2}$ we have $\gb_1\gb_2\gr_i=\gb_2\gb_1\gr_i$. Consider $\gb_1,\gb_2\in\gr_{i+2}$. If $\gb_1\in\gr_{i+1}$ or $\gb_2\in\gr_{i+1}$ then $\gb_1\gb_2\gr_i=\gb_2\gb_1\gr_i$ follows by \Cref{ecuacion 1 lema cociente no abeliano} when only one of the conditions holds, or by the fact that $\gr_{i+1}/\gr_i$ is abelian when both holds. Hence, we can assume that $\gb_1,\gb_2\in\gr_{i+2}\settminus \gr_{i+1}$. As $\gr_{i+2}/\gr_{i+1}$ is torsion-free abelian of rank $1$, we know that there exists $\gb_0\in \gr_{i+2}\settminus \gr_{i+1}$ and $r,s\in\ent\settminus \{0\}$ such that $\gb_0^r\gr_{i+1}=\gb_1\gr_{i+1}$ and $\gb_0^s\gr_{i+1}=\gb_2\gr_{i+1}$. Then, there exists $\ga_1,\ga_2\in\gr_{i+1}$ such that $\gb_0^r\ga_1=\gb_1$ and $\gb_0^s\ga_2=\gb_2$. By \Cref{ecuacion 1 lema cociente no abeliano}, we deduce that $\gb_0\gr_i$ commutes with $\ga_1\gr_i$ and with $\ga_2\gr_i$ on the group $\gr_{i+2}/\gr_i$. Also notice that $\ga_1\gr_i$ and $\ga_2\gr_i$ commute as they are elements of the abelian group $\gr_{i+1}/\gr_i$. Hence, we have $$\gb_1\gb_2\gr_i=\gb_0^r\ga_1\gb_0^s\ga_2\gr_i=\gb_0^{r+s}\ga_1\ga_2\gr_i=\gb_0^{r+s}\ga_2\ga_1\gr_i=\gb_0^s\ga_2\gb_0^r\ga_1\gr_i=\gb_2\gb_1\gr_i,$$ as we wanted. Therefore, $\gr_{i+2}/\gr_i$ is an abelian group.\\
Hence, we have shown that $\gr_i$ is normal on $\gr_{i+2}$ and that the group $\gr_{i+2}/\gr_i$ is abelian, as we wanted.
\end{proof}
\end{lema}

Assuming the notation of the previous lemma, if $\ga\in\gr_i$, then $\gb\inv\ga\gb$ is an element in $\gr_{i+1}\settminus\gr_i$ satisfying the same properties. In particular, by changing $\ga$ by  $\gb\inv\ga\gb$ if necessary, we deduce the following version of the previous lemma.

\begin{lema}\label{corolario lema cociente no abeliano}
Let $\gr$ be a group. Take $$\gr_0 \nor \gr_1\nor \cdots \nor \gr_n$$ a rational series on $\gr$ without abelian jumps. Then, for each $0\leq i<n-2$ there exists $\ga\in\gr_{i+1}\settminus\gr_i$ and $\gb\in\gr_{i+2}\settminus \gr_{i+1}$ such that $$\gb\inv\ga\gb\gr_i\neq \ga\gr_i$$ on the group $\gr_{i+1}/\gr_i$.
\end{lema}

\begin{prop}\label{convexidad de series}
Let $\gr$ be a group and $\ce$ be a proper subgroup. Assume that $\gr$ admits a rational series $$\ce=\gr_0 \nor \gr_1\nor \cdots \nor \gr_n=\gr$$ on $\gr$ starting on $\gr$, ending on $\ce$ and without abelian jumps. Then, $\gr_j$ is $\lleq$-convex for all $\lleq$ Conradian left-preorder on $\gr$ relative to $\ce$ and for all $0\leq j\leq n$.
\begin{proof}
Suppose that $\gr$ admits a rational series $$\ce=\gr_0 \nor \gr_1\nor \cdots \nor \gr_n=\gr$$ without abelian jumps. Let $\lleq$ be a Conradian left-preorder on $\gr$ relative to $\ce$. We want to prove that $\gr_j$ is $\lleq$-convex for all $0\leq j\leq n$. If we prove that $\gr_j$ is $\lleq_{|\gr_{j+1}}$-convex for all $0\leq j<n$, then the fact that $\gr_j$ is $\lleq$-convex for all $0\leq j\leq n$ follows directly by applying finitely many times \Cref{convexo en convexo}. Hence, we are going to prove that $\gr_j$ is $\lleq_{|\gr_{j+1}}$-convex for all $0\leq j<n$ in order to prove the statement. We are going to prove it by induction on $j$.\\

For $j=0$ it is trivial, as $\gr_0=\ce$ is always convex on every left-preorder relative to $\ce$. Given $0< j<n$ and assuming  that $\gr_{j-1}$ is $\lleq_{|\gr_j}$-convex as inductive hypothesis, we are going to prove that $\gr_j$ is $\lleq_{|\gr_{j+1}}$-convex. We do it by contradiction, so suppose that $\gr_j$ is not $\lleq_{|\gr_{j+1}}$-convex. Hence, there exists $\ga_1,\ga_2\in\gr_j$ and $\gb_0\in \gr_{j+1}\settminus\gr_j$ such that
\begin{equation}\label{eq hipotesis convexidad de series}
\ga_1\ce\lleq\gb_0\ce\lleq\ga_2\ce
\end{equation}\\

By \Cref{corolario lema cociente no abeliano}, we know that there exists $\tilde{\ga}\in\gr_j\settminus\gr_{j-1}$ and $\tilde{\gb}\in\gr_{j+1}\settminus \gr_j$ such that
$$\gb\inv\tilde{\ga}\tilde{\gb}\gr_{j-1}\neq \tilde{\ga}\gr_{j-1}.$$
Then, using that $\tilde{\gb}\inv\tilde{\ga}\tilde{\gb},\tilde{\ga}\in\gr_j$ and that $\gr_j/\gr_{j-1}$ is an torsion-free abelian group of rank $1$, we deduce that 
\begin{equation}\label{eq previa convexidad de series}
\tilde{\gb}\inv\tilde{\ga}^r\tilde{\gb}\gr_{j-1}= \tilde{\ga}^s\gr_{j-1} \mbox{ for some } s,r\in\ent \mbox{ such that } s\neq r
\end{equation}
Consider the subgroup $\sgrc=\gen{\ga_1,\ga_2,\tilde{\ga},\gb_0,\tilde{\gb}}$.  Until the end of the proof, as an abuse of notation, we will identify $\lleq$ with its restriction to $\sgrc$, and this restriction is Conradian left-preorder on $\sgrc$. By \Cref{lema convexo maximal}, since $\sgrc$ is finitely generated, we can consider $\sgrb$ a maximal proper $\lleq$-convex subgroup of $\sgrc$.  In other words, we have that $(\sgrb,\sgrc)$ is a $\lleq$-convex jump. Then, by \Cref{Conradiano como cocientes de saltos}, we know that $\sgrb$ is normal in $\sgrc$ and there exists a $\hat{\lleq}$-preserving injective morphism $\df{\tau}{\sgrc/\sgrb}{\real}$ for $\hat{\lleq}$ the left-preorder induced by $\lleq$ on $\sgrc$ relative to $\sgrb$. Notice that
\begin{equation}\label{eq 0 convexidad de series}
\gc_1\ce\lleq\gc_2\ce \quad\Longrightarrow\quad \tau(\gc_1\sgrb )\leq \tau(\gc_2\sgrb )\quad \mbox{ for all } \gc_1\gc_2\in\sgrc
\end{equation}
as a consequence of \Cref{observacion ordenes inducidos} and $\tau$ being  $\hat{\lleq}$-preserving.\\

Firstly, we claim that
\begin{equation}\label{eq 1 convexidad de series}
\tau(\gc\sgrb )=0 \mbox{ for all } \gc\in\gr_{j-1}\cap\sgrc
\end{equation}
Indeed, take $\gc\in\gr_{j-1}\cap\sgrc$. Define $\abs{\tilde{\ga}}$ to be $\tilde{\ga}$ if $\ce\lle\tilde{\ga}\ce$ and $\tilde{\ga}\inv$ if $\tilde{\ga}\ce\lle \ce$. By inductive hypothesis, $\gr_{j-1}$ is $\lleq_{|\gr_j}$-convex, so using the fact that $\gc^m\in\gr_{j-1}$ for all $m\in \ent$ and that $\abs{\tilde{\ga}},\abs{\tilde{\ga}}\inv\in\gr_j\settminus\gr_{j-1}$, we deduce that $$\abs{\tilde{\ga}}\inv\ce\lleq\gc^m\ce\lleq \abs{\tilde{\ga}}\ce \mbox{ for all } m\in\ent.$$ Then, applying \Cref{eq 0 convexidad de series} for $\abs{\tilde{\ga}},\gc\in\sgrc$, we obtain $$-\tau(\abs{\tilde{\ga}}\sgrb )\leq m\tau(\gc\sgrb )\leq \tau(\abs{\tilde{\ga}}\sgrb ) \mbox{ for all } m\in\ent.$$ It follows that $\tau(\gc\sgrb )=0$, as we claimed.\\

Secondly, we claim that 
\begin{equation}\label{eq 2 convexidad de series}
\tau(\ga\sgrb )=0 \mbox{ for all } \ga\in\gr_j\cap\sgrc
\end{equation}
We are now going to prove it. Consider $\ga\in\gr_j\cap\sgrc$. By \Cref{eq previa convexidad de series}, we deduce that $\tilde{\gb}\inv\tilde{\ga}^r\tilde{\gb}=\tilde{\ga}^s\tilde{\gc}$ for some $\tilde{\gc}\in\gr_{j-1}$ where $r\neq s$. Notice that $\tilde{\gc}\in\gr_{j-1}\cap\sgrc$, so by \Cref{eq 1 convexidad de series} we obtain $\tau(\tilde{\gc}\sgrb )=0$. Therefore, we deduce $$r\tau(\tilde{\ga}\sgrb )=-\tau(\tilde{\gb}\sgrb )+r\tau(\tilde{\ga}\sgrb )+\tau(\tilde{\gb}\sgrb )=\tau(\tilde{\gb}\inv\tilde{\ga}^r\tilde{\gb}\sgrb )=\tau(\tilde{\ga}^s\tilde{\gc}\sgrb )=s\tau(\tilde{\ga}\sgrb )+\tau(\tilde{\gc}\sgrb )=s\tau(\tilde{\ga}\sgrb ),$$ but $r\neq s$, so 
\begin{equation}\label{eq 3 convexidad de series}
\tau(\tilde{\ga}\sgrb )=0.
\end{equation}
Since $\ga\in\gr_j$, $\tilde{\ga}\in\gr_j\settminus\gr_{j-1}$ and since $\gr_j/\gr_{j-1}$ is an torsion-free abelian group of rank $1$, we deduce that there exists $l,m\in\ent$ with $l\neq 0$ such that $$\ga^l\gr_{j-1}= \tilde{\ga}^m\gr_{j-1}.$$ Hence, there exists $\gc_0\in\gr_{j-1}$ such that $\ga^l= \tilde{\ga}^m\gc_0$. Notice that $\gc_0\in\gr_{j-1}\cap\sgrc$, so applying \Cref{eq 1 convexidad de series} we deduce $\tau(\gc_0\sgrb )=0$. From \Cref{eq 3 convexidad de series}, from $\tau(\gc_0\sgrb )=0$ and from $\ga^l= \tilde{\ga}^m\gc_0$ we deduce that $$l\tau(\ga\sgrb )=\tau(\ga^l\sgrb )=\tau(\tilde{\ga}^m\gc_0\sgrb )=m\tau(\tilde{\ga}\sgrb )+\tau(\gc_0\sgrb )=0.$$ Since $l\neq 0$, we obtain $\tau(\ga\sgrb )=0$, as we claimed.\\

From \Cref{eq hipotesis convexidad de series} and from \Cref{eq 0 convexidad de series}, we deduce that $$\tau(\ga_1\sgrb )\leq\tau(\gb_0\sgrb )\leq\tau(\ga_2\sgrb ).$$ Notice that $\tau(\ga_1\sgrb )=0$ and $\tau(\ga_2\sgrb )=0$ follows from \Cref{crossing segunda ecuacion} since $\ga_1,\ga_2\in\gr_j\cap\sgrc$. Therefore, we deduce that 
\begin{equation}\label{eq 4 convexidad de series}
\tau(\gb_0\sgrb )=0.
\end{equation}
Since $\gb_0,\tilde{\gb}\in \gr_{j+1}$ and $\gr_{j+1}/\gr_j$ is an torsion-free abelian group of rank $1$, then there exists $l,m\in\ent$ with $l,m\neq 0$ such that $$\tilde{\gb}^l \gr_j=\gb_0^m \gr_j.$$ Then, there is $\hat{\ga}\in \gr_j$ such that $\tilde{\gb}^l=\gb_0^m\hat{\ga}$. In particular, we deduce that $\hat{\ga}\in \gr_j\cap \sgrc$ and that so then $$l\tau(\tilde{\gb}\sgrb )=\tau(\tilde{\gb}^l\sgrb )=\tau(\gb_0^m\hat{\ga}\sgrb )=m\tau(\gb_0\sgrb )+\tau(\hat{\ga}\sgrb ).$$ Since $\tau(\gb_0\sgrb)=0$ from \Cref{eq 4 convexidad de series} and since $\tau(\hat{\ga}\sgrb )=0$ from \Cref{eq 2 convexidad de series} and $\hat{\ga}\in\gr_j\cap\sgrc$, we obtain $l\tau(\tilde{\gb}\sgrb )=0$. As $l\neq 0$, we obtain 
\begin{equation}\label{eq 5 convexidad de series}
\tau(\tilde{\gb}\sgrb )=0.
\end{equation}
Finally, recall that $\sgrc=\gen{\ga_1,\ga_2,\tilde{\ga},\gb_0,\tilde{\gb}}$. From \Cref{eq 3 convexidad de series} we have $$\tau(\tilde{\ga}\sgrb )=0,$$ from \Cref{eq 4 convexidad de series} we have $$\tau(\gb_0\sgrb )=0,$$ from \Cref{eq 2 convexidad de series} we have $$\tau(\ga_1\sgrb )=\tau(\ga_2\sgrb )=\tau(\tilde{\ga}\sgrb )=0$$ since $\ga_1,\ga_2,\tilde{\ga}\in \gr_j\cap\sgrc$, and from \Cref{eq 5 convexidad de series} we have $$\tau(\tilde{\gb}\sgrb )=0.$$ Therefore, the morphism $\df{\tau}{\sgrc/\sgrb}{\real}$ sends a generating set of elements to $0$, so it is the trivial group morphism. Hence, as $\tau$ is injective, we obtain $\sgrb=\sgrc$. This is a contradiction with the fact that $\sgrb\neq\sgrc$ as $(\sgrb,\sgrc)$ is a $\lleq$-convex jump.
\end{proof}
\end{prop}

\begin{prop}\label{numero de ordenes Conradianos desde serie}
Let $\gr$ be a group and $\ce$ be a proper subgroup. Assume that $\gr$ admits a rational series $$\ce=\gr_0 \nor \gr_1\nor \cdots \nor \gr_n=\gr$$  on $\gr$ starting on $\gr$, ending on $\ce$, having length $n$ and without abelian jumps. Then, $$\card{\cord_{\gr_l}(\gr_k)}=2^{k-l} \quad\mbox{ for all }\; 0\leq l<k\leq n.$$ In particular, $$\card{\cord_\ce(\gr)}=2^n.$$
\begin{proof}
Fix $i,j\in\nat$ such that $0\leq l<k\leq n$. As a direct consequence of \Cref{biyeccion Conradiana general para todo n}, we obtain
\begin{equation*}
\card{\{\lleq\in \cord_{\gr_l}(\gr_k) \tq \gr_i \mbox{ are }\lleq-\mbox{convex for all } i=l+1,\dots,k-1\}}=\card{\prod_{j=l+1}^{k}\cord_{\gr_{j-1}}(\gr_j)}.
\end{equation*}
Then, we deduce that
\begin{equation*}
\card{\{\lleq\in \cord_{\gr_l}(\gr_k) \tq \gr_i \mbox{ are }\lleq-\mbox{convex for all } i=l+1,\dots,k-1\}}=\prod_{j=l+1}^{k}\card{\cord_{\gr_{j-1}}(\gr_j)}.
\end{equation*}
Notice that $\card{\cord_{\gr_{j-1}}(\gr_j)}=2$ follows by \Cref{ejemplo estructura en cocientes abelianos} as $\gr_{j+1}/\gr_j$ is torsion-free abelian of rank $1$ since the series considered above is rational. Hence, we obtain
\begin{equation}\label{numero de ordenes Conradianos desde serie eq 1}
\card{\{\lleq\in \cord_{\gr_l}(\gr_k) \tq \gr_i \mbox{ are }\lleq-\mbox{convex for all } i=l+1,\dots,k-1\}}=2^{k-l}.
\end{equation}
By \Cref{convexidad de series}, we know that, for all $\lleq\in\card{\cord_{\gr_l}(\gr_k)}$, we have $\gr_j$ is $\lleq$-convex for all $j=l+1,\dots,k-1$. Then, we obtain $$\card{\cord_{\gr_l}(\gr_k)}=\{\lleq\in \cord_{\gr_l}(\gr_k) \tq \gr_i \mbox{ are }\lleq-\mbox{convex for all } i=l+1,\dots,k-1\}.$$ Hence, by \Cref{numero de ordenes Conradianos desde serie eq 1}, we conclude that $$\card{\cord_{\gr_l}(\gr_k)}=2^{k-l}.$$ In particular, taking $k=0$ and $l=n$, we also deduce $$\card{\cord_\ce(\gr)}=2^n.$$
\end{proof}
\end{prop}

\begin{prop}\label{aislado implica serie de convexos finita}
Let $\gr$ be a group and $\ce$ be a proper subgroup. Let $\lleq$ be a Conradian left-preorder on $\gr$ relative to $\ce$. If $\lleq$ is isolated in $\cord_\ce(\gr)$ then there are finitely many $\lleq$-convex subgroups.
\begin{proof}
We prove this by contradiction. In this way, assume that $\lleq$ is isolated in $\cord_\ce(\gr)$ but there are infinitely many $\lleq$-convex subgroups. Since $\lleq$ is isolated in $\cord_\ce(\gr)$, there are $\ga_1,\dots,\ga_k\in \gr$ for $k\in\natp$ such that
\begin{equation}\label{aislado implica serie de convexos finita eq1}
\{\lleq_*\in \cord_\ce(\gr)\tq \ce \lle_* \ga_i\ce \;\mbox{ for }  0<i\leq k\}=\{\lleq\}.
\end{equation}
We know that there are infinitely many $\lleq$-convex subgroup and we know that the set of $\lleq$-convex subgroups is linearly ordered by \Cref{convexos estan linealmente ordenados}, so we can consider $$\sgr_1\subsetneq\sgr_2\subsetneq\dots\subsetneq\sgr_{k+1}\subsetneq\sgr_{k+2}$$ $k+2$ different $\lleq$-convex subgroups. Then, we obtain $k+1$ disjoint sets $\sgr_{j+1}\settminus\sgr_j$ for $j=1,\dots,k+1$. Then, since we have that $\ga_1,\dots,\ga_k$ are $k$ different elements, by cardinality reasons, we deduce that there is a $j_0\in \{1,\dots,k+1\}$ such that no element in $\{\ga_1,\dots,\ga_k\}$ belongs to $\sgr_{j_0+1}\settminus\sgr_{j_0}$. For simplicity on the notation, assume that $j_0=1$, so that $$\{\ga_1,\dots,\ga_k\}\cap(\sgr_2\settminus\sgr_1)=\emptyset.$$ We now consider the mutually inverse bijections $\mu$ and $\rho$ from \Cref{biyeccion Conradiana general para todo n} for the finite chain of subgroups $$\ce \subsetneq \sgr_1\subsetneq \sgr_2 \subsetneq \gr.$$ As we know that $\sgr_1$ and $\sgr_2$ are $\lleq$-convex subgroups, we can consider 
\begin{equation}\label{aislado implica serie de convexos finita eq2}
\rho(\lleq)=(\lleq_1,\lleq_2,\lleq_3) \;\mbox{ where }\; \lleq_1\in \cord_\ce(\sgr_1)\;\mbox{, }\; \lleq_2\in \cord_{\sgr_1}(\sgr_2) \;\mbox{ and }\; \lleq_3\in \cord_{\sgr_2}(\gr).
\end{equation}
Define $P$ the positive cone of $\lleq$, $P_1$ the positive cone of $\lleq_1$, $P_2$ the positive cone of $\lleq_2$ and $P_3$ the positive cone of $\lleq_3$. Then, by \Cref{aislado implica serie de convexos finita eq1}, we have $$\ga_1,\dots,\ga_k\in P$$ By the definition of the mutually inverse bijections $\mu$ and $\rho$ (see \Cref{biyeccion Conradiana general para todo n}), we deduce that 
\begin{equation*}
\ga_1,\dots,\ga_k\in P_1\cup P_2\cup P_3 \quad \mbox{and} \quad P_2\subseteq \sgr_2\settminus\sgr_1.
\end{equation*}
Since $\{\ga_1,\dots,\ga_k\}\cap(\sgr_2\settminus\sgr_1)=\emptyset$, we deduce that
\begin{equation}\label{aislado implica serie de convexos finita eq3}
\ga_1,\dots,\ga_k\in P_1\cup P_3.
\end{equation}
Consider $\lleq_2^-$ the element of $ \cord_{\sgr_1}(\gr_2)$ whose positive cone is $P_2\inv$, so that $\lleq_2^-$ is not equal to $\lleq_2$. Define 
\begin{equation}\label{aislado implica serie de convexos finita eq4}
\lleq_*=\mu(\lleq_1,\lleq_2^-,\lleq_3),
\end{equation}
so we have that $\lleq_*$ is not equal to $\lleq$, as $\mu$ is a bijection. We are going to reach a contradiction with the fact that $\lleq_*$ is not equal to $\lleq$. By \Cref{biyeccion Conradiana general para todo n}, the positive cone $P_*$ of $\lleq_*$ is given by $$P_*=P_1\cup P_2\inv\cup P_3.$$ From this and from \Cref{aislado implica serie de convexos finita eq3} we obtain $$\ga_1,\dots,\ga_k\in P_*,$$ which is equivalent to say that $$\ce \lle_* \ga_i\ce \;\mbox{ for }  0<i\leq k.$$ Applying \Cref{aislado implica serie de convexos finita eq1}, we deduce that $\lleq$ is equal to $\lleq_*$, which is a contradiction, as we wanted.
\end{proof}
\end{prop}

\begin{defi}
Let $\gr$ be a group and $\ce$ be a proper subgroup. Let $\lleq$ be a Conradian left-preorder on $\gr$ relative to $\ce$. By \Cref{convexos estan linealmente ordenados}, we know that the set of $\lleq$-convex subgroups is linearly ordered, so we can consider the subgroups series given by all $\lleq$-convex subgroups. We call it the \textit{series of} $\lleq$\textit{-convex subgroups}.
\end{defi}

\begin{prop}\label{aislado implica que sus convexos son serie racional finita sin saltos abelianos}
Let $\gr$ be a group and $\ce$ be a proper subgroup. Let $\lleq$ be a Conradian left-preorder on $\gr$ relative to $\ce$. If $\lleq$ is isolated in $\cord_\ce(\gr)$, then the series of $\lleq$-convex subgroups is a finite rational series on $\gr$ starting on $\gr$, ending on $\ce$ and without abelian jumps.
\begin{proof}
Assume that $\lleq$ is isolated in $\cord_\ce(\gr)$. By \Cref{aislado implica serie de convexos finita}, the set of $\lleq$-convex subgroups is finite, so the series of $\lleq$-convex subgroups is finite. Then, we can write the series of $\lleq$-convex subgroups as $$\ce=\gr_0 \subseteq \gr_1\subseteq \dots \subseteq \gr_n=\gr$$ where $\{\gr_0,\dots,\gr_n\}$ is the set of $\lleq$-convex subgroups. In particular, $(\gr_i,\gr_{i+1})$ is a $\lleq$-convex jump for all $0\leq i<n$ since there are not any other  $\lleq$-convex subgroups strictly in between. Using the fact that $\lleq$ is Conradian, by \Cref{Conradiano como cocientes de saltos}, we deduce that $\gr_i$ is normal on $\gr_{i+1}$ and the group $\gr_{i+1}/\gr_i$ is torsion-free abelian and non-trivial. We claim that $\gr_{i+1}/\gr_i$ has abelian rank equal to $1$. Indeed, if $\gr_{i+1}/\gr_i$ has abelian rank greater than $1$, then, by \Cref{ejemplo estructura en cocientes abelianos}, $\cord_{\gr_i}(\gr_{i+1})$ has no isolated orders, a contradiction with \Cref{lema aislado con convexos implica aislado en un cociente} and the fact that $\gr_i$ and $\gr_{i+1}$ are $\lleq$-convex and $\lleq$ is isolated in $\cord_\ce(\gr)$. Hence, we have shown that the series of convex subgroups given by $$\ce=\gr_0 \subseteq \gr_1\subseteq \dots \subseteq \gr_n=\gr$$ satisfies that $\gr_i$ is normal on $\gr_{i+1}$ and the group $\gr_{i+1}/\gr_i$ is torsion-free abelian of abelian rank one for all $0\leq i<n$. Hence, the series of convex subgroups given by $$\ce=\gr_0 \nor \gr_1\nor \cdots \nor \gr_n=\gr$$ is a finite rational series. We claim that this series has no abelian jump. We prove it by contradiction. Assume that $$\gr_i\nor\gr_{i+1}\nor\gr_{i+2}$$ is an abelian jump. Then, $\gr_i$ is normal on $\gr_{i+2}$ and the group $\gr_{i+2}/\gr_i$ is torsion-free abelian, and it must have abelian rank $2$ since we know that $\gr_{i+1}/\gr_i$ and $\gr_{i+2}/\gr_{i+1}$ both have abelian rank equal to $1$. Then, by \Cref{ejemplo estructura en cocientes abelianos}, $\cord_{\gr_i}(\gr_{i+2})$ has no isolated orders. This is a contradiction with \Cref{lema aislado con convexos implica aislado en un cociente} and the fact that $\gr_i$ and $\gr_{i+1}$ are $\lleq$-convex and $\lleq$ is isolated in $\cord_\ce(\gr)$. Hence, we have shown that $$\ce=\gr_0 \nor \gr_1\nor \cdots \nor \gr_n=\gr$$ is a finite rational series without abelian jumps.
\end{proof}
\end{prop}

We now have the partial version of the main theorem of the section.

\begin{teo}\label{teorema lista de equivalencias relativas}
Let $\gr$ be a group and $\ce$ be a proper subgroup. Then, the following conditions are equivalent:
\begin{itemize}
\item[a)] All the elements on $\cord_\ce(\gr)$ are isolated and $\cord_\ce(\gr)$ is non-empty.
\item[b)] There exists an element on $\cord_\ce(\gr)$ that is isolated.
\item[c)] There exists a finite rational series on $\gr$ starting on $\gr$, ending on $\ce$ and without abelian jumps.
\item[d)] $\cord_\ce(\gr)$ is finite and non-empty.
\end{itemize}
\begin{proof}
We are going to prove that each condition implies the next one. The fact that condition a) implies condition b) is direct. The fact that condition b) implies condition c) follows from \Cref{aislado implica que sus convexos son serie racional finita sin saltos abelianos}. The fact that  condition c) implies condition d) is a direct consequence of \Cref{numero de ordenes Conradianos desde serie}. We finally prove that condition d) implies condition a). By \Cref{espacio de ordenes es totalmente disconexo}, $\cord_\ce(\gr)$ is Hausdorff, so $\cord_\ce(\gr)$ being finite implies that $\cord_\ce(\gr)$ is discrete, so then all the elements on $\cord_\ce(\gr)$ are isolated.
\end{proof}
\end{teo}

\begin{lema}\label{lema acotando usando extensiones}
Let $\gr$ be a group and let $\ce$ be a proper subgroup. Let $\sgr$ be a subgroup of $\gr$ such that $\ce\subseteq \sgr$. If $\cord_\ce(\gr)\neq \emptyset$, then $$\card{\cord_{\sgr}(\gr)}\leq\card{\cord_\ce(\gr)}.$$
\begin{proof}
Assume that $\cord_\ce(\gr)\neq \emptyset$. Then, we can consider $\lleq$ a Conradian left-preorder on $\gr$ relative to $\ce$. Take $\lleq_0$ the restriction of $\lleq$ to $\sgr$. The map 
\begin{equation*}
\begin{array}{r@{\hspace{0pt}}c@{\hspace{0pt}}
c@{\hspace{4pt}}l}
& \cord_{\sgr}(\gr) &\longrightarrow& \cord_\ce(\gr) \\

& \lleq^*\ \ &\mapsto& \mu_{\gr,\sgr,\ce}(\lleq_0,\lleq^*)
\end{array}
\end{equation*}
from the set $\cord_{\sgr}(\gr)$ to the set $\cord_\ce(\gr)$ is injective by \Cref{ordenar desde dos ordenes}. We deduce that $\card{\cord_{\sgr}(\gr)}\leq\card{\cord_\ce(\gr)}$.
\end{proof}
\end{lema}

\begin{prop}\label{unicidad serie racional relativa}
Let $\gr$ be a group and $\ce$ be a proper subgroup. Let $\lleq$ be a Conradian left-preorder on $\gr$ relative to $\ce$. If $\cord_\ce(\gr)$ is finite, then the series of $\lleq$-convex subgroups is the unique finite rational series on $\gr$ starting on $\gr$ and ending on $\ce$.
\begin{proof}
Assume $\cord_\ce(\gr)$ is finite. By \Cref{teorema lista de equivalencias relativas}, we deduce that $\lleq$ is isolated in $\cord_\ce(\gr)$. Hence, by \Cref{aislado implica que sus convexos son serie racional finita sin saltos abelianos}, we obtain that the series of $\lleq$-convex subgroups is a finite rational series on $\gr$ starting on $\gr$ and ending on $\ce$. Therefore, it only remains to check the uniqueness. In order to do that, consider $$\ce=\gr_0 \subseteq \gr_1\subseteq \dots \subseteq \gr_n=\gr$$ the series of $\lleq$-convex subgroups. We are going to prove the uniqueness by induction in the length $n\in\natp$ of the previous series (the case $n=0$ is impossible as $\ce$ is proper). For $n=1$ it is trivial.\\

We now do the inductive case. We assume as inductive hypothesis that it is true for all length smaller than $n$, and we are going to prove it for $n>0$. Let $$\ce=\grb_0 \nor \grb_1\nor \cdots \nor \grb_m=\gr$$ be another finite rational series on $\gr$ starting on $\gr$ and ending on $\ce$. Denote $\grc=\gr_{n-1}$ and $\grb=\grb_{m-1}$. If we had that 
\begin{equation}\label{unicidad serie racional relativa eq1}
\grc=\grb,
\end{equation}
then applying inductive hypothesis to the series $$ \ce=\gr_0 \nor \gr_1\nor \cdots \nor \gr_{n-1}=\grc \quad \mbox{ and } \quad \ce=\grb_0 \nor \grb_1\nor \cdots \nor \grb_{m-1}=\grb=\grc$$
we could conclude that $n=m$ and the series
$$ \ce=\gr_0 \nor \gr_1\nor \cdots \nor \gr_n=\gr \quad \mbox{ and } \quad \ce=\grb_0 \nor \grb_1\nor \cdots \nor \grb_m=\gr$$
are equal, as we wanted. Therefore, we only need to prove \Cref{unicidad serie racional relativa eq1}.\\

Before proving it, we claim that $\gr/(\grc\cap\grb)$ is a torsion-free abelian group of abelian rank $1$. We now prove this. By definition of being rational, we have that $\gr/\grc$ and $\gr/\grb$ are torsion-free abelian groups. Then, the group  $\gr/\grc\times\gr/\grb$ is torsion-free abelian. Consider the group morphism
\begin{equation*}
\begin{array}{r@{\hspace{0pt}}c@{\hspace{0pt}}
c@{\hspace{4pt}}l}
& \gr &\longrightarrow& \gr/\grc\times\gr/\grb \\

& \ga\ \ &\mapsto& (\ga\grc,\ga\grb)
\end{array}.
\end{equation*}
It has kernel $\grc\cap\grb$, so the first isomorphism theorem give us an injective morphism from $\gr/(\grc\cap\grb)$ to $\gr/\grc\times\gr/\grb$, then $\gr/(\grc\cap\grb)$ is a torsion-free abelian group. By \Cref{lema acotando usando extensiones}, since $\lleq\in \cord_\ce(\gr)$ and since $\card{\cord_\ce(\gr)}$ is finite, we have that $\cord_{\grc\cap\grb}(\gr)$ is finite. Hence, we have that $\gr/(\grc\cap\grb)$ is a torsion-free abelian group such that $\cord_{\grc\cap\grb}(\gr)$ is finite, or equivalently, $\gr/(\grc\cap\grb)$ is a torsion-free abelian group having finitely many Conradian left-orders. This implies, by \Cref{ejemplo estructura en cocientes abelianos}, that $\gr/(\grc\cap\grb)$ is a torsion-free abelian group of abelian rank $1$ or $0$. Moreover, $\gr/(\grc\cap\grb)$ is a torsion-free abelian group of abelian rank $1$ since $\grc\cap\grb\neq\gr$. This proves the claim, as we wanted.\\

We are now able to prove \Cref{unicidad serie racional relativa eq1}. We prove this by showing that $\ga\notin \grc$ implies $\ga\notin \grb$ for each $\ga\in\gr$ and that $\ga\notin \grc$ implies $\ga\notin \grb$ for each $\ga\in\gr$. We just prove the first implication as the second one is analogous. Consider $\ga\in\gr$ such that $\ga\notin \grc$, we want to prove that $\ga\notin \grb$. As $\grb\neq\gr$, we can take $\gb\in\gr$ such that 
\begin{equation}\label{unicidad serie racional relativa eq2}
\gb\grb\neq \grb.
\end{equation}
Then we have that $\gb\notin \grb$ and $\ga\notin \grb$, so then $$\gb\notin \grc\cap\grb\quad\mbox{ and }\quad\ga\notin \grc\cap\grb.$$ Then, using the fact that $\gr/(\grc\cap\grb)$ is a torsion-free abelian group of abelian rank $1$, we deduce that there exists $r,s\in\ent$ non-zero numbers such that $$\ga^r\grc\cap\grb=\gb^s\grc\cap\grb.$$ Then, we have $$\ga^{-r}\gb^s\in\grc\cap\grb\subseteq \grb.$$ Then, we deduce $$\ga^r\grb=\gb^s\grb.$$ Applying \Cref{unicidad serie racional relativa eq2} and the fact that $r$ and $s$ are non-zero, we deduce $$\ga\grb\neq \grb.$$ Hence, $\ga\notin \grb$. This finishes the proof of \Cref{unicidad serie racional relativa eq1}.
\end{proof}
\end{prop}

We are able to prove the main theorem of the section.

\begin{teo}\label{teorema lista de equivalencias relativas con unicidad}
Let $\gr$ be a group and $\ce$ be a proper subgroup. Then, the following conditions are equivalent:
\begin{itemize}
\item[a)] All the elements on $\cord_\ce(\gr)$ are isolated and $\cord_\ce(\gr)$ is non-empty.
\item[b)] There exists an element on $\cord_\ce(\gr)$ that is isolated.
\item[c)] There exists a finite rational series on $\gr$ starting on $\gr$, ending on $\ce$ and without abelian jumps.
\item[d)] There exists a unique finite rational series on $\gr$ starting on $\gr$ and ending on $\ce$, and this series does not have abelian jumps.
\item[e)] $\cord_\ce(\gr)$ is finite and non-empty.
\end{itemize}
Moreover, under these equivalent conditions, we have that $\card{\cord_\ce(\gr)}=2^n$ for $n$ the length of any series as in c) and d).
\begin{proof}
The conditions a), b), c) and e) are equivalent by \Cref{teorema lista de equivalencias relativas}. The condition d) obviously implies the condition c). The condition c) implies the condition d) by \Cref{unicidad serie racional relativa} since $\cord_\ce(\gr)$ is finite by \Cref{teorema lista de equivalencias relativas}. By \Cref{numero de ordenes Conradianos desde serie}, we have that $\card{\cord_\ce(\gr)}=2^n$ for $n$ the length of any series as in c) and d)
\end{proof}
\end{teo}

\begin{obs}
Following the notation of the condition d) from \Cref{teorema lista de equivalencias relativas con unicidad}, we have that $\ce$ is not normal in general. However, when $\ce\nor \gr$, then the unique finite rational series on $\gr$ starting on $\gr$ and ending on $\ce$ satisfies that all the subgroups appearing on this series are normal on $\gr$. This is direct from the uniqueness and the fact that conjugating a finite rational series on $\gr$ starting on $\gr$ and ending on $\ce$ provides us another finite rational series on $\gr$ starting on $\gr$ and ending on $\ce$, since $\ce$ is normal on $\gr$. 
\end{obs}

\begin{coro}\label{numero de relativos es siempre potencia de dos}
Let $\gr$ be a group and $\ce$ be a proper subgroup. If $\cord_\ce(\gr)$ is finite and non-empty, then  $\card{\cord_\ce(\gr)}=2^n$ for some $n\in\nat$ such that $n>1$.
\begin{proof}
Apply \Cref{teorema lista de equivalencias relativas con unicidad}.
\end{proof}
\end{coro}

\begin{coro}\label{numero de relativos es finito o incontable}
Let $\gr$ be a group and $\ce$ be a proper subgroup. Then $\cord_\ce(\gr)$ is either finite or uncountable.
\begin{proof}
If $\cord_\ce(\gr)$ is not finite, by \Cref{teorema lista de equivalencias relativas con unicidad}, we deduce that $\cord_\ce(\gr)$  has no isolated points. Then, by \Cref{espacio de ordenes son cerrados y compactos 1} and \Cref{espacio de ordenes es totalmente disconexo}, we deduce that $\cord_\ce(\gr)$ is a non-empty compact Hausdorff space without isolated points. Then, we deduce that  $\cord_\ce(\gr)$ is uncountable (see Theorem 27.7 from \cite{topology-munkres} for example). Notice that, in fact, when $\gr$ is countable, then $\cord_\ce(\gr)$ is a Cantor set, as it is done in \Cref{observacion cantor}. 
\end{proof}
\end{coro}

\section{Groups having finitely many Conradian left-preorders}

In this section we want to prove \Cref{teorema lista de equivalencias absolutas con unicidad introduccion}, which is stated here as \Cref{teorema lista de equivalencias absolutas con unicidad}. The structure of this section is very similar to the one of the last section. We prove some technical result, and then we prove, \Cref{teorema lista de equivalencias absolutas}, which is a partial version of the main theorem. After that, we exhibit some auxiliary results, we prove the main theorem, and we conclude with some corollaries and some examples. 

\begin{lema}\label{lema interseccion finita de relativos en Conradianos}
Let $\gr$ be a group and $\ce_1,\dots,\ce_k$ be proper subgroups of $\gr$ for $k\in\natp$. If $\cord_{\ce_i}(\gr)\neq \emptyset$ for all $i=1,\dots,k$, then $\cord_{\ce_1\cap\dots\cap\ce_k}(\gr)\neq \emptyset$.
\begin{proof}
We prove the statement by induction on $k\in\natp$. For $k=1$ it is trivial. We now do the inductive case. Take $k\in\natp$ such that $k>1$ and assume that the statement is true for $k-1$. We want to prove it for $k$, so consider  $\ce_1,\dots,\ce_k$ proper subgroups of $\gr$ such that $\cord_{\ce_i}(\gr)\neq \emptyset$ for all $i=1,\dots,k$. By inductive hypothesis, we know that  $\cord_{\ce_1\cap\dots\cap\ce_{k-1}}(\gr)\neq \emptyset$. Then, we can take $\lleq_1$ a Conradian left-preorder on $\gr$ relative to $\ce_1\cap\dots\cap\ce_{k-1}$. Since we have $\cord_{\ce_k}(\gr)\neq \emptyset$, we can also take $\lleq_2$ a Conradian left-preorder on $\gr$ relative to $\ce_k$. Consider $\lleq_1'$ the restriction of $\lleq_1$ to $\ce_k$, so that $\lleq_1'$ is a Conradian left-preorder on $\ce_k$ relative to $\ce_1\cap\dots\cap\ce_k$. By \Cref{ordenar desde dos ordenes}, we have that $\mu_{\gr,\ce_1\cap\dots\cap\ce_k,\ce_k}(\lleq_1',\lleq_2)$ is a Conradian left-preorder on $\gr$ relative to $\ce_1\cap\dots\cap\ce_k$. Hence, $\cord_{\ce_1\cap\dots\cap\ce_k}(\gr)\neq \emptyset$. This finishes the proof for the inductive case.
\end{proof}
\end{lema}

\begin{prop}\label{interseccion arbitraria de relativos en Conradianos}
Let $\gr$ be a group and $(\ce_i)_{i\in I}$ be a family of proper subgroups of $\gr$. If $\cord_{\ce_i}(\gr)\neq \emptyset$ for all $i\in I$, then $$\cord_{\bigcap_{i\in I}\ce_i}(\gr)\neq \emptyset.$$
\begin{proof}
Consider the subsets $V_{\ga}^{r}$ and $\overline{V}_{\ga}^{r}$ of $\{\pm1,0\}^\gr$  we have defined in Section 3, which are open and closed at the same time. 
Define
\begin{equation*}
\begin{split}
A_0=&\bigcap_{\ga\in \gr}\left((V_{\ga}^{1}\cap V_{\ga\inv}^{-1})\cup (V_{\ga}^{-1}\cap V_{\ga\inv}^{1})\cup (V_{\ga}^{0}\cap V_{\ga\inv}^{0})\right),\\
A_1=&\{\pm 1,0\}^\gr\settminus\bigcup_{\ga,\gb\in \gr}\left(V_{\ga}^{1}\cap V_{\gb}^{1}\cap \overline{V}_{\ga\gb}^{1}\right),\\
A_2=&\bigcap_{\ga\in \bigcap_{i\in I}\ce_i}V_{\ga}^{0},\\
A_3=&\{\pm 1,0\}^\gr\settminus\bigcup_{\ga,\gb\in \gr}\left(V_{\ga}^{1}\cap V_{\gb}^{1}\cap \overline{V}_{\ga\inv \gb \ga^2}^{1}\right).
\end{split}
\end{equation*}
Notice that, by following the definitions, we have
\begin{equation}\label{interseccion arbitraria de relativos en Conradianos eq1}
\begin{split}
\chi\in A_0 \sii& \left(\chi\inv(\{1\})\right)\inv=\chi\inv(\{-1\}),\\
\chi\in A_1 \sii& \chi\inv(\{1\})\; \textrm{ is a subsemigroup of } \gr,\\
\chi\in A_2 \sii& \bigcap_{i\in I}\ce_i\subseteq \chi\inv(\{0\}),\\
\chi\in A_3 \sii& \mbox{ for each } \ga,\gb\in \chi\inv(\{1\})\; \mbox{ we have } \ga\inv\gb\ga^2\in \chi\inv(\{1\}).
\end{split}
\end{equation}
Define the subsets
\begin{equation}\label{interseccion arbitraria de relativos en Conradianos eq2}
\tilde{A}_{\{\gc_1,\dots,\gc_l\}}=\overline{V}_{\gc_1}^{0}\cap\dots \cap \overline{V}_{\gc}^{0} \quad \textrm{for each} \; \gc_1,\dots,\gc_l\in \gr\settminus (\bigcap_{i\in I}\ce_i) \; \textrm{ and }\; l\in\natp.
\end{equation}
By definition, for every $\gc_1,\dots,\gc_l\in \gr\settminus (\bigcap_{i\in I}\ce_i)$ we have
$$\chi\in \tilde{A}_{\{\gc_1,\dots,\gc_l\}} \sii \gc_r\notin\chi\inv(\{0\}) \; \textrm{for each } \; r=1,\dots,l\;.$$
We also define
\begin{equation*}
B_{\{\gc_1,\dots,\gc_l\}}=A_0\cap A_1\cap A_2\cap A_3\cap \tilde{A}_{\{\gc_1,\dots,\gc_l\}}\quad \textrm{for each} \; \gc_1,\dots,\gc_l\in \gr\settminus (\bigcap_{i\in I}\ce_i)\; \textrm{ and }\; l\in\natp.
\end{equation*}

We now show that
\begin{equation}\label{interseccion arbitraria de relativos en Conradianos eq3}
B_{\{\gc_1,\dots,\gc_l\}}\neq \emptyset  \;\mbox{ for all }\gc_1,\dots,\gc_l\in \gr \settminus (\bigcap_{i\in I}\ce_i) \mbox{ where }l\in\natp.
\end{equation}
Indeed, consider $\gc_1,\dots,\gc_l\in \gr\settminus (\bigcap_{i\in I}\ce_i)$ where $l\in\natp$. Then, for each $t=1,\dots, l$ we can take $i_t\in I$ such that $\gc_t\notin \ce_{i_t}$. In particular, we have $\gc_1,\dots,\gc_l\notin \ce_{i_1}\cap\dots\cap\ce_{i_l}$. By hypothesis, we have that $\cord_{\ce_{i_t}}(\gr)\neq \emptyset$ for all  $t=1,\dots, l$, so by \Cref{lema interseccion finita de relativos en Conradianos}, we deduce that $$\cord_{\ce_{i_1}\cap\dots\cap\ce_{i_l}}(\gr)\neq \emptyset.$$ Take $\lleq\in \cord_{\ce_{i_1}\cap\dots\cap\ce_{i_l}}(\gr)$. By the fact that $\gc_1,\dots,\gc_l\notin \ce_{i_1}\cap\dots\cap\ce_{i_l}$ and by applying the definitions, \Cref{interseccion arbitraria de relativos en Conradianos eq1}, \Cref{interseccion arbitraria de relativos en Conradianos eq2}, \Cref{orden relativo como semigrupo} and \Cref{Conradiano como cocientes de saltos}, we deduce that the element $\phi_\lleq$, corresponding to $\lleq$ seen as an element in  $\{\pm1,0\}^\gr$, belongs to $B_{\{\gc_1,\dots,\gc_l\}}$. This proves \Cref{interseccion arbitraria de relativos en Conradianos eq3}, as we wanted. \\

We consider the family $\Omega$ of non-empty closed subsets of $\{\pm 1, 0\}^\gr$ given by the sets $B_{\{\gc_1,\dots,\gc_l\}}$ for $\gc_1,\dots,\gc_l\in \gr\settminus (\bigcap_{i\in I}\ce_i)$ and $l\in\natp$. We claim that 
\begin{equation}\label{interseccion arbitraria de relativos en Conradianos eq4}
\bigcap_{B\in \Omega}B\neq \emptyset.
\end{equation}
Since all the elements in $\Omega$ are closed sets of the compact space $\{\pm 1, 0\}^\gr$, \Cref{interseccion arbitraria de relativos en Conradianos eq4} will be proven if we show that $\Omega$ has the finite intersection property. We now prove that $\Omega$ has the finite intersection property. Take $B_{\{\gc_{1,1},\dots,\gc_{1,l_1}\}}, \dots, B_{\{\gc_{t,1},\dots,\gc_{t,l_t}\}}\in \Omega$. By following the definitions, we have $$B_{\{\gc_{1,1},\dots,\gc_{1,l_1}, \dots, \dots,\gc_{t,1},\dots,\gc_{t,l_t}\}}\subseteq B_{\{\gc_{1,1},\dots,\gc_{1,l_1}\}}\cap \dots\cap B_{\{\gc_{t,1},\dots,\gc_{t,l_t}\}},$$ but $B_{\{\gc_{1,1},\dots,\gc_{1,l_1}, \dots \dots,\gc_{t,1},\dots,\gc_{t,l_t}\}}$ is non-empty by \Cref{interseccion arbitraria de relativos en Conradianos eq3}, so the intersection above is non-empty. Then, we have shown that $\Omega$ has the finite intersection property, so \Cref{interseccion arbitraria de relativos en Conradianos eq4} holds.\\

We are now going to prove that there exists a $\chi\in \{\pm 1, 0\}^\gr$ that will provide us a Conradian left-preorder with the needed properties. By \Cref{interseccion arbitraria de relativos en Conradianos eq3}, we can consider $$\chi\in\bigcap_{B\in \Omega}B.$$ By definition of $\Omega$, we have $\chi\in B_{\{\gc_1,\dots,\gc_l\}}$ for all $\gc_1,\dots,\gc_l\in \gr\settminus (\bigcap_{i\in I}\ce_i)$ and all $l\in\natp$. By \Cref{interseccion arbitraria de relativos en Conradianos eq1}, we have that $\chi\in A_i$ for $i=0,1,2,3$ and that $\chi\in \tilde{A}_{\{\gc_1,\dots,\gc_l\}}$ for all $\gc_1,\dots,\gc_l\in \gr\settminus (\bigcap_{i\in I}\ce_i)$ and all $l\in\natp$. We now prove that 
\begin{equation}\label{interseccion arbitraria de relativos en Conradianos eq5}
\bigcap_{i\in I}\ce_i= \chi\inv(\{0\}).
\end{equation}
Since $\chi\in A_2$, we deduce that $\bigcap_{i\in I}\ce_i\subseteq \chi\inv(\{0\})$. For the other inclusion, given $\gc\in \gr\settminus (\bigcap_{i\in I}\ce_i)$, we know that $\chi\in \tilde{A}_{\{\gc\}}$, so $\gc\notin\chi\inv(\{0\})$. Hence, we have shown \Cref{interseccion arbitraria de relativos en Conradianos eq5}.\\

We are now ready use $\chi$ to obtain the Conradian left-preorder. We can decompose $\gr$ as $$\gr=\chi\inv(\{1\})\sqcup\chi\inv(\{-1\})\sqcup\chi\inv(\{0\}).$$ Define $P=\chi\inv(\{1\})$, and notice that $P$ is a subsemigroup of $\gr$ since $\chi\in A_1$. Notice also that $P\inv=\chi\inv(\{-1\})$ as a consequence of $\chi\in A_0$. Then, we can write the the decomposition bellow as $$\gr=P\sqcup P\inv \sqcup \chi\inv(\{0\})$$ where $P$ is a subsemigroup of $\gr$. Therefore, using \Cref{interseccion arbitraria de relativos en Conradianos eq5}, the decomposition bellow is in fact $$\gr=P\sqcup P\inv \sqcup (\bigcap_{i\in I}\ce_i)$$ where $P$ is a subsemigroup of $\gr$. Hence, by \Cref{orden relativo como semigrupo}, there is a left-preorder $\lleq$ on $\gr$ relative to $\bigcap_{i\in I}\ce_i$ with positive cone $P$. Notice that $\lleq$ is Conradian because $\chi\in A_3$. We conclude that $$\lleq\in \cord_{\bigcap_{i\in I}\ce_i}(\gr).$$
\end{proof}
\end{prop}

\begin{lema}\label{interseccion series}
Let $\gr$ be a group. Let $$\gr_0 \nor \gr_1\nor \cdots \nor \gr_n=\gr \quad \mbox{ and }\quad\grb_0 \nor \grb_1\nor \cdots \nor \grb_k=\gr$$ be two finite rational series on $\gr$ both starting on $\gr$. Then, there exists $$\gr_0\cap \grb_0=\grc_0 \nor \grc_1\nor \cdots \nor \grc_l =\grb_0 \nor \grb_1\nor \cdots \nor \grb_k=\gr$$ a finite rational series on $\gr$ starting on $\gr$ and ending on $\gr_0\cap \grb_0$ of length $k+l$ for some $l\leq n$.
\begin{proof}
Notice that for every $i=1,\dots, n$ we have $$\gr_{i-1}\cap \grb_0\nor \gr_{i}\cap \grb_0$$ since $\gr_{i-1}\nor \gr_{i}$. First, we claim that for every $i=1,\dots, n$ we have 
\begin{equation}\label{interseccion series eq1}
\gr_{i-1}\cap \grb_0=\gr_{i}\cap \grb_0 \quad \mbox{ or } \quad (\gr_{i}\cap \grb_0)/(\gr_{i-1}\cap \grb_0) \mbox{ is torsion-free abelian of abelian rank } 1.
\end{equation}
We now prove it. Consider the morphism 
\begin{equation*}
\begin{array}{r@{\hspace{0pt}}c@{\hspace{0pt}}
c@{\hspace{4pt}}l}
\phi\colon & \gr_i\cap \grb_0 &\longrightarrow& \gr_i/\gr_{i-1} \\

& \ga\ \ &\mapsto& \ga \gr_{i-1}
\end{array}
\end{equation*}
and notice that its kernel is $$\{\ga\in \gr_i\cap \grb_0 \tq \ga \gr_{i-1}=\gr_{i-1}\}=\gr_{i-1}\cap\gr_{i}\cap \grb_0=\gr_{i-1}\cap \grb_0.$$ Then, by the isomorphism theorem, we can identify $(\gr_{i}\cap \grb_0)/(\gr_{i-1}\cap \grb_0)$ with a subgroup of $\gr_i/\gr_{i-1}$. Then, since $\gr_i/\gr_{i-1}$ is a torsion-free abelian group of abelian rank $1$ as the corresponding series is rational, we deduce that $(\gr_{i}\cap \grb_0)/(\gr_{i-1}\cap \grb_0)$ is a torsion-free abelian group of abelian rank $1$ or $0$. If $(\gr_{i}\cap \grb_0)/(\gr_{i-1}\cap \grb_0)$ has abelian rank $1$ we are done, and if it has abelian rank $0$, then the group is trivial, so we have $\gr_{i-1}\cap \grb_0=\gr_{i}\cap \grb_0$. Hence, we have proven \Cref{interseccion series eq1}.\\

Now, consider the finite series given by $$\grb_0\cap \gr_0 \nor \grb_0\cap \gr_1\nor \cdots \nor \grb_0\cap \gr_n =\grb_0 \nor \grb_1\nor \cdots \nor \grb_k=\gr.$$ When removing the repeated subgroups, using \Cref{interseccion series eq1} and the fact that $$\grb_0 \nor \grb_1\nor \cdots \nor \grb_k$$ is rational, we obtain $$\gr_0\cap \grb_0=\grc_0 \nor \grc_1\nor \cdots \nor \grc_l =\grb_0 \nor \grb_1\nor \cdots \nor \grb_k=\gr$$ a finite rational series on $\gr$ starting on $\gr$ and ending on $\gr_0\cap \grb_0$ of length $k+l$ for some $l\leq n$.
\end{proof}
\end{lema}

\begin{prop}\label{finito si y solo si el relativo al minimal es finito}
Let $\gr$ be a group. Assume that $\cord(\gr)\neq \emptyset$. Consider $\cce=\bigcap_{\sgr\in \hacheaux}\sgr $ where $\hacheaux=\{\sgr \mbox{ subgroup with } \cord_\sgr(\gr)\neq \emptyset \}$. If $\card{\cord_\cce(\gr)}<\infty$, then $\card{\cord(\gr)}<\infty$. Moreover, if $n\in\nat$ is such that $\card{\cord_\cce(\gr)}=2^n$, then $$\card{\cord(\gr)}=2^{n+1}-2.$$
\begin{proof}
Assume that $\card{\cord_\cce(\gr)}<\infty$. Take $n\in\nat$ such that $\card{\cord_\cce(\gr)}=2^n$, which is possible by \Cref{numero de relativos es siempre potencia de dos}. Then, by \Cref{teorema lista de equivalencias relativas}, we can take $$\cce=\gr_0 \nor \gr_1\nor \cdots \nor \gr_n=\gr$$ a finite rational series on $\gr$ starting on $\gr$, ending on $\cce$ and without abelian jumps.\\

We claim that for each $\sgr$ subgroup of $\gr$ we have
\begin{equation}\label{finito si y solo si el relativo al minimal es finito eq1}
\card{\cord_\sgr(\gr)}<\infty
\end{equation}
We now prove it. Consider $\sgr$ subgroup of $\gr$ such that $\cord_\sgr(\gr)\neq \emptyset$. By applying \Cref{interseccion arbitraria de relativos en Conradianos} and looking at the definition of $\cce$, we deduce that $\cord_\cce(\gr)\neq \emptyset$. Then, by \Cref{lema acotando usando extensiones} and since $\card{\cord_\cce(\gr)}<\infty$, we deduce that $\card{\cord_\sgr(\gr)}<\infty$. This proves \Cref{finito si y solo si el relativo al minimal es finito eq1}, as we claimed.\\

We claim that for each $\sgr$ subgroup of $\gr$ we have
\begin{equation}\label{finito si y solo si el relativo al minimal es finito eq2}
\cord_\sgr(\gr)\neq \emptyset \quad \Rightarrow \quad \sgr=\gr_i \;\mbox{ for some } i\in \{0,\dots,n-1\}.
\end{equation}
We now prove it. Consider $\sgr$ subgroup of $\gr$ such that $\cord_\sgr(\gr)\neq \emptyset$. By \Cref{finito si y solo si el relativo al minimal es finito eq1}, we know that $\card{\cord_\sgr(\gr)}<\infty.$. Then, applying \Cref{teorema lista de equivalencias relativas}, we deduce that we can take $$\sgr=\sgr_0 \nor \sgr_1\nor \cdots \nor \sgr_k=\gr$$ a finite rational series on $\gr$ starting on $\gr$, ending on $\sgr$ and without abelian jumps. Then, applying \Cref{interseccion series} to $$\cce=\gr_0 \nor \gr_1\nor \cdots \nor \gr_n=\gr \quad \mbox{ and } \quad \sgr=\sgr_0 \nor \sgr_1\nor \cdots \nor \sgr_k=\gr,$$ we deduce that there exists $$\cce\cap \sgr=\sgrc_0 \nor \sgrc_1\nor \cdots\nor \sgr =\sgr_0 \nor \sgr_1\nor \cdots \nor \sgr_k=\gr$$ a finite rational series on $\gr$ starting on $\gr$ and ending on $\cce\cap \sgr$. Notice that $\cce\cap \sgr=\cce$ by definition of $\cce$ and since $\cord_\sgr(\gr)\neq \emptyset$. Hence, there exists $$\cce=\sgrc_0 \nor \sgrc_1\nor \cdots\nor \sgr =\sgr_0 \nor \sgr_1\nor \cdots \nor \sgr_k=\gr$$ a finite rational series on $\gr$ starting on $\gr$ and ending on $\cce$, but we know that $\card{\cord_\cce(\gr)}<\infty$ and that $$\cce=\gr_0 \nor \gr_1\nor \cdots \nor \gr_n=\gr$$ is another  finite rational series on $\gr$ starting on $\gr$ and ending on $\cce$, so then the uniqueness given in \Cref{unicidad serie racional relativa}  implies that $ \sgr=\gr_i$ for some $i\in \{0,\dots,n-1\}$. This proves \Cref{finito si y solo si el relativo al minimal es finito eq2}.\\
By definition of $\cord(\gr)$, we have $$\cord(\gr)=\bigcup_{\sgrc \mbox{ subgroup of }\gr}\cord_{\sgrc}(\gr).$$ Applying \Cref{finito si y solo si el relativo al minimal es finito eq2}, we obtain $$\cord(\gr)=\bigcup_{i=0}^{n-1}\cord_{\gr_i}(\gr).$$ Since this union is clearly dijoint, we deduce that $$\card{\cord(\gr)}=\sum_{i=0}^{n-1}\card{\cord_{\gr_i}(\gr)}.$$ Then, applying \Cref{numero de ordenes Conradianos desde serie}, we obtain $$\card{\cord(\gr)}=\sum_{i=0}^{n-1}2^{n-i}=\sum_{j=1}^{n}2^j=\left(\sum_{j=0}^{n}2^j\right)-1=\left(2^{n+1}-1\right)-1= 2^{n+1}-2.$$ Therefore, $\card{\cord(\gr)}=2^{n+1}-2$ and, in particular, $\card{\cord(\gr)}<\infty$.
\end{proof}
\end{prop}

\begin{defi}
Let $\gr$ be a group. A rational series starting on $\gr$ of the form $$\gr_0 \nor \gr_1\nor \cdots \nor \gr_n=\gr$$ is of \textit{maximal length} if for any other $$\grb_0 \nor \grb_1\nor \cdots \nor \grb_k=\gr$$ rational series starting on $\gr$, we have $k\leq n$.
\end{defi}

\begin{lema}\label{lema tecnico lista de equivalencias absolutas}
Let $\gr$ be a group. Assume that $\cord(\gr)$ is non-empty. Let $\cce=\bigcap_{\sgr\in \hacheaux}\sgr $ for $\hacheaux=\{\sgr \mbox{ subgroup with } \cord_\sgr(\gr)\neq \emptyset \}$. Suppose that there exists a finite rational series on $\gr$ starting on $\gr$ of maximal length, and this series ends on $\gr_0$ satisfying $\cord(\gr_0)=\emptyset$ and has no abelian jump. Then, $\gr_0=\cce$.
\begin{proof}
The inclusion $\cce\subseteq \gr_0$ is clear by definition of $\cce$ and since $\cord_{\gr_0}(\gr)\neq \emptyset$ as a consequence of \Cref{teorema lista de equivalencias relativas}. We prove the other inclusion $\gr_0\subseteq \cce$ by contradiction, so assume that $\cce\subsetneq \gr_0$. By \Cref{interseccion arbitraria de relativos en Conradianos} and by the definition of $\cce$, we deduce that $\cord_\cce(\gr)\neq \emptyset$. Hence, we can take $\lleq$ a Conradian left-preorder on $\gr_0$ relative to $\cce$. Since $\cce\subsetneq \gr_0$, the restriction of $\lleq$ on $\gr_0$ give us a Conradian left-preorder on $\gr$ relative to $\cce$. Hence, we have obtained Conradian a left-preorder on $\gr_0$ relative to $\cce$, so we have an element on $\cord(\gr_0)$, a contradiction with the fact that $\cord(\gr_0)=\emptyset$. This shows that $\gr_0\subseteq \cce$. Therefore, we conclude that $\gr_0=\cce$.
\end{proof}
\end{lema}

We have a partial version of the main theorem.

\begin{teo}\label{teorema lista de equivalencias absolutas}
Let $\gr$ be a group. Assume that $\cord(\gr)$ is non-empty. Let $\cce=\bigcap_{\sgr\in \hacheaux}\sgr $ for $\hacheaux=\{\sgr \mbox{ subgroup with } \cord_\sgr(\gr)\neq \emptyset \}$. Then, the following conditions are equivalent:
\begin{itemize}
\item[a)] All the elements on $\cord(\gr)$ are isolated.
\item[b)] There exists a finite rational series on $\gr$ starting on $\gr$ of maximal length, and this series has no abelian jumps and satisfies that $\cord(\gr_0)=\emptyset$ for $\gr_0$ the subgroup on which the series ends.
\item[c)] There exists a finite rational series on $\gr$ starting on $\gr$ and ending on $\cce$, and this series has no abelian jumps.
\item[d)] $\cord_\cce(\gr)$ is finite.
\item[e)] $\cord(\gr)$ is finite.
\end{itemize}
\begin{proof}
We are going to prove that a) $\Rightarrow$ d), b) $\Leftrightarrow$ c), c) $\Rightarrow$ d), d) $\Rightarrow$ e), e) $\Rightarrow$ a) and e) $\Rightarrow$ c).\\

We prove that a) implies d). Assume that all the elements on $\cord(\gr)$ are isolated. Since all the elements on $\cord(\gr)$ are isolated and $\cord_\cce(\gr)$ is a subspace of $\cord(\gr)$, we deduce that all the elements on $\cord_\cce(\gr)$ are isolated. Then, by \Cref{teorema lista de equivalencias relativas}, we deduce that $\cord_\cce(\gr)$ is finite. This proves that a) implies d).\\

We show that  b) implies c). Assume that there exists a finite rational series on $\gr$ starting on $\gr$ of maximal length, and this series has no abelian jumps and satisfies that $\cord(\gr_0)=\emptyset$ for $\gr_0$ the subgroup on which the series ends. By \Cref{lema tecnico lista de equivalencias absolutas},  we have $\gr_0=\cce$. Then, the previous series is a finite rational series on $\gr$ starting on $\gr$ and ending on $\cce$ without abelian jumps, as we need. This shows that that  b) implies c).\\

We prove that c) implies b). Suppose that there exists a finite rational series on $\gr$ starting on $\gr$ and ending on $\cce$, and this series has no abelian jumps. If we show that the previous series has maximal length and $\cord(\cce)=\emptyset$, we deduce there exists a finite rational series on $\gr$ starting on $\gr$ of maximal length, and this series has no abelian jumps and satisfies that $\cord(\gr_0)=\emptyset$ for $\gr_0$ the subgroup on which the series ends, as we need. We are going to prove each fact independently. Firstly, we prove that $\cord(\cce)=\emptyset$ by contradiction. If $\cord(\cce)$ is non-empty, then there exists a subgroup $\sgr\subsetneq\cce$ such that $\cord_\sgr(\cce)$ is non-empty. Then, as $\cord_\cce(\gr)$ is non-empty by \Cref{teorema lista de equivalencias relativas}, applying \Cref{lema acotando usando extensiones} we obtain that $\cord_\sgr(\gr)$ is non-empty. This contradicts the fact that $\sgr\subsetneq\cce$ and the definition of $\cce$. Secondly, we prove that the finite rational series on $\gr$ starting on $\gr$ and ending on $\cce$ has maximal length. We do this by contradiction. Denote by $n$ the length of the finite rational series on $\gr$ starting on $\gr$ and ending on $\cce$, an assume that it is not the maximal length. Then, there exists a finite rational series on $\gr$ starting on $\gr$ and ending on some subgroup $\sgrb$ and having length $k$ for $k>n$. Then, applying \Cref{interseccion series} to these two series, we obtain a finite rational series on $\gr$ starting on $\gr$ and ending on $\sgrb\cap \cce$ and having length $r$ for $r>n$. Notice that $\sgrb\cap \cce=\cce$ by definition of $\cce$ and since $\cord_\sgr(\gr)\neq \emptyset$ by \Cref{teorema lista de equivalencias relativas}. Hence, we have a finite rational series on $\gr$ starting on $\gr$ and ending on $\cce$ and having length $r$ for $r>n$, but we also have  a finite rational series on $\gr$ starting on $\gr$ and ending on $\cce$ having length $n$, so we have a contradiction with the uniqueness provided by \Cref{unicidad serie racional relativa} (we know that $\cord_\cce(\gr)$ is finite by \Cref{teorema lista de equivalencias relativas}). This finishes the proof of c) implies b).\\

Now, c) implies d) is direct from \Cref{teorema lista de equivalencias relativas}. The fact that d) implies e) is clear by \Cref{finito si y solo si el relativo al minimal es finito}. The fact that e) implies a) follows since $\cord(\gr)$ is Hausdorff by \Cref{espacio de ordenes es totalmente disconexo} and every finite Hausdorff space is discrete, so all its elements are isolated. The fact that e) implies c) follows by noticing that $\cord_\cce(\gr)$ is s subset of $\cord(\gr)$.
\end{proof}
\end{teo}

\begin{obs}
Following the notation of \Cref{teorema lista de equivalencias absolutas}, if we consider $$\cce=\bigcap_{\sgr\in \hacheaux}\sgr  \quad \mbox{ for }\quad \hacheaux=\{\sgr \mbox{ subgroup with } \cord_\sgr(\gr)\neq \emptyset \},$$ then $\cce$ is a normal subgroup of $\gr$. This is direct from the fact that $\cord_\sgr(\gr)\neq \emptyset$ implies that $ \cord_{\ga\sgr\ga\inv}(\gr)\neq \emptyset$ for all $\ga\in\gr$ by just conjugating the positive cone. Then, one can consider the group $\gr/\cce$. One can build a map from $\cord(\gr)$ to $\cord(\gr/\cce)$ induced by the projection and, using the definition of $\cce$, one can show that it is in fact an homeomorphism. We do not prove this and we left it as an exercise since this is not necessary for our work, but we find it illustrative. In this way, knowing this fact, it is easy to prove that $\cord(\gr)$ being finite is equivalent to saying that there is $\ce$ a normal subgroup   in $\gr$ satisfying that $\cord(\ce)$ is empty and $\cord_{\{1\}}(\gr/\ce)$ is finite.
\end{obs}

\begin{prop}\label{unicidad serie racional maximal}
Let $\gr$ be a group. If there exists a finite rational series on $\gr$ starting on $\gr$ of maximal length, then it is unique.
\begin{proof}
Take $$\gr_0 \nor \gr_1\nor \cdots \nor \gr_n=\gr \quad \mbox{ and } \quad \grb_0 \nor \grb_1\nor \cdots \nor \grb_n=\gr$$ two finite rational series on $\gr$ starting on $\gr$ of maximal length $n$. We want to prove that both series are equal. If $n=0$ it is trivial, so assume $n>0$. In particular, we have that $\gr_0$ and $\grb_0$ are proper subgroups of $\gr$. By applying \Cref{interseccion series}, we deduce that there exists $$\gr_0\cap \grb_0=\grc_0 \nor \grc_1\nor \cdots \nor \grc_l =\grb_0 \nor \grb_1\nor \cdots \nor \grb_n=\gr$$ a finite rational series on $\gr$ starting on $\gr$ and ending on $\gr_0\cap \grb_0$ of length $n+l$ for some $l\leq n$. Since $n$ is the maximal length of such series, we deduce that $l=0$  and, in particular, by looking at the last series we have $\gr_0\cap \grb_0=\grb_0$. Hence, we have shown that $\grb_0\subseteq \gr_0$. Analogously and interchanging the two series of the beginning of the proof, we can also show that $\gr_0\subseteq \grb_0$. Therefore, we conclude that $\gr_0=\grb_0$. Hence, the series $$\gr_0 \nor \gr_1\nor \cdots \nor \gr_n=\gr \quad \mbox{ and } \quad \gr_0=\grb_0 \nor \grb_1\nor \cdots \nor \grb_n=\gr$$ are two finite rational series on $\gr$ starting on $\gr$ and ending on $\gr_0=\grb_0$. Also notice that $\cord_{\gr_0}(\gr)$ is finite by \Cref{teorema lista de equivalencias relativas}. Then, by the uniqueness from \Cref{unicidad serie racional relativa}, we deduce that these two series are equal.
\end{proof}
\end{prop}

We now present the main theorem.

\begin{teo}\label{teorema lista de equivalencias absolutas con unicidad}
Let $\gr$ be a group. Assume that $\cord(\gr)$ is non-empty. Let $\cce=\bigcap_{\sgr\in \hacheaux}\sgr $ for $\hacheaux=\{\sgr \mbox{ subgroup with } \cord_\sgr(\gr)\neq \emptyset \}$. Then, the following conditions are equivalent:
\begin{itemize}
\item[a)] All the elements on $\cord(\gr)$ are isolated.
\item[b)] There exists a finite rational series on $\gr$ starting on $\gr$ of maximal length, and this series has no abelian jumps and satisfies that $\cord(\gr_0)=\emptyset$ for $\gr_0$ the subgroup on which the series ends.
\item[c)] There exists a unique finite rational series on $\gr$ starting on $\gr$ of maximal length, and this series has no abelian jumps and satisfies that $\cord(\gr_0)=\emptyset$ for $\gr_0$ the subgroup on which the series ends.
\item[d)] There exists a finite rational series on $\gr$ starting on $\gr$ and ending on $\cce$, and this series has no abelian jumps.
\item[e)] There exists a unique finite rational series on $\gr$ starting on $\gr$ and ending on $\cce$, and this series has no abelian jumps.
\item[f)] $\cord_\cce(\gr)$ is finite.
\item[g)] $\cord(\gr)$ is finite.
\end{itemize}
Moreover, under these equivalent conditions, we have that $\card{\cord(\gr)}=2^{n+1}-2$ and $\card{\cord_\cce(\gr)}=2^n$ for $n$ the length of any series as in b), c), d) and e).
\begin{proof}
The fact that a), b), d), f) and g) are equivalent between them is by \Cref{teorema lista de equivalencias absolutas}. The fact that b) and c) are equivalent follows from \Cref{unicidad serie racional maximal}, and the fact that d) and e) are equivalent follows from \Cref{unicidad serie racional relativa} since $\cord_\cce(\gr)$ is finite as a consequence of \Cref{teorema lista de equivalencias relativas}. Finally, we prove that under these equivalent conditions, we have that $\card{\cord(\gr)}=2^{n+1}-2$ and $\card{\cord_\cce(\gr)}=2^n$ for $n$ the length of any series as in b), c), d) and e). Assume that we are under all the equivalent conditions of the statement. Since b) holds, \Cref{lema tecnico lista de equivalencias absolutas} implies that the rational series appearing on conditions b) and c) ends on $\cce$. Then, by the uniqueness given by \Cref{unicidad serie racional relativa}, we deduce that there is a unique rational series satisfying some of the equivalent conditions b), c), d) and e). Hence, they all have the same length, say $n$. By \Cref{numero de ordenes Conradianos desde serie}, we have $\card{\cord_\cce(\gr)}=2^n$. Then, by \Cref{finito si y solo si el relativo al minimal es finito}, we also deduce that $\card{\cord(\gr)}=2^{n+1}-2$. This shows that  $\card{\cord(\gr)}=2^{n+1}-2$ and $\card{\cord_\cce(\gr)}=2^n$ for $n$ the length of any series as in b), c), d) and e).
\end{proof}
\end{teo}

\begin{obs}
Following the notation of condition c) from \Cref{teorema lista de equivalencias absolutas con unicidad}, the unique finite rational series on $\gr$ starting on $\gr$ of maximal length satisfies that all the subgroups appearing on the series are normal on $\gr$. This is a consequence of the uniqueness and the fact that conjugating a finite rational series on $\gr$ starting on $\gr$ of maximal length gives us another finite rational series on $\gr$ starting on $\gr$ of maximal length. 
\end{obs}

\begin{ej}
If one tries to establish an analogy between \Cref{teorema lista de equivalencias absolutas con unicidad} and \Cref{teorema lista de equivalencias relativas con unicidad}, one may think that it is reasonable to say that in the list of equivalent conditions of \Cref{teorema lista de equivalencias absolutas con unicidad} the following condition is missing: there is an isolated element on $\cord(\gr)$. However, it is false that it is equivalent to the other conditions. Clearly, the condition a) implies that there is an isolated element on $\cord(\gr)$, but the reciprocal is not true. We show this with an example of a group having some isolated element and some non-isolated element. The example is the group $\ent^2$. Since it is abelian, all the subgroups are normal and all the left-preorders are Conradian.  Consider $\ce=\gen{(0,1)}$. Since $\ent^2/\ce$ is cyclic, there are only $2$ Conradian left-preorders on $\ent^2$ relative to $\ce$, say $\lleq_\ce^+$ and $\lleq_\ce^-$. Then, we have $$\{\lleq'\in \cord(\ent)\tq \phi_{\lleq}(0,1)=0,\; \phi_{\lleq}(1,0)=1 \}=\{\lleq_\ce^+\} \;\mbox{ and }$$ $$\{\lleq'\in \cord(\ent)\tq \phi_{\lleq}(0,1)=0,\; \phi_{\lleq}(1,0)=-1 \}=\{\lleq_\ce^-\}.$$ Hence, both Conradian left-preorders are isolated in $\cord(\ent^2)$. However, by \Cref{ejemplo estructura en cocientes abelianos} we know that $\cord_{\{(0,0)\}}(\ent^2)$ has no isolated elements, so all the elements in $\cord_{\{(0,0)\}}(\ent^2)$ are not isolated in $\cord(\ent^2)$.
\end{ej}

\begin{coro}\label{numero de absolutos es finito o incontable}
Let $\gr$ be a group. Then $\cord(\gr)$ is either finite or uncountable.
\begin{proof}
Assume that $\cord(\gr)$ is not finite. By \Cref{teorema lista de equivalencias absolutas con unicidad}, we deduce that $\cord_\ce(\gr)$ is not finite for a certain $\ce$ subgroup. Then, by \Cref{numero de relativos es finito o incontable}, we deduce that $\cord_\ce(\gr)$ is uncountable. Then, since $\cord_\ce(\gr)$ is a subset of $\cord(\gr)$, we deduce that $\cord(\gr)$ is uncountable.
\end{proof}
\end{coro}

We now present a proposition that allows us to skip the requirement of $\cord(\gr_0)=\emptyset$ on conditions b) and c) of \Cref{teorema lista de equivalencias absolutas con unicidad} when we know that $\gr_0$ is finitely generated (for example, for groups whose subgroups are all finitely generated). 

\begin{prop}\label{para finitamente generado me puedo saltar una condicion}
Let $\gr$ be a group. Assume that there exists a finite rational series on $\gr$ starting on $\gr$ of maximal length. Denote by $\gr_0$ the subgroup where this series ends. If $\gr_0$ is finitely generated, then $\cord(\gr_0)=\emptyset$.
\begin{proof}
Suppose that exists $$\gr_0 \nor \gr_1\nor \cdots \nor \gr_n=\gr$$ a finite rational series on $\gr$ starting on $\gr$ of maximal length $n$, so that this series ends on the finitely generated subgroup $\gr_0$. We do this by contradiction, so assume $\cord(\gr_0)\neq\emptyset$. Then, we can consider $\lleq$ a Conradian left-preorder on $\gr_0$ relative to $\ce$ for $\ce$ a proper subgroup of $\gr_0$. By \Cref{lema convexo maximal}, there exists a maximal proper $\lleq$-convex subgroup $\sgr$ of $\gr_0$. This implies that $(\sgr,\gr_0)$ is a $\lleq$-convex jump. Then, since $\lleq$ is Conradian, by \Cref{Conradiano como cocientes de saltos}, we deduce that $\sgr$ is normal on $\gr_0$ and the group $\gr_0/\sgr$ is torsion-free abelian and non-trivial. In fact $\gr_0/\sgr$ is finitely generated as $\gr_0$ is so by hypothesis, so $\sgr$ is normal on $\gr_0$ and the group $\gr_0/\sgr$ is finitely-generated torsion-free abelian and non-trivial. Then, we can take $\df{\phi}{\gr_0/\sgr}{\ent}$ a surjective group morphism, and we can also take the surjective group morphism $\df{\pi}{\gr_0}{\gr_0/\sgr}$ given by the canonical projection. Hence, the group morphism $\df{\phi\co\pi}{\gr_0}{\ent}$ is surjective. Then, if we consider $\sgrc$ the kernel of $\phi\co\pi$, we deduce that $\sgrc$ is normal on $\gr_0$ and the group $\gr_0/\sgrc$ is isomorphic to $\ent$. In particular, $\sgrc$ is normal on $\gr_0$ and $\gr_0/\sgrc$ is a torsion-free abelian group of abelian rank $1$. Hence, $$\sgrc\nor\gr_0 \nor \gr_1\nor \cdots \nor \gr_n=\gr$$ is a finite rational series on $\gr$ starting on $\gr$ of length $n+1$. This contradicts that th maximal length is $n$.
\end{proof}
\end{prop}

In the previous proposition it is necessary to put the condition of $\gr_0$ being finitely generated, as we can see in the following example inspired on the collection of groups $(T_n)_{n\in\natp}$ from \cite{Conradianfinito-rivas}.

\begin{ej}\label{ejemplo T}
Consider $T$ to be the group generated by the set of generators $\{a_i \tq i\in\natp\}$ and by the relations $$a_{i+1}a_ia_{i+1}\inv= a_i\inv \;\; \mbox{ for all } \; i\in\natp \quad \mbox{ and } \quad a_ia_j=a_ja_i \;\; \mbox{ for all } \; \abs{i-j}\geq 2.$$ For each $n\in\natp$ we consider the subgroup $$T_n=\gen{a_1,\dots,a_n},$$ and this $T_n$ corresponds to the one given in \cite{Conradianfinito-rivas}, and in this reference it is shown that $T_n$ admits a Conradian left-order (in fact, it has $2^n$ Conradian left-orders). Notice that we have $$\{1\}\nor T_1\nor T_2\nor \cdots\nor T_n\nor \cdots\; \subseteq T \quad \mbox{ and } \quad T=\bigcup_{n\in\natp} T_n,$$ where the series $\{1\}\nor T_1\nor T_2\nor \cdots\nor T_n\nor \cdots$ is rational. Clearly, each finitely generated subgroup of $T$ is contained in some $T_n$, which we know admits a Conradian left-order. Then, all finitely generated subgroup of $T$ admits a Conradian left-order so, by Theorem 9.15 from \cite{orderedgroupstopology-clay}, we deduce that $T$ admits a Conradian left-order. Hence, $\cord(T)\neq \emptyset$.\\

Notice that $T$ has trivial abelianization, so there is no normal proper subgroup on $T$ whose quotient is abelian. Hence, the unique finite rational series on $T$ starting on $T$ of maximal length is the one of length $0$, that is, the series with just $T$ and no other subgroups. As we have shown before $\cord(T)\neq \emptyset$. This example shows us that in \Cref{para finitamente generado me puedo saltar una condicion} the condition of $\gr_0$ being finitely generated cannot be avoided.
\end{ej}

\begin{coro}\label{teorema lista de equivalencias absolutas con unicidad para grupos Conradiano-ordenables}
Let $\gr$ be a group. Assume that $\gr$ admits a Conradian left-order. Then, the following conditions are equivalent:
\begin{itemize}
\item[a)] There is an isolated element on $\cord_{\{1\}}(\gr)$.
\item[b)] All the elements on $\cord_{\{1\}}(\gr)$ are isolated.
\item[c)] All the elements on $\cord(\gr)$ are isolated.
\item[d)] There exists a finite rational series on $\gr$ starting on $\gr$ and ending on $\{1\}$, and this series has no abelian jumps.
\item[e)] There exists a unique finite rational series on $\gr$ starting on $\gr$ and ending on $\{1\}$, and this series has no abelian jumps.
\item[f)] $\cord_{\{1\}}(\gr)$ is finite.
\item[g)] $\cord(\gr)$ is finite.
\end{itemize}
Moreover, under these equivalent conditions, we have that $\card{\cord(\gr)}=2^{n+1}-2$ and $\card{\cord_{\{1\}}(\gr)}=2^n$ for $n$ the length of any series as in d) and e).
\begin{proof}
The set of Conradian left-orders on $\gr$ is equal to $\cord_{\{1\}}(\gr)$. Then, the fact that $\gr$ admits a Conradian left-order implies that $\cord_{\{1\}}(\gr)$ is non-empty. This implies that $\cord(\gr)$ is non-empty and that $$\bigcap_{\sgr \mbox{ subgroup with } \cord_\sgr(\gr)\neq \emptyset}\sgr=\{1\}.$$ Knowing this, we deduce from \Cref{teorema lista de equivalencias absolutas con unicidad} that the conditions c), d), e), f) and g) are equivalent and that, under these conditions, we have that $\card{\cord(\gr)}=2^{n+1}-2$ and $\card{\cord_{\{1\}}(\gr)}=2^n$ for $n$ the length of any series as in d) and e). From \Cref{teorema lista de equivalencias relativas} we deduce that the conditions a), b) and d) are equivalent.
\end{proof}
\end{coro}

We use all this machinery to compute the number of left-preorders and the number of Conradian left-preorders in the collection of groups $(T_n)_{n\in\natp}$ from Section 4.2 of \cite{Conradianfinito-rivas} and from Example 2.2.14 of \cite{god}, that has also appeared in \Cref{ejemplo T}.

\begin{ej}\label{ejemplo Tn}
For each $n\in\natp$, define $T_n$ as the group generated by $\{a_1,\dots,a_n\}$ with the relations $$a_{i+1}a_ia_{i+1}\inv= a_i\inv \;\; \mbox{ for all } \; i=1,\dots,n-1\quad \mbox{ and } \quad a_ia_j=a_ja_i \;\; \mbox{ for all } \; \abs{i-j}\geq 2.$$ 
Fix $n\in\natp$. One can prove that, if $\lleq$ is a left-preorder on $T_n$ relative to $\ce$, then $\ce=\gen{a_1,\dots,a_k}$ for some $k=0,\dots,n-1$ (this follows from some easy computations using the relations and basic properties of left-preorders). Then, since $T_n$ is left-orderable, one can extend it to a left-order. As in Example 2.2.14 of \cite{god}, one can check that $\lleq$ is determined by the sign (which is $1$ if it is on the positive cone or $-1$ if its inverse is on the positive cone) of the elements $a_{k+1},\dots,a_n$, so there are at most $2^{n-k}$ left-preorders of this form. When considering all possible values of $k$, we deduce that there are at most $$\sum_{k=0}^{n-1}2^{n-k}= 2^{n+1}-2$$ left-preorders on $T_n$. As in Example 2.2.14, notice that $$\{1\}\nor \gen{a_1}\nor \gen{a_1,a_2}\nor\dots\nor \gen{a_1,\dots,a_n}=T_n$$ is a rational series, and it is easy to see that it has no abelian jumps. By \Cref{teorema lista de equivalencias absolutas con unicidad para grupos Conradiano-ordenables}, $T_n$ has $2^{n+1}-2$ Conradian left-preorders. Since we know that $T_n$ has at most $2^{n+1}-2$ left-preorders, we conclude that the space of left preorders on $T_n$ and the space of Conradian left preorders on $T_n$ are equal and consists on exactly $2^{n+1}-2$ left-preorders.
\end{ej}

The previous example can be used to determine the set of left-preorders on the Klein-bottle group.

\begin{ej}\label{ejemplo botella klein}
Consider $K$ the Klein-bottle group, which is generated by $\{a,b\}$ with the relation $a b a\inv= b\inv$. As it is done in \Cref{ejemplo Tn}, we deduce that it has $6$ left-preorders. In fact, from \Cref{ejemplo Tn} one can deduce that the possible left-preorders are given by the following six different positive cones: $\semgen{a,b}$, $\semgen{a\inv,b}$, $\semgen{a,b\inv}$, $\semgen{a\inv,b\inv }$, $\semgen{a}$ and $\semgen{a\inv}$.
\end{ej}

In the last two examples we have that the space of left-preorders and  the space of Conradian left-preorders are both finite. However, this is not always the case, as we see in the following example.

\begin{ej}\label{ejemplo Cn}
In Section 4.2 of \cite{Conradianfinito-rivas}, Rivas considers a collection of groups $(C_n)_{n\in\natp}$ such that $C_n$ has $2^n$ Conradian left-orders but infinitely many left-orders. Since all left-orders are left-preorders, we deduce that $C_n$ has infinitely many left-preorders for all $n\in\natp$. Applying \Cref{teorema lista de equivalencias absolutas con unicidad para grupos Conradiano-ordenables}, we deduce that this collection $(C_n)_{n\in\natp}$ satisfies that $C_n$ has $2^{n+1}-2$ Conradian left-preorders. This, for each  $n\in\natp$, provides us an example of a group having $2^{n+1}-2$ Conradian left-preorders but infinitely many left-preorders.
\end{ej}

\section{Conradian left-preorders from a dynamical point of view}

We now provide a dynamical version of Conradian left-preorders whose statement and proof follow the ideas on Theorem 1.4 from \cite{caracterizacionconrad-navasrivasclay}. Notice that these actions also appear on \cite{morris-antolinrivas}.

\begin{defi}
Let $\gr$ be a group acting on a totally ordered set $(\espacio, \men)$ by order-preserving bijections. A \textit{crossing} for this action consists on $(\ga, \gb; \uu, \vv, \ww)$, where $\ga,\gb\in \gr$ and $\uu,\vv,\ww\in\espacio$, and satisfying the following conditions:
\begin{enumerate}
\item $\uu\men \ww \men \vv$.
\item $\gb^n(\uu)\men \vv$ and $\uu\men \ga^n(\vv)$ for all $n\in \nat$.
\item There exists $N,M\in \nat$ such that $\ga^N(\vv)\men \ww\men \gb^M(\uu)$.
\end{enumerate}
\end{defi}

\begin{teo}\label{Conradiano como dinamica}
Let $\gr$ be a group and $\ce$ be a proper subgroup. Let $\lleq$ be a left-preorder relative to $\ce$. Then, $\lleq$ is Conradian if and only if the action of $\gr$ on $(\gr/\ce,\lleq)$ by left multiplication admits no crossing.
\begin{proof}
We start by proving that $(\gr/\ce,\lleq)$ admits no crossing implies that $\lleq$ is Conradian. We do it by contraposition. Assume that $\lleq$ is not Conradian, so there exist $\ga,\gb\in\gr$ such that $\ce\lle\ga\ce$ and $\ce\lle\gb\ce$ and satisfying
\begin{equation}\label{crossing ecuacion hipotesis}
\ga\gb^n\ce\lle\gb\ce\quad \mbox{ for all } \quad n\in\nat.
\end{equation}
Define $\uu=\ce$, $\vv=\ga\inv \gb\ce$ and $\ww=\gb^2 \ce$. We prove that $(\ga, \gb; \uu, \vv, \ww)$ is a crossing for the action by checking each condition from the definition of a crossing. \\

We start by checking the first condition. Using \Cref{crossing ecuacion hipotesis}, we have $\ga\gb^2 \ce\lle \gb\ce$. We deduce $\gb^2 \ce\lle\ga\inv \gb\ce$. Then, using the fact that $\ce\lle\gb^2\ce$ since $\ce\lle\gb\ce$, we obtain $\ce\lle \gb^2\ce\lle \ga\inv \gb\ce$. Hence, we have proven $\uu\men \ww \men \vv$, so first condition of the definition of crossing holds. \\

We now check the second condition. Fix an arbitrary $n\in\nat$. On the one hand, \Cref{crossing ecuacion hipotesis} implies $\gb^n\ce\lle\ga\inv\gb\ce$, which is equivalent to $\gb^n(\uu)\men \vv$. On the other hand, \Cref{crossing ecuacion hipotesis} implies $\ga\ce\lle\gb\ce$, so then $\ce\lle\ga\inv\gb\ce$. Then, since $\ce\lle\ga\ce$ and since the product of $\lleq$-positive elements is $\lleq$-positive, we deduce that $\ce\lle\ga^n\ga\inv\gb\ce$, which is the same as saying $\uu\men \ga^n(\vv)$. Hence, the second condition of the definition holds.\\

We check the third condition. Firstly, since $\ce\lle\gb\ce$, we have $\gb\ce\lle\gb^2\ce$. Then, since $\ga\ga\inv\gb\ce=\gb\ce$, we have $\ga\ga\inv\gb\ce\lle\gb^2\ce$. This is the same as saying $\ga^N(\vv)\men \ww$ for $N=1$. Secondly, we know $\ce\lle\gb\ce$, so then $\gb^2\ce\lle\gb^3\ce$. This is equivalent to $\ww\men \gb^M(\uu)$ for $M=3$. Therefore, the third condition of the definition holds. We conclude that  $(\ga, \gb; \uu, \vv, \ww)$ is a crossing for the corresponding action, as we wanted.\\

Conversely, we want to prove that if $\lleq$ is Conradian then $(\gr/\ce,\lleq)$ admits no crossing, and we are going to do it by contraposition. Assume that $(\ga, \gb; \uu, \vv, \ww)$ is a crossing for $(\gr/\ce,\lleq)$, so we are going to prove that $\lleq$ is not Conradian. Take $\uu=\urep\ce$, $\vv=\vrep\ce$ and $\ww=\wrep\ce$. By definition of crossing, we have that
\begin{equation}\label{crossing primera ecuacion}
\urep\ce\lle \wrep\ce \lle \vrep\ce
\end{equation}
\begin{equation}\label{crossing segunda ecuacion}
\gb^n\urep\ce\lle \vrep\ce \mbox{ and }\urep\ce\lle \ga^n\vrep\ce \mbox{ for all }n\in \nat
\end{equation}
\begin{equation}\label{crossing tercera ecuacion}
\mbox{There exists }N,M\in \nat \mbox{ such that }\ga^N\vrep\ce\lle \wrep\ce\lle \gb^M\urep\ce
\end{equation}
holds. We are going to show that $\lleq$ is not Conradian by showing that the elements $\wrep\inv \gb^M\ga^N\wrep$ and $\wrep\inv \gb^M\wrep$ satisfy that $\ce\lle\wrep\inv \gb^M\ga^N\wrep\ce$, that $\ce\lle\wrep\inv \gb^M\wrep\ce$ and that $(\wrep\inv \gb^M\ga^N\wrep)(\wrep\inv \gb^M\wrep)^n\ce\lle \wrep\inv \gb^M\wrep\ce$ for all $n\in\nat$.\\

Firstly, we are going to prove that $\ce\lle\wrep\inv \gb^M\ga^N\wrep\ce$. On the one hand, from \Cref{crossing tercera ecuacion} we obtain $\ga^N\vrep\ce\lle \wrep\ce$, so then $\gb^M\ga^{2N}\vrep\ce\lle \gb^M\ga^N\wrep\ce$. On the other hand, from \Cref{crossing segunda ecuacion} we have $\urep\ce\lle \ga^{2N}\vrep\ce$, so then $\gb^M\urep\ce\lle \gb^M\ga^{2N}\vrep\ce$. Combining these two conclusions, we have $$\gb^M\urep\ce\lle \gb^M\ga^{2N}\vrep\ce\lle \gb^M\ga^N\wrep\ce,$$ but $\wrep\ce\lle \gb^M\urep\ce$ by \Cref{crossing tercera ecuacion}, so then $\wrep\ce\lle \gb^M\ga^N\wrep\ce$, and hence $\ce\lle \wrep\inv\gb^M\ga^N\wrep\ce$, as we wanted.

Secondly, we are going to show that $\ce\lle\wrep\inv \gb^M\wrep\ce$. We notice that $\wrep\ce\lle\gb\wrep\ce$ because, if not, we would have $\gb\wrep\ce\lleq\wrep\ce$, and it would imply $$\gb^M\wrep\ce\lleq\gb^{M-1}\wrep\ce\lleq\cdots \lleq\gb\wrep\ce\lleq\wrep\ce,$$ but using $\gb^M\urep\ce\lleq\gb^M\wrep\ce$ which is direct from \Cref{crossing primera ecuacion}, we deduce $\gb^M\urep\ce\lleq\wrep\ce$, but by \Cref{crossing tercera ecuacion} we have $\wrep\ce\lle\gb^M\urep\ce$, hence $\wrep\ce\lle\wrep\ce$, which is a contradiction. Hence, we have $\wrep\ce\lle\gb\wrep\ce$, so then $$\wrep\ce\lle\gb\wrep\ce\lle\cdots \lle\gb^{M-1}\wrep\ce\lle\gb^M\wrep\ce,$$ and by left multiplying by $\wrep\inv$, we obtain $\ce\lle\wrep\inv \gb^M\wrep\ce$, as we needed.\\

Thirdly, it only remains to prove $(\wrep\inv \gb^M\ga^N\wrep)(\wrep\inv \gb^M\wrep)^n\ce\lle \wrep\inv \gb^M\wrep\ce$ for all $n\in\nat$. Consider an arbitrary $n\in\nat$. From \Cref{crossing tercera ecuacion} we have $\wrep\ce\lle \gb^M\urep\ce$, so then $\gb^{Mn}\wrep\ce\lle \gb^{Mn+M}\urep\ce$, and using $\gb^{Mn+M}\urep\ce\lle\vrep\ce$ which is direct from \Cref{crossing segunda ecuacion}, we deduce $\gb^{Mn}\wrep\ce\lle\vrep\ce$, and hence $(\gb^{M})^n\wrep\ce\lle\vrep\ce$. Then, we obtain $$(\gb^M\ga^N)(\gb^{M})^n\wrep\ce\lle\gb^M\ga^N\vrep\ce,$$ but we have $\gb^M\ga^N\vrep\ce\lle\gb^M\wrep\ce$ from \Cref{crossing tercera ecuacion} and left multiplying by $\gb^M$, so then $$(\gb^M\ga^N)(\gb^{M})^n\wrep\ce\lle\gb^M\wrep\ce,$$ and hence $$\wrep\inv(\gb^M\ga^N)(\gb^{M})^n\wrep\ce\lle\wrep\inv\gb^M\wrep\ce.$$ Notice that $$\wrep\inv(\gb^M\ga^N)(\gb^{M})^n\wrep\ce=(\wrep\inv \gb^M\ga^N\wrep)(\wrep\inv \gb^M\wrep)^n\ce\lle \wrep\inv \gb^M\wrep\ce,$$ so we conclude $$(\wrep\inv \gb^M\ga^N\wrep)(\wrep\inv \gb^M\wrep)^n\ce\lle \wrep\inv \gb^M\wrep\ce$$ as we wanted.
\end{proof}
\end{teo}

\noindent {\bf Acknowledgements}. The author thanks Yago Antolín for the supervision and for all the crucial suggestions and comments about the paper. The author thanks Cristóbal Rivas for reading the first draft and providing comments and suggestions. The author thanks Enric Ventura the careful reading of the PhD thesis that has motivated this publication, whose feedback has been crucial to improve the paper and to correct some errors.

\end{document}